\documentclass[10pt,a4paper]{article}

\usepackage{amssymb}
\usepackage{amsmath}
\usepackage{algorithm}     
\usepackage{algpseudocode} 
\usepackage{subfig}
\usepackage{csquotes}
\usepackage{tikz}
\usetikzlibrary{arrows}
\usepackage{mathtools}
\usepackage{bbm}
\usepackage{amsthm}
\usepackage{pgfplots}
\usepackage{float}
\usepackage{hyperref}
\usepackage{dutchcal}
\usepgfplotslibrary{external}

\newcommand{\vertiii}[1]{{\left\vert\kern-0.25ex\left\vert\kern-0.25ex\left\vert #1 
   \right\vert\kern-0.25ex\right\vert\kern-0.25ex\right\vert}}

\newtheorem{theorem2}{Theorem}[section]

\newtheorem{remark}[theorem2]{Remark}

\usepackage{geometry}
\makeatletter
\def\blfootnote{\xdef\@thefnmark{$\star$}\@footnotetext}
\makeatother
\newenvironment{Authors}%
  {\begin{center}\begin{bfseries}}%
  {\end{bfseries}\end{center}}
\newenvironment{Addresses}%
  {\begin{flushleft}\begin{itshape}}%
  {\end{itshape}\end{flushleft}}
  \newcommand{\email}[1]{\hspace*{\stretch{1}}\emph{\texttt{#1}}}

 \usepackage{fancyhdr}  
  \fancypagestyle{plain}{
\fancyhead{}
\fancyhead[C]{\hfill Submitted to Elsevier, October 2023}
 }
\begin{document}

\thispagestyle{plain}

\title{Model order reduction  by  convex displacement interpolation}
 \date{}
 \maketitle

 \maketitle
\vspace{-50pt} 
 
\begin{Authors}
Simona Cucchiara$^{1}$, Angelo Iollo$^{2}$, Tommaso Taddei$^{2}$, Haysam Telib$^{1}$
\end{Authors}

\begin{Addresses}
$^1$
Optimad Engineering, 10134 Torino, Italy \\  \email{simona.cucchiara@optimad.it,haysam.telib@optimad.it} \\
$^2$
Univ. Bordeaux, CNRS, Bordeaux INP, IMB, UMR 5251, F-33400 Talence, France\\ Inria Bordeaux Sud-Ouest, Team MEMPHIS, 33400 Talence, France,
\email{angelo.iollo@inria.fr,tommaso.taddei@inria.fr} \\
\end{Addresses}

\begin{abstract}
We present a   nonlinear interpolation technique for parametric fields that exploits optimal transportation of coherent structures of the solution to achieve accurate performance.
The approach generalizes the nonlinear interpolation procedure introduced in  [Iollo, Taddei, J. Comput. Phys., 2022] 
to multi-dimensional parameter domains 
and to datasets of several snapshots.
Given a library of high-fidelity simulations, we rely on a scalar testing function  and on a point set registration method to identify coherent structures of the solution field in the form of sorted point clouds.
Given a new parameter value, we exploit a regression method to predict the new point cloud;
then, we resort to a boundary-aware registration technique to define bijective mappings that deform the new point cloud into the point clouds of the neighboring elements of the dataset, while preserving the boundary of the domain;
finally, we define the estimate as a weighted combination of modes obtained by composing the   neighboring snapshots with the previously-built mappings.
We present several numerical examples
for compressible and incompressible, viscous and inviscid flows to demonstrate the accuracy of the method.
Furthermore, we employ
 the nonlinear interpolation procedure to augment the dataset of simulations for linear-subspace projection-based model reduction: our data augmentation procedure is designed to reduce offline costs --- which are dominated by snapshot generation --- of model reduction techniques for nonlinear advection-dominated problems.
\end{abstract}

\medskip

 \section{Introduction}
\label{sec:introduction}

\subsection{Nonlinear interpolation of parametric fields}

Despite the many recent contributions, model order reduction (MOR) of parametric problems with compactly-supported features  --- 
such as shocks or 
shear layers
---
remains an outstanding task  for state-of-the-art techniques due to the fundamental inadequacy of linear approximations.
 The aim of this work is to devise a general --- i.e., independent of the underlying equations ---
 interpolation technique  for steady-state parametric problems, with emphasis on fluid dynamics applications.

During the past decade, several authors have proposed 
\emph{Lagrangian  methods} to deal with this class of problems
\cite{ching2022model,iollo2014advection,mirhoseini2023model,mojgani2017arbitrary,mojgani2021low,ohlberger2013nonlinear,taddei2015reduced,taddei2020registration}.
We denote by $\mu$ the vector of model parameters in the region $\mathcal{P}\subset \mathbb{R}^P$ and  we denote by $\Omega\subset \mathbb{R}^d$ the open computational domain; then, 
we introduce the parametric field of interest $u:\Omega\times \mathcal{P} \to \mathbb{R}^D$
and the solution manifold
$\mathcal{M}=\{  u_{\mu}:=u(\cdot; \mu): \mu\in \mathcal{P} \}$.
Lagrangian approximations rely on the ansatz
\begin{subequations}
\label{eq:lagrangian}
\begin{equation}
\label{eq:lagrangian_a}
\widehat{u}_{\mu} = 
\widetilde{u}_{\mu}
\circ \Phi_{\mu}^{-1},
\end{equation}
 where $\widetilde{u}_{\mu}$ is a linear (or affine) approximation of the form
 \begin{equation}
 \label{eq:lagrangian_b}
 \widetilde{u}_{\mu}(x) =  \sum_{i=1}^n \widehat{\omega}_{\mu}^i  \zeta_i(x) ,
 \;\;
 x\in \Omega, \mu\in \mathcal{P},
 \end{equation}
for proper choices of the weights $\widehat{\omega}_{\mu}^1,\ldots, \widehat{\omega}_{\mu}^n$ and the parameter-independent fields $\zeta_1,\ldots,\zeta_n: \Omega \to \mathbb{R}^D$, and
$\Phi: \Omega\times \mathcal{P} \to \Omega$ is a suitably-chosen bijection that tracks the coherent structures of the solution; here, $D$ denotes the number of state variables, while $d$ is the spatial dimension.
\end{subequations}  
Lagrangian approaches are motivated by the observation that in many problems of interest coherent structures that are  troublesome for linear approximations vary smoothly with the parameter and they hence can be tracked through a low-rank parameter-dependent mapping $\Phi$.

In \cite{iollo2022mapping}, we (Iollo, Taddei) proposed a general method dubbed  convex displacement interpolation (CDI) 
that relies on optimal transportation to perform accurate nonlinear interpolations between solution snapshots; the approach was proposed  for databases of two snapshots and one-dimensional parameter domains.
Similarly to Lagrangian approaches, CDI relies on the assumption that the location of coherent features of the solution field depends smoothly on the parameter; however, unlike in \eqref{eq:lagrangian}, it does not rely on the definition of a reference configuration where the location of the coherent features is (approximately) freezed.
In this work, we discuss an extension of CDI to multi-dimensional parameter domains and to datasets of several snapshots.

\subsection{Convex displacement interpolation}

\begin{figure}[H]
\input{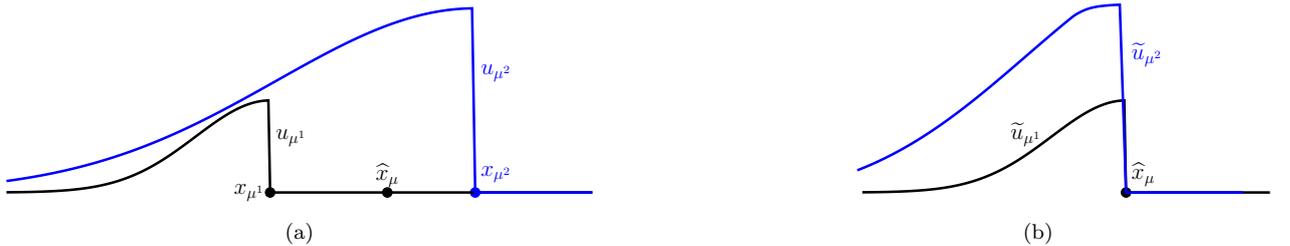}
 \caption{schematic of the CDI procedure.
 (a) dataset of snapshots $\{u_{\mu^1}, u_{\mu^2}\}$ and corresponding sensors (shock location),
 $\widehat{x}_{\mu}$ denotes the predicted shock location for the parameter $\mu$.
 (b) mapped snapshots $\{\widetilde{u}_{\mu^1}, \widetilde{u}_{\mu^2}\}$ in \eqref{eq:CDI}.
}
 \label{fig:nozzle_vis}
 \end{figure}

CDI exploits  the standard offline/online decomposition to predict the state map $\mu\in \mathcal{P} \mapsto u_{\mu}\in \mathcal{M}$;
Algorithm \ref{alg:CDI_overview} summarizes the general procedure and highlights the key steps of the methodology that will be discussed in the next sections.
Given the solution manifold $\mathcal{M}$, we introduce the training set $\mathcal{P}_{\rm train}=\{ \mu^k \}_{k=1}^{n_{\rm train}} \subset \mathcal{P}$ and the dataset of solutions 
$\mathcal{D}_{\rm train}=\{ u_{\mu}: \mu\in  \mathcal{P}_{\rm train}\}$. During the offline (or learning) stage, for each element of 
$\mathcal{D}_{\rm train}$, we determine a set of points
$X_{\mu}^{\rm raw}=\{  x_{i,\mu}^{\rm raw} \}_{i=1}^{N_{\mu}}\subset \overline{\Omega}$ that describe local features of the solution field that we wish to track: note that the sets
$\{ X_{\mu}^{\rm raw}: \mu\in \mathcal{P}_{\rm train}\}$ do not necessarily have the same number of elements and are not matched with each other.
Then, we define the reference points $X^{\rm ref}:= \{ x_i^{\rm ref}  \}_{i=1}^N$ and the new sensors 
$X_{\mu}:= \{ x_{i,\mu} = x_i^{\rm ref} + v_{i,\mu} \}_{i=1}^N$
 such that the set $X_{\mu}$ approximates --- in a sense to be defined --- $X_{\mu}^{\rm raw}$ for all $\mu\in \mathcal{P}_{\rm train}$.
 During the online stage, given the new parameter value $\mu\notin \mathcal{P}_{\rm train}$, first, we  apply a regression method to predict the new sensor locations $\widehat{X}_{\mu} =\{ \widehat{x}_{i,\mu} = x_i^{\rm ref} + \widehat{v}_{i,\mu} \}_{i=1}^N$;
second, we identify a set of $\kappa$ nearest neighbors $\mathcal{P}_{\rm nn}^{\mu}=\{  \nu^i \}_{i=1}^{\kappa}\subset \mathcal{P}_{\rm train}$; 
third, we define the mappings $\{ \Phi_{\nu} : \nu \in \mathcal{P}_{\rm nn}^{\mu}  \}$ such that each mapping $\Phi_{\nu}$ is bijective in $\Omega$ and $\Phi_{\nu}(  \widehat{x}_{i,\mu}  ) \approx \widehat{x}_{i,\nu}$ for $i=1,\ldots,N$ and $\nu\in \mathcal{P}_{\rm nn}^{\mu} $;
fourth, we return the (generalized) CDI as
\begin{equation}
\label{eq:CDI}
\widehat{u}_{\mu} = \sum_{\nu  \in \mathcal{P}_{\rm nn}^{\mu}}
\omega_{\mu}^{\nu} 
\widetilde{u}_{\nu},
\quad
{\rm where} \;\;
\widetilde{u}_{\nu} = 
u_{\nu}\circ \Phi_{\nu},
\end{equation}
for a proper choice of the weights 
$\{ \omega_{\mu}^{\nu}:   \nu \in \mathcal{P}_{\rm nn}^{\mu} \}$.
Figure \ref{fig:nozzle_vis} illustrates the computational procedure for a simple one-dimensional problem.

\begin{algorithm}[H]                      
\caption{: offline/online decomposition.}     
\label{alg:CDI_overview}     
 \normalsize 

\textbf{Offline stage} \hfill \emph{performed once} 
\begin{algorithmic}[1]
\State
Generate the dataset 
$\mathcal{D}_{\rm train}=\{ u_{\mu}: \mu\in  \mathcal{P}_{\rm train}\}$.
\medskip

\State
Identify the point clouds $\{X_{\mu}^{\rm raw} : \mu\in   \mathcal{P}_{\rm train} \}$.
\hfill section \ref{sec:feature_selection}
\medskip

\State
Define the template set $X^{\rm ref}$ and the sorted point clouds
$\{X_{\mu}: \mu\in   \mathcal{P}_{\rm train} \}$.

\hfill section \ref{sec:PSR}

\end{algorithmic}

\medskip

\textbf{Online stage} 
\hfill
\emph{performed for any} $\mu\in \mathcal{P}$
\begin{algorithmic}[1]

\State
Estimate the new points 
$\widehat{X}_{\mu}=\{ \widehat{x}_{i,\mu}  \}_{i=1}^N$.
\hfill section \ref{sec:regression}
\medskip

\State
Select the neighboring parameters
 $\mathcal{P}_{\rm nn}^{\mu}=\{ \nu^i \}_{i=1}^{\kappa}
\subset  \mathcal{P}_{\rm train}$.
\hfill section \ref{sec:nearest_neighbors}
\medskip

\State
Compute the mappings 
 $\Phi_{\nu}$ based on $\widehat{X}_{\mu}$ and 
 $\widehat{X}_{\nu}$ for all
 $\nu\in \mathcal{P}_{\rm nn}^{\mu}$.
 
 \hfill section \ref{sec:registration}
\medskip

\State
Compute the weights $\{ \omega_{\mu}^{\nu} : \nu \in \mathcal{P}_{\rm nn}^{\mu} \}$ and return the estimate \eqref{eq:CDI}.

 \hfill section \ref{sec:weights}
\end{algorithmic}
\end{algorithm}

The practical implementation of  CDI requires to carefully choose appropriate methods for each step of the offline and online stages  of Algorithm \ref{alg:CDI_overview}. 

\paragraph{Sensor selection}
There is a broad literature in fluid mechanics and scientific computing on the definition of quantitative criteria for the identification of coherent flow features such as shear layers, shocks, or vortices: representative examples include shock sensors for high-order methods in CFD
\cite{nicoud1999subgrid,persson2006sub}, and
local criteria for vortex identification
\cite{jeong1995identification}.

\paragraph{Point cloud matching}
The problem of determining sorted point clouds based on a dataset  of unmatched point clouds and a template set is referred to as 
 point-set registration 
(PSR, \cite{horaud2010rigid,ma2015robust,ma2018nonrigid,maiseli2017recent,zitova2003image}):
PSR methods implicitly introduce a global smooth mapping that coherently deforms the points of the template set $X^{\rm ref}$ to match --- in a sense to be defined --- the target cloud $X_{\mu}$, for all $\mu\in \mathcal{P}_{\rm train}$.

\paragraph{Regression}
The problem of estimating the  sensor locations for new parameter values
is a standard multi-target regression problem that can be handled through off-the-shelves machine learning routines including radial basis function (RBF, \cite{wendland2004scattered}) regression.
As in  non-intrusive MOR methods (e.g., \cite{guo2018reduced}), we first apply 
proper orthogonal decomposition (POD, 
\cite{sirovich1987turbulence,volkwein2011model}) to determine a low-rank representation of the sorted point clouds and then we apply RBF to the POD coefficients.

\paragraph{Registration}
The mapping generated by PSR methods does not fulfill the bijectivity condition $\Phi(\Omega)=\Omega$;  in order to avoid extrapolation of full-order solutions outside the domain of definition, it is hence necessary to devise an efficient boundary-aware registration technique that ensures bijectivity in $\Omega$. Note that the registration step should be significantly less expensive than  full-order solves in order for the CDI method to be beneficial.

\paragraph{Choice of the weights}
The weights 
$\{ \omega_{\mu}^{\nu} : \nu \in \mathcal{P}_{\rm nn}^{\mu} \}$
should be chosen to ensure relevant properties of the interpolation method (cf. section \ref{sec:properties}). 
Since our emphasis is on datasets of very modest size, we rely on 
the inverse distance weighting (IDW, \cite{shepard1968two}) method.

\subsection{Contributions and relation to previous works}

The past decade has witnessed a surge in the development of nonlinear MOR methods
to mitigate the Kolmogorov barrier 
\cite{amsallem2008interpolation,barnett2022quadratic,barnett2022neural,lee2020model,peherstorfer2020model}.
As discussed in the previous section, our method shares important features with Lagrangian methods based on the ansatz \eqref{eq:lagrangian}; the ansatz \eqref{eq:CDI} is also formally equivalent to  front tracking methods proposed in \cite{krah2023front,reiss2018shifted}.
Note also that our approach relies on a piecewise-nonlinear approximation, which is motivated by the need to reduce the mesh interpolation and registration costs that scale linearly with $\kappa$ and also by the particular strategy employed to define the weights $\{ \omega_{\mu}^{\nu}: \nu\in \mathcal{P}_{\rm nn}^{\mu}\}$. In this respect, the method is related to previously-proposed localized MOR techniques \cite{amsallem2012nonlinear,eftang2010hp}.

The CDI procedure outlined in  Algorithm \ref{alg:CDI_overview}
generalizes the method in \cite{iollo2022mapping}: 
in more detail, it reduces to the two-field CDI estimate for one-dimensional parameters and $\kappa=2$, 
 and a proper choice of the mappings 
 $\{ \Phi_{\nu} : \nu\in \mathcal{P}_{\rm nn}^{\mu} \}$
(cf.   Remark \ref{remark:equivalenceCDI}).
We observe that, similarly to Lagrangian approximations, CDI relies on geometric mappings to track local features of the solution field and ultimately improve state predictions;
on the other hand, CDI does not rely on the definition of a reference domain and the estimate is inherently Eulerian. We note that the definition  of the template $X^{\rm ref}$ is instrumental to the construction of the sorted point clouds $\{ X_{\mu}:\mu\in \mathcal{P}_{\rm train} \}$, which enable the application of standard regression algorithms during the online stage for the prediction of $\widehat{X}_{\mu}$.
We also observe that during the online stage each mapping  $ \Phi_{\nu} $ is computed based on the point clouds $\widehat{X}_{\mu}$ and $\widehat{X}_{\nu}$ (in lieu of $X_{\nu}$): 
thanks to this choice, we prove that CDI is an interpolation method (cf. section \ref{sec:properties}).

In this work, we rely on well-established techniques for sensor selection, point cloud matching, multi-target regression, and weight selection;
on the other hand, we emphasize that the development of boundary-aware registration techniques to determine bijective mappings in the domain $\Omega$ of interest is an active research area with many recent contributions.
In this work, we propose  a new elasticity-based method that borrows ideas from penalization methods for PDE discretization on Cartesian meshes
(see, e.g., 
 \cite{chantalat2009level,angot1999penalization});
we also note that the method shares important features with morphing techniques
for mesh adaptation \cite{tezduyar1992new,tonon2021linear} and image processing
\cite{camion2001geodesic,cao2005large}.
We further consider the optimization-based method first proposed in
 \cite{taddei2020registration} and further developed in several subsequent works.

The ultimate goal of this work is to devise a nonlinear interpolation strategy that can cope with datasets of very modest size.
The evaluation of the CDI estimate \eqref{eq:CDI} is not 
necessarily \emph{online-efficient} --- that is, 
the prediction cost scales with the cost of the underlying high-fidelity (HF) mesh --- but it is significantly less expensive  than the solution to the HF model.
In this respect, CDI might be employed  to initialize nonlinear solvers for nonlinear problems
\cite{kadeethum2022enhancing}, or to enrich the snapshot dataset that can later be used to train a projection-based reduced-order model (PROM) --- in effect, \emph{data augmentation} \cite{bernard2018reduced}.
We discuss the application of CDI
to data augmentation through the vehicle of a parametric two-dimensional   inviscid flow 
problem.

The outline of the paper is as follows.
In section \ref{sec:methods}, we present the methods for each step of the offline and online stages of Algorithm \ref{alg:CDI_overview}, and we show relevant properties of the interpolation procedure.
In section \ref{sec:motivating_example}, 
we  present a motivating example that illustrates the many features of the approach in a simplified setting.
In section \ref{sec:nonlinear_interpolation_numerics} we present extensive numerical results for three different test cases that cover compressible and incompressible, viscous and inviscid flows.
In section \ref{sec:data_augmentation_numerics},
we illustrate the application to data augmentation.
Section \ref{sec:conclusions} completes the paper.

\section{Methodology}
\label{sec:methods}

\subsection{Sensor (or feature) selection}
\label{sec:feature_selection}
We rely on the definition of a problem-dependent \emph{scalar testing function} $\mathcal{T}$  to identify unmatched point clouds that can be exploited to derive  out-of-sample  predictions (cf. section \ref{sec:regression}).
Given the domain $\Omega$, we define a discrete set of points $P_{\rm hf} \subset \overline{\Omega}$ and we compute $\{\mathcal{T}(x, u_{\mu}^{\rm hf}): x\in P_{\rm hf} \}$ for all $\mu\in \mathcal{P}_{
\rm train}$, where $u_{\mu}^{\rm hf}$ denotes the high-fidelity (HF) estimate of the solution field for $\mu\in \mathcal{P}$; then, we set 
\begin{equation}
\label{eq:raw_point_clouds}
X_{\mu}^{\rm raw} = \left\{
x_{i,\mu}^{\rm raw} 
\right\}
=
\left\{
x\in P_{\rm hf} : 
\mathcal{T}(x, u_{\mu}^{\rm hf}) \geq 
\mathcal{t}_{\mu}
\right\}
\end{equation}
for some threshold $\mathcal{t}_{\mu}$ that is specified below.

The choice of the scalar testing function is highly problem-dependent and requires a deep understanding of the underlying physical model and of the features we wish to track.  For compressible flow problems in the transonic or supersonic regime, we rely on the Ducros sensor \cite{modesti2017low,nicoud1999subgrid} to identify the shocks:
we denote by 
$P_{\rm hf} = \{ x_k^{\rm hf}  \}_{k=1}^{N_{\rm e}}$ the  cell  centers of the HF mesh and we compute
\begin{subequations}
\label{eq:ducros_sensor}
\begin{equation}
\mathcal{T}(x_k^{\rm hf}; u_{\mu}^{\rm hf})
=\sup_{x\in \texttt{D}_k} 
\big| \phi(x;  u_{\mu}^{\rm hf})  \big|
\end{equation}
where $\texttt{D}_k$ denotes the $k$-th cell of the mesh with  center  $x_k^{\rm hf}$, $\phi$ is the Ducros sensor
\begin{equation}
\phi(x;  u)
:=
\dfrac{(- \nabla \cdot v)^{+}}{\sqrt{(\nabla \cdot v)^2 + \Vert \nabla \times v \Vert_2^2 + a^2}}\dfrac{\Vert \nabla p \Vert_2 }{p + \epsilon}\Vert v \Vert_2,
\end{equation}
$v$ is the velocity field, $p$ is the pressure,
$a$ is the sound speed and $\epsilon>0$ is a user-defined tolerance that is set equal to $0.01$.
Finally, the threshold $\mathcal{t}_{\mu}$ in \eqref{eq:raw_point_clouds} is chosen equal to the $\gamma_{\rm thr}$  quantile  over the training set:
\begin{equation}
\label{eq:quantile}
\mathcal{t}_{\mu} =\texttt{quantile} \left(
\{\mathcal{T}(x, u_{\mu}^{\rm hf}): x\in P_{\rm hf} \}, 
\gamma_{\rm thr}
\right).
\end{equation}
\end{subequations}
For the incompressible flow problem of section \ref{sec:nonlinear_interpolation_numerics}, we rely on the isolines of the streamfunction to identify the recirculation area, for two-dimensional channel flows.
In more detail, we define the streamline function 
\begin{equation}
\label{eq:streamline_function}
\Psi_{\mu}\left(x = [x_1,x_2] \right)
\, = \,
\int_{y_{\rm btm}(x_1)}^{x_2} \, \left( u_{\mu}(x_1,s) \right)_1 ds,
\end{equation}
where 
$y_{\rm btm}(x_1)=\inf \{x_2 \, : \,  [x_1,x_2]\in \Omega   \}$ denotes the bottom boundary of the channel; then, we set 
$\mathcal{T}(x, u_{\mu}^{\rm hf}) = -
\Psi_{\mu}\left(x  \right) $ and $\mathfrak{t}_{\mu} \equiv 0$; that is,
$$
X_{\mu}^{\rm raw} = \left\{
x\in P_{\rm hf} \, : \,
\Psi_{\mu}\left(x  \right)  \leq 0
\right\}.
$$
We remark that we also tried to exploit the Q-criterion and the triple decomposition of relative motion (TDM) 
\cite{kolavr2007vortex}: however, in our experience, these two indicators were not able to distinguish between   pure shearing motions and the actual swirling motion of the vortex.

\subsection{Point cloud matching}
\label{sec:PSR}

The point clouds in \eqref{eq:raw_point_clouds} are not necessarily of the same size and are not sorted; therefore, we resort to PSR to find matched point clouds:
\begin{equation}
\label{eq:matched_point_clouds}
X_{\mu} = \texttt{PSR} \left( X^{\rm ref}, X_{\mu}^{\rm raw}  \right),
\quad
X^{\rm ref}:= 
 X_{\mu^{\rm ref}}^{\rm raw}, 
 \;\;   \mu^{\rm ref} \in \mathcal{P}_{\rm train}.
\end{equation}
Given the point clouds $X=\{ x_i \}_{i=1}^N$ and
$Y=\{ y_j \}_{j=1}^M$, the problem of PSR consists in finding a map $T:\mathbb{R}^d \to \mathbb{R}^d$ that (approximately) minimizes the distance 
$$
{\rm dist} \left( Y, T(X) \right)
=
\max_{y\in Y} 
\left(
\min_{x\in X} \| y - T(x)  \|_2
\right)
$$
in a suitable model class; the output of the algorithm is the deformed set of points $T(X)$. A broad range of PSR methods
--- including the approach considered in this work ---
  relies on a probabilistic interpretation of the problem: these methods rely on the assumption that $\{ x_i \}_{i=1}^N$ and
$\{ y_j \}_{j=1}^M$ are independent identically distributed (iid) samples from a given probability distribution with unknown parameters.
As discussed in the introduction, Algorithm \ref{alg:CDI_overview}  can cope with any PSR  technique  that can be written in the form \eqref{eq:matched_point_clouds}; we refer to a future work for the development of a specialized PSR procedure for MOR.
In the remainder of this section, we review  the PSR algorithm employed in the numerical experiments of sections \ref{sec:nonlinear_interpolation_numerics} and \ref{sec:data_augmentation_numerics}.

We rely on Gaussian-based PSR.  Given the point cloud $X=\{ x_i \}_{i=1}^N$, we resort to maximum likelihood estimation
(MLE, see, e.g., 
\cite[Chapter 8]{rice2006mathematical}) to estimate mean and covariance matrix:
\begin{equation}
\label{eq:gauss_model}
\mu_{X} = \frac{1}{N} \;
\sum_{i=1}^{N}  x_i,
\qquad
\Sigma_{X} = \frac{1}{N} \;
\sum_{i=1}^{N}   \,  ( x_i - \mu_{X}) \, ( x_i - \mu_{X})^{\top},
\end{equation}
where $(\cdot)^{\top}$ is the transpose operator, and we define the Gaussian model
$g_X= \mathcal{N}( \mu_{X},  \Sigma_{X} )$.
Similarly, given the point cloud $Y=\{ y_j \}_{j=1}^M$, we define the MLE estimates
$\mu_{Y},\Sigma_{Y}$ and the corresponding Gaussian model
$g_Y= \mathcal{N}( \mu_{Y},   \Sigma_{Y} )$.
Finally, we apply the optimal transport map
(cf. \cite{mccann1997convexity})
 from $g_X$ to $g_Y$
\begin{equation}
\label{eq:OT_gaussian}
T_{X,Y}(x)
=
\mu_Y \,+ \,
\Sigma_X^{-1/2} 
\left(
\Sigma_X^{1/2} 
\Sigma_Y
\Sigma_X^{1/2} 
\right)^{1/2}
\Sigma_X^{-1/2} 
(x-\mu_X).
\end{equation}
 to define the deformed point cloud
 $\widetilde{X} = 
\left\{
\widetilde{x}_i = 
T_{X,Y}(x_i)
\right\}_{i=1}^N$.

\subsection{Regression of sorted point clouds}
\label{sec:regression}

We rely on RBF regression to determine the approximation $\mu\in \mathcal{P} \mapsto \widehat{X}_{\mu} \in \mathbb{R}^{N\times d}$. Given the sorted point clouds
$\{ {X}_{\mu} : \mu\in \mathcal{P}_{\rm train}  \}$, we first apply POD to determine the equivalent representation
$$
X_{\mu^k}
=
\sum_{i=1}^M Z_i \beta_k^i,
\quad
{\rm with} \;
Z_1,\ldots,Z_M\in \mathbb{R}^{N\times d}, \;\;
 \beta_1^1,\ldots,
 \beta_{n_{\rm train}}^M \in \mathbb{R},
$$
and $M={\rm min} \{N\cdot d, n_{\rm train} \}$. Then, we apply RBF regression to each coefficient of the expansion separately,
\begin{equation}
\label{eq:RBF_regression}
\widehat{\beta}^i
=
{\rm arg} \min_{\beta\in \mathfrak{H}_{\phi}}
\lambda \|  \beta \|_{\mathfrak{H}_{\phi}}^2
\,+\,
\sum_{k=1}^{n_{\rm train}}
\left(
\beta(\mu^k) - \beta_k^i
\right)^2,
\quad
i=1,\ldots,M.
\end{equation}
Here, $\mathfrak{H}_{\phi}$ denotes the native space associated with the kernel $\phi$ and $\lambda>0$ is a regularization coefficient.
We rely on  
the Gaussian kernel, and we exploit cross-validation to select $\lambda$ and the kernel width; furthermore, we estimate the out-of-sample R-squared indicator 
and we keep the coefficients of the expansion for which the latter exceeds the threshold $0.5$. 
 In conclusion, we obtain the estimate
 \begin{equation}
 \label{eq:RBF_estimator_all}
 \widehat{X}_{\mu}
=
\sum_{i=1}^m Z_i \widehat{\beta}_{\mu}^i,
 \end{equation}
with $m\leq M$.

\subsection{Choice of the nearest neighbors}
\label{sec:nearest_neighbors}
Given the parameter value $\mu\in \mathcal{P}$, we select the parameters $\nu^1,\ldots,\nu^{\kappa}\in \mathcal{P}_{\rm train}$ that minimize the Euclidean distance 
${\rm dist}(\mu,\nu)= \|  \mu - \nu \|_2$.
The Euclidean distance is invariant under rototranslations of the parameters: as discussed in section \ref{sec:properties}, this property is key to prove the frame indifference of the CDI procedure.

\subsection{Boundary-aware registration}
\label{sec:registration}
The mapping $T_{X,Y}$  in \eqref{eq:OT_gaussian} does not preserve the boundary of the domain:
for each $\nu \in \mathcal{P}_{\rm nn}^{\mu}$, we  should hence find  a mapping $\Phi_{\nu}: \Omega\to \Omega$ that is bijective in $\Omega$ and satisfies
$\Phi_{\nu}(x_i^{\rm ref}) \approx \widehat{x}_{i,\nu}$ for $i=1,\ldots,N$.
Below, we present two registration techniques to accomplish this task. Note  that we here decouple the problem of PSR (\emph{point cloud matching}) from  the problem  of determining a bijective mapping in $\Omega$: this enables the use of off-the-shelves PSR routines for the matching problem.
The development of a general registration procedure for arbitrary domains is the subject of ongoing research and a thorough review of the subject is beyond the scope of the present paper.
It is straightforward to verify that both registration methods return the identity map if
$x_i^{\rm ref}=\widehat{x}_{i,\mu}$ for $i=1,\ldots,N$:
this observation is instrumental to prove the interpolation property of CDI in section \ref{sec:properties}.

\subsubsection{Elasticity-based registration}
\label{sec:elasticity_based}
We denote by $\widetilde{v} = T_{X,Y} - \texttt{id}$ the spatial displacement field associated with the PSR algorithm
where
$\texttt{id}$ is the identity map; we assume that $\widetilde{v}$ is defined over $\mathbb{R}^d$ --- this assumption holds for several PSR algorithms including the probabilistic method \eqref{eq:OT_gaussian} and its generalizations based on Gaussian mixture models (e.g.,
\cite{myronenko2010point}).
Given the reference point cloud $\{  x_i \}_{i=1}^N$, we define the  region
$A_{\eta}$  such that
\begin{subequations}
\label{eq:corrected_velocity_field}
\begin{equation}
\label{eq:corrected_velocity_field_Aeta}
A_{\eta} = \bigcup_{i=1}^N \mathcal{B}_{\eta}(x_i),
\end{equation}
where $\mathcal{B}_{\eta}(x)$ denotes the ball of radius $\eta$ centered in $x$; we introduce the regularized signed distance function
\begin{equation}
\label{eq:signed_distance}
H_{\eta,\delta}(x)
=
f_{\delta}
\left(
\phi(x, A_{\eta})
\right), \;  {\rm with} \;
\left\{
\begin{array}{l}
\displaystyle{
f_{\delta}(t) = 
\frac{1}{2} \left(
\tanh \left( \frac{t}{\delta}  \right) + 1
\right)}, \\[3mm]
\phi(x, A_{\eta})
=
\left\{
\begin{array}{ll}
{\rm dist}(x, \partial A_{\eta})
& x\in A_{\eta} \\
-{\rm dist}(x, \partial A_{\eta}),
& x\in \partial A_{\eta}. \\
\end{array}
\right.
\end{array}
\right.
\end{equation}
Finally,  we define the modified field
\begin{equation}
\label{eq:corrected_velocity_field_b}
\left\{
\begin{array}{ll}
\displaystyle{
\Delta v +   \nabla \left(  \nabla \cdot v \right)
=
\frac{1}{\epsilon} 
(v - \widetilde{v} )
H_{\eta,\delta} 
}
& x \in \Omega,  \\[3mm]
\displaystyle{
\mathbf{t} \cdot \left( \nabla v \mathbf{n}  \right)
= 0,
\;\;
v \cdot \mathbf{n} = 0
}
& x \in \partial \Omega,  \\
\end{array}
\right.
\end{equation}
where
$\mathbf{n}$ is the normal vector and $\mathbf{t}$ denotes the tangent vector to $\partial \Omega$. 
As in \cite{angot1999penalization}, the right-hand-side of \eqref{eq:corrected_velocity_field}$_1$ weakly enforces the condition $v=\widetilde{v}$ in the region $A_{\eta}$ in \eqref{eq:corrected_velocity_field_Aeta}; 
on the other hand, the condition $v \cdot \mathbf{n} = 0$ ensures the non-penetration condition, for $\epsilon\ll 1$. 
\end{subequations}

We observe that the field $v$ depends on the hyper-parameters $\epsilon,\eta,\delta$; we discuss their choice in the numerical experiments of section \ref{sec:nonlinear_interpolation_numerics}. Finally, we define the mapping $\Phi$ as follows:
\begin{equation}
\label{eq:mappingPhi}
\Phi(\xi) :=\mathfrak{X}(\xi, 1), \;\; {\rm where} \;\;
\left\{
\begin{array}{ll}
\partial_t \mathfrak{X}(\xi, t) = v( \mathfrak{X}(\xi, t) )
&
t\in (0,1),  \\
\mathfrak{X}(\xi, 0) = \xi .& \\
\end{array}
\right.
\end{equation}

We observe that the field $v$ corresponds to the Eulerian flow velocity, while 
$\mathfrak{X}(\xi, \cdot):(0,1) \to \mathbb{R}^d$ describes the temporal evolution of the particle located at $\xi$ at time $t=0$.
The condition 
$v \cdot \mathbf{n} = 0$ ensures that
$\mathfrak{X}(\xi, t)\in \Omega$ for all $t>0$ at the continuous level:
in practice, we discretize \eqref{eq:mappingPhi} using an explicit Euler scheme and we properly choose a time step $\Delta t$ to ensure that the geometry error is below a given threshold.
We remark that our approach does not ensure that
$\Phi(x_i)=y_i$ for $i=1,\ldots,N$ nor it ensures the bijectivity of the mapping  in $\Omega$; nevertheless, the numerical experiments of this manuscript show that in practice the method provides proper deformations of the domain $\Omega$ that approximately satisfy  the interpolation condition $\Phi(x_i) = y_i$ for $i=1,\ldots,N$. The development of more effective elasticity-based registration techniques as well as their rigorous analysis is the subject of ongoing research.

\subsubsection{Optimization-based registration}
\label{sec:optimization_registration}
An alternative approach consists in formulating the problem of registration as a minimization problem of the form
\begin{equation}
\label{eq:optimization_based_registration}
\min_{\Phi \in \mathcal{W}_{\Omega}} \;
\frac{1}{N} \sum_{i=1}^N \|  \Phi(x_i^{\rm ref}) - \widehat{x}_{i,\mu}  \|_2^2 \; + \;
\mathfrak{P}(\Phi),
\end{equation}
that can be solved using a gradient-based (quasi-Newton) method. The optimization statement in  \eqref{eq:optimization_based_registration} depends on the choice of the penalty term 
$\mathfrak{P}$ and of the search space
 $\mathcal{W}_{\Omega}$:
 the former  should  enforce local bijectivity in $\Omega$ and also promote the smoothness  (in a Sobolev sense) of  the mapping;
 the latter should   ensure that locally-bijectivity  --- i.e.,  ${\rm det}(\nabla \Phi) > 0$ in $\overline{\Omega}$ --- implies global  bijectivity in $\Omega$ and should be minimal for $\Phi=\texttt{id}$ to ensure the interpolation property of CDI.
We refer to 
\cite{taddei2020registration,ferrero2022registration}  for a detailed presentation of the registration procedure for two-dimensional domains: the generalization of the approach to arbitrary domains is the subject of ongoing research.

\subsection{Choice of the weights}
\label{sec:weights}
We determine the weights 
$\{ \omega_{\mu}^{\nu} : \nu\in \mathcal{P}_{\rm nn}^{\mu}  \}$
using the  inverse distance weighting (IDW, \cite{shepard1968two}) method. Given the nearest-neighbor training parameters $\nu^1,\ldots,$ $\nu^{\kappa} \in \mathcal{P}_{\rm train}$ and the new parameter $\mu\in \mathcal{P}$, we consider
\begin{subequations}
\label{eq:IDW}
\begin{equation}
\label{eq:IDWa}
\left\{
\begin{array}{ll}
\displaystyle{
\omega_{\mu}^{\nu} = 
\frac{  \overline{\omega}_{\mu}^{\nu}   }{
\sum_{\nu'\in    \mathcal{P}_{\rm nn}^{\mu}   }  \overline{\omega}_{\mu}^{\nu'}
},
\quad
{\rm with} \;\;
\overline{\omega}_{\mu}^{\nu'} = \frac{1}{{\rm dist}^p(\mu, \nu')},
}
&
{\rm if} \;\; \mu\notin \mathcal{P}_{\rm train};
\\[5mm]
\displaystyle{
\omega_{\mu}^{\nu} = 
\left\{
\begin{array}{ll}
1 & {\rm if} \; \mu=\nu \\[3mm]
0 & {\rm otherwise} \\
\end{array}
\right.
}
&
{\rm if} \;\; \mu\in \mathcal{P}_{\rm train}.
\\
\end{array}
\right.
\end{equation}
The algorithm can cope with arbitrary choices of the distance ${\rm dist}:\mathcal{P} \times \mathcal{P} \to \mathbb{R}_+$. We here consider the Euclidean distance
\begin{equation}
\label{eq:euclidean_distance}
{\rm dist}(\mu, \nu) = \| \mu - \nu \|_2.
\end{equation}
\end{subequations}
The power parameter $p\geq 1$ influences the relative importance of neighboring parameters: large values of $p$ assign greater influence to values closest to the interpolated parameter. Since we consider relatively small values of $\kappa$, we expect that the sensitivity with respect to $p$ is limited.
In the numerical experiments, we set $p=1$.

\begin{remark}
\label{remark:equivalenceCDI}
Let $\mathcal{P}\subset \mathbb{R}$ and let $\kappa=2$. Given $\mu\in [\nu^0,\nu^1]$ with $\nu_0,\nu_1\in \mathcal{P}_{\rm train}$, we define
$s=\frac{\mu-\nu^0}{\nu^1-\nu^0}$;
then, exploiting 
\eqref{eq:IDWa} with $p=1$, 
we obtain
$$
\widehat{u}_{\mu} = (1-s) \widetilde{u}_{\nu^0} + s \widetilde{u}_{\nu^1},
$$
which corresponds to the two-field CDI proposed in 
\cite{iollo2022mapping}, for a proper choice of the mappings.
\end{remark}

\subsection{Properties of the interpolation procedure}
\label{sec:properties}
The CDI \eqref{eq:CDI} with the choice of the weights \eqref{eq:IDW} has four  noteworthy  properties that are of interest for interpolation of fluid mechanics fields.
\begin{enumerate}
\item
\emph{Interpolation.}
Let $\mu\in \mathcal{P}_{\rm train}$. Then, $\widehat{u}_{\mu}=u_{\mu}$.
\begin{proof}
The proof  exploits the discussion of the previous sections.
Recalling the choice of $\mathcal{P}_{\rm nn}^{\mu}$ in section \ref{sec:nearest_neighbors}, we find $\mu \in \mathcal{P}_{\rm nn}^{\mu}$. Then, exploiting \eqref{eq:IDWa}, we find 
$\widehat{u}_{\mu} = \widetilde{u}_{\mu} = u_{\mu}\circ \Phi_{\mu}$.
Recalling the correction in section \ref{sec:regression}, we have $
\widehat{X}_{\mu}=X_{\mu}$  and hence (cf. section \ref{sec:registration}) $\Phi_{\mu}=\texttt{id}$.
\end{proof}
\item
\emph{Maximum principle.}
Let $b\in \mathbb{R}^D$;
then,  $\sup_{x\in \Omega} \big|  b\cdot \widehat{u}_{\mu}(x) \big|
\leq
\sup_{x\in \Omega,\nu\in \mathcal{P}} \big|  b\cdot u_{\nu}(x) \big|$ for all $\mu\in \mathcal{P}$.
\begin{proof}
We have
$$
\begin{array}{rl}
\displaystyle{\big|  b\cdot \widehat{u}_{\mu}(x) \big|}
&
\displaystyle{
\overset{\rm (i)}{\leq}
\sum_{\nu\in \mathcal{P}_{\rm nn}^{\mu}}
 \omega_{\mu}^{\nu}    
\big|  b\cdot \widetilde{u}_{\nu}(x) \big|
\overset{\rm (ii)}{\leq}
\max_{\nu\in \mathcal{P}_{\rm nn}^{\mu}}
\left(
\sup_{x\in \Omega}
\big|  b\cdot \widetilde{u}_{\nu}(x) \big|
\right)
}
\\[3mm]
&
\displaystyle{
\overset{\rm (iii)}{=}
\max_{\nu\in \mathcal{P}_{\rm nn}^{\mu}}
\left(
\sup_{x\in \Omega}
\big|  b\cdot  {u_{\nu}}(x) \big|
\right)
\leq
\sup_{\nu\in \mathcal{P}, x\in \Omega}
\big|  b\cdot  {u_{\nu}}(x) \big|
}.
\\
\end{array}
$$
In (i), we applied the definition of CDI, the triangle inequality and the positivity of the weights $\{ \omega_{\mu}^{\nu}  \}_{\nu}$;
in (ii), we used the fact that $\sum_{\nu\in \mathcal{P}_{\rm nn}^{\mu}}
   \omega_{\mu}^{\nu}      = 1$;
   in (iii), we used that 
   $\Phi_{\nu}(\Omega)=\Omega$ for all $\nu\in \mathcal{P}$.
 \end{proof}
\item
\emph{Minimum principle.}
Let $b\in \mathbb{R}^D$;
then, 
$\inf_{x\in \Omega}  b\cdot \widehat{u}_{\mu}(x)   
\geq 
\inf_{x\in \Omega, \nu\in \mathcal{P}}  b\cdot {u}_{\nu}(x)$ for all $\mu\in \mathcal{P}$.
\begin{proof}
Given $x\in \Omega, \mu\in \mathcal{P}$, 
we have  indeed
that 
$\omega_{\mu}^{\nu} \geq 0$,
$ \sum_{\nu  \in \mathcal{P}_{\rm nn}^{\mu}}
\omega_{\mu}^{\nu} = 1$, 
  and 
$ b \cdot   \widetilde{u}_{\nu} (x)    \geq c :=
\inf_{x\in \Omega, \nu\in \mathcal{P}}  \left( b\cdot {u}_{\nu}(x) \right) $
for all $\nu\in \mathcal{P}_{\rm nn}^{\mu}$. Therefore, 
$$
 b \cdot  \widehat{u}_{\mu}(x)  \ 
=
\sum_{\nu  \in \mathcal{P}_{\rm nn}^{\mu}}
\omega_{\mu}^{\nu} 
\left( b \cdot   \widetilde{u}_{\nu} (x)   \right) 
\geq
c \sum_{\nu  \in \mathcal{P}_{\rm nn}^{\mu}}
\omega_{\mu}^{\nu} 
= c.
$$
\end{proof}
\item
\emph{Frame indifference.}
Given the parametric field 
$\mu\mapsto u_{\mu}$, consider the rototranslation $p(\mu) = R \mu + b$ for some rotation matrix $R$ and vector $b\in \mathbb{R}^P$, and the parametric field
$p\mapsto v_p = u_{\mu(p)}$.
Then, 
$\widehat{u}_{p(\mu)} = \widehat{u}_{\mu}$ for all $\mu\in \mathcal{P}$.
\begin{proof}
It suffices to check that all the steps of the CDI procedure are invariant under roto-translations.
First, we observe that RBF interpolation/regression is invariant under rototranslations (cf. \cite{wendland2004scattered}).
Second, since the choice of the neighboring elements
and the choice of the CDI weights \eqref{eq:IDW}
are both based on  the Euclidean distance,  we find that both steps also invariant  under roto-translations.
Finally, the registration procedure is independent of the parameterization considered and in particular is 
invariant under roto-translations.
\end{proof}
\end{enumerate}

The maximum and minimum principles can be used to ensure that the CDI prediction is physically meaningful for out-of-sample parameters.
To provide a concrete reference, if we apply  CDI to  
the solution to the Euler or Navier-Stokes equations in primitive variables, we can guarantee that density and pressure are strictly positive and do not exceed the maximum value attained in the training set.
Similarly, since the mappings 
$\{\Phi_{\nu} \}_{\nu}$ in \eqref{eq:CDI} are bijective in $\Omega$, if the far-field conditions are parameter-independent and constant on $\partial \Omega$, we can exploit the same argument to prove that the CDI is consistent with the far-field conditions.

\section{Motivating example}
\label{sec:motivating_example}

\subsection{Model problem}
We consider the parametric problem:
\begin{equation}
\label{eq:laplace_equation}
\left\{
\begin{array}{ll}
\displaystyle{
-\frac{\partial^2 u_{\mu}}{\partial x^2}
=
f_{\mu}
}
&
{\rm in} \, \Omega=(-1,1),  \\[3mm]
 u_{\mu}(-1)=u_{\mu}(1)=0; &
\\
\end{array}
\right.
\;
f_{\mu}(x) = 
\frac{1}{\sigma}
{\rm exp} \left(
\frac{(x-\mu)^2}{\sigma^2}
\right),
\end{equation}
with $\mu\in \mathcal{P} = [-0.9,0.9]$ and $\sigma>0$. 
We introduce the  space $\mathcal{U} = H_0^1(\Omega)$  endowed with the inner product
\begin{subequations}
\begin{equation}
(w, v) = \int_{-1}^1  \partial_x w \partial_x v   \, +\,w v \,  dx,
\end{equation}
and the induced norm $\|  \cdot \| = \sqrt{(\cdot,\cdot)}$.
Exploiting integration by part, we find the variational problem associated with   \eqref{eq:laplace_equation}:
\label{eq:laplace_equation_variational}
\begin{equation}
{\rm find} \; u_{\mu} \in \mathcal{U} \; :  \;
a(u_{\mu}, v) = F_{\mu}(v)
\;\; \forall \, v \in  \mathcal{U}, 
\end{equation}
where
\begin{equation}
 a(w, v) = \int_{-1}^1  \partial_x w  \, \partial_x v \, dx,
 \quad
F_{\mu}(v) = \int_{-1}^1   f_{\mu} v \, dx.
\end{equation}
\end{subequations}
By tedious but straightforward calculations, we can show that
\begin{equation}
\label{eq:continuity_laplace}
\frac{\pi}{2+\pi} \| u \|^2
\leq
a(u,u)
\leq 
\| u \|^2
\quad
\forall \, u\in \mathcal{U}.
\end{equation}

\subsection{Linear method}

We denote by $\mathcal{Z}\subset \mathcal{U}$ an $n$-dimensional linear subspace of $\mathcal{U}$; we further introduce the Galerkin reduced-order model:
\begin{equation}
\label{eq:galerkinROM}
{\rm find} \; \widehat{u}_{\mu} \in \mathcal{Z} \; :  \;
a( \widehat{u}_{\mu}, v) = F_{\mu}(v)
\;\; \forall \, v \in  \mathcal{Z}.
\end{equation}
Exploiting C{\`e}a Lemma for symmetric coercive problems 
(e.g., \cite[Lemma 4.3]{quarteroni2009numerical}) and the estimate \eqref{eq:continuity_laplace}, we find
\begin{equation}
\label{eq:cea_lemma}
\|  \widehat{u}_{\mu} - {u}_{\mu} \|
\leq
\sqrt{1+\frac{2}{\pi}}
\min_{v\in \mathcal{Z}} \|  v - {u}_{\mu} \|,
\quad
\forall \, \mu \in \mathcal{P}.
\end{equation}
We conclude that for this model problem Galerkin projection achieves near optimal performance over the entire parameter domain, for any choice of the reduced space $\mathcal{Z}$.

\begin{figure}[H]
\centering

\subfloat[] 
{\includegraphics[width=0.5\textwidth]
 {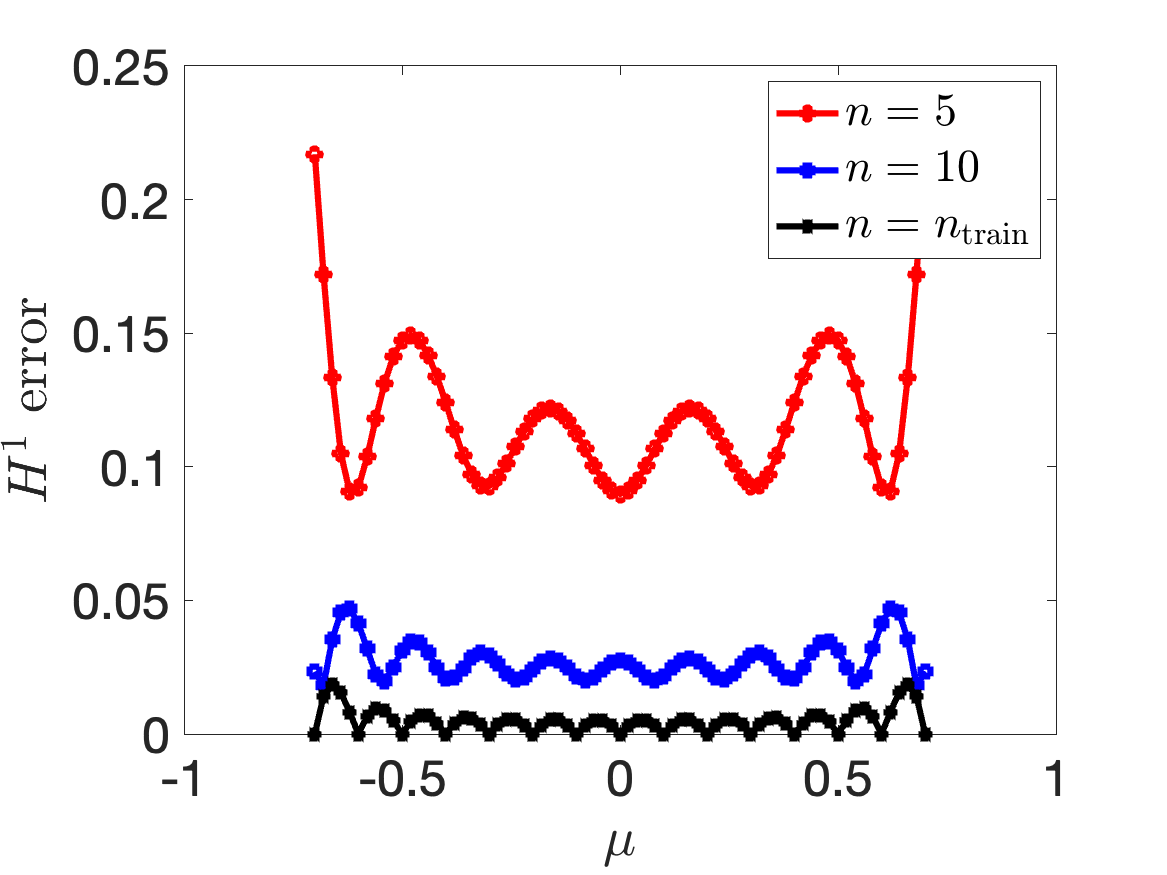}}
 ~~
   \subfloat[] 
{\includegraphics[width=0.5\textwidth]
 {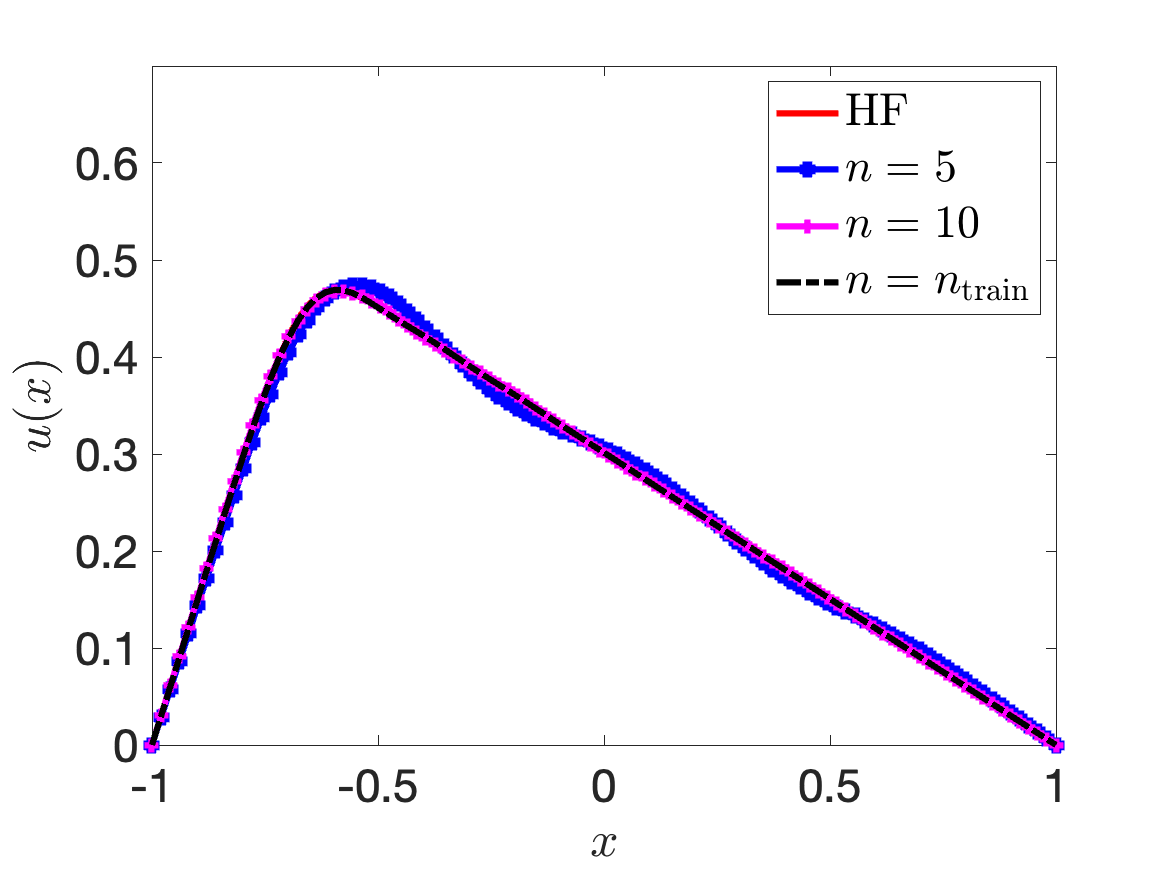}}
 
 \caption{motivating example; linear ROM for $\sigma=10^{-1}$. 
 (a) $H^1$ relative error.
 (b) comparison of HF and ROM solutions for the parameter that maximizes the error   for $n=n_{\rm train}$.
}
 \label{fig:motivating_linear_sigmam1}
 \end{figure}

\begin{figure}[H]
\centering

\subfloat[] 
{\includegraphics[width=0.5\textwidth]
 {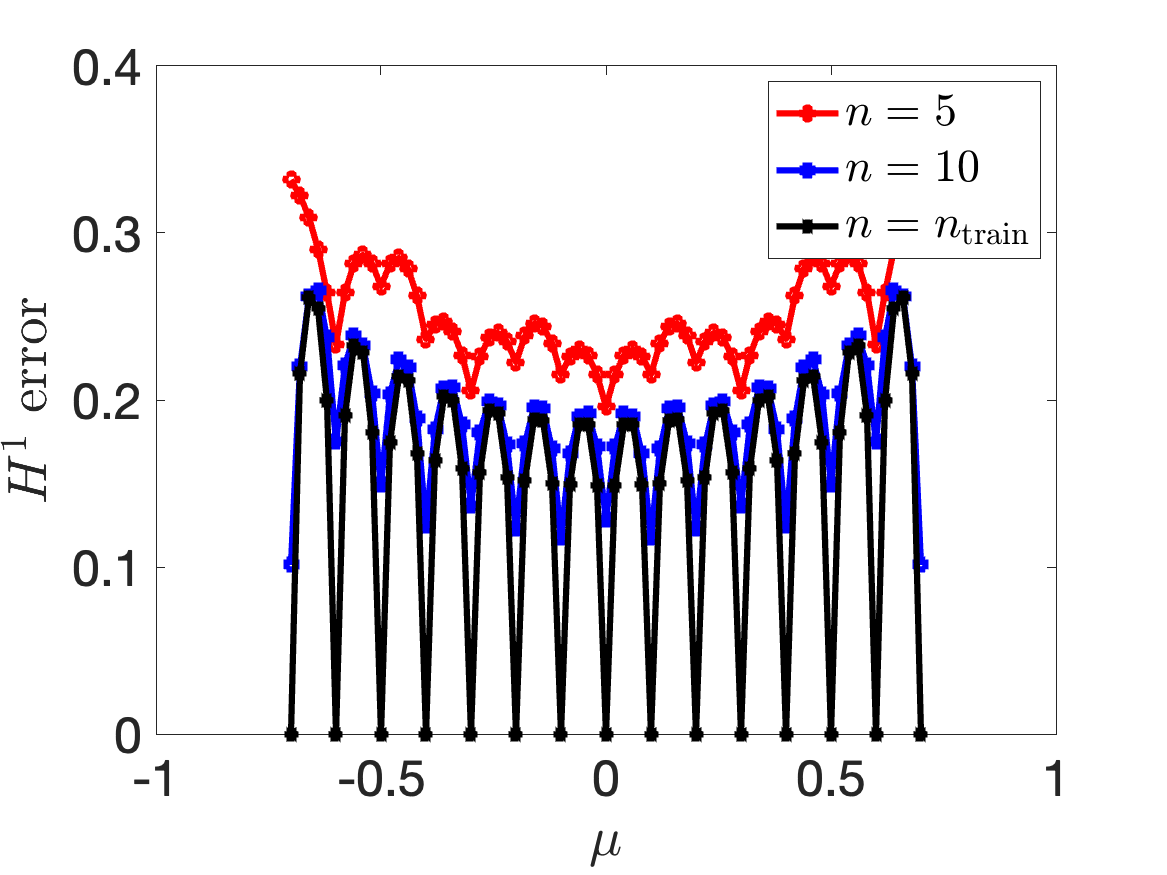}}
 ~~
   \subfloat[] 
{\includegraphics[width=0.5\textwidth]
 {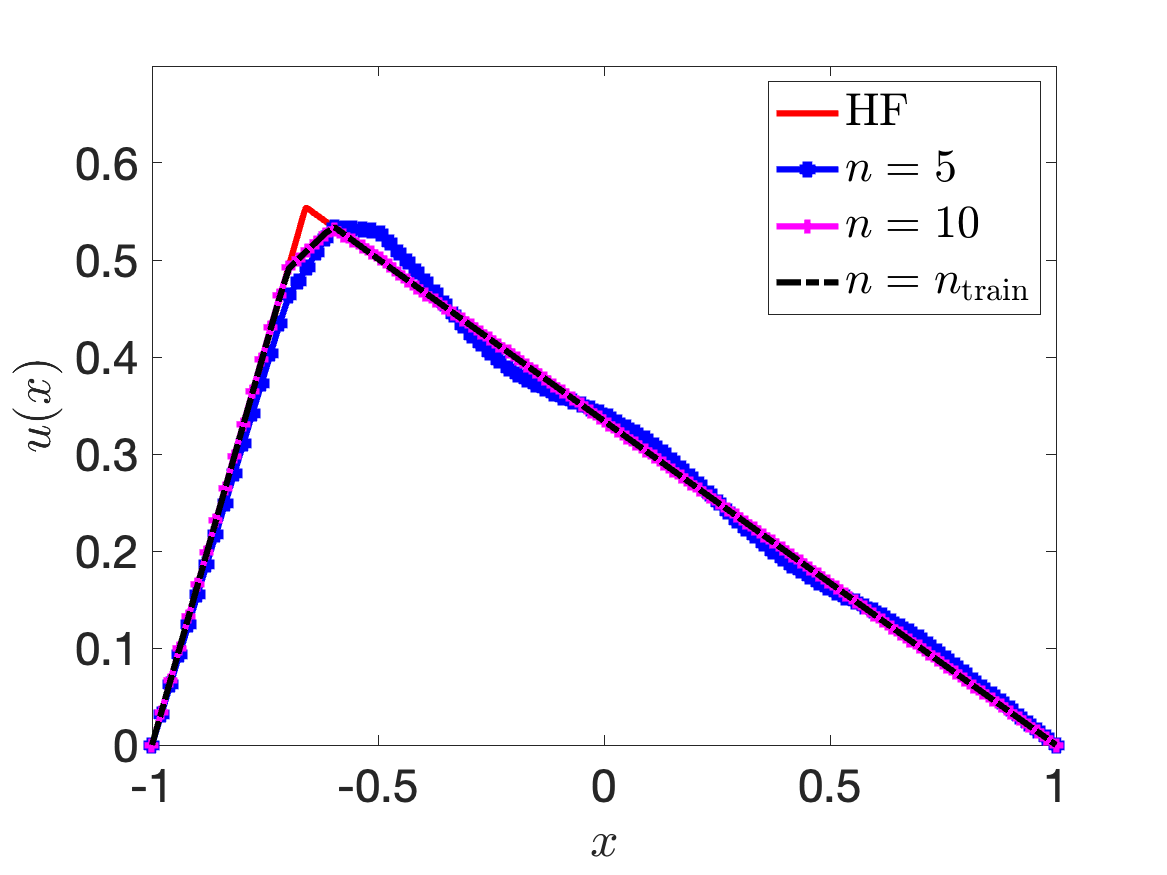}}
 
\caption{motivating example; linear ROM for $\sigma=10^{-3}$. 
(a) $H^1$ relative error.
(b) comparison of HF and ROM solutions for the parameter that maximizes the error   for $n=n_{\rm train}$.
}
 \label{fig:motivating_linear_sigmam3}
 \end{figure}

We generate $n_{\rm train}=15$ snapshots for equispaced parameter values  $\mathcal{P}_{\rm train}  \subset \mathcal{P}$ and for two values of $\sigma$, $\sigma=10^{-1}$ and $\sigma=10^{-3}$. For each value of $\sigma$, we assess performance of the POD Galerkin ROM \eqref{eq:galerkinROM} for three values of the reduced basis size $n$, $n=5,10,15$. 
Figure \ref{fig:motivating_linear_sigmam1}(a) shows the behavior of the relative $H^1$ error over $\mathcal{P}$ for three choices of $n$ for $\sigma=10^{-1}$, while Figure \ref{fig:motivating_linear_sigmam1}(b)
shows the solution and the Galerkin approximations for the parameter that maximizes the prediction error for $n=15$.
Figure \ref{fig:motivating_linear_sigmam3} replicates the results for $\sigma=10^{-3}$.

We observe that the ROM for $n=n_{\rm train}$ provides exact reconstructions for in-sample parameters: this immediately follows from \eqref{eq:cea_lemma}.
We note, however, that for $\sigma=10^{-3}$ the error is large --- in the order of $20\%$ --- and does not decrease as we increase $n$ from $n=10$ to $n=15$. 
We hence conclude that  the available snapshot set is not sufficient to achieve an accurate \emph{linear} representation of the elements of the solution manifold.

We draw two important conclusions. First, the inadequacy of linear approximations primarily depends on the geometry of the solution manifold --- which depends on the PDE model and on the parameterization --- rather than on the class of PDEs of interest (hyperbolic, elliptic, parabolic); on the other hand, the PDE model affects the sub-optimality (cf. \eqref{eq:cea_lemma}) of Galerkin projection.
Second, we observe that several nonlinear MOR methods (e.g.,
\cite{barnett2022quadratic,barnett2022neural}) seek state estimates that are nonlinear functions of the generalized coordinates but, nevertheless, belong to the span of the training snapshots:
 these methods hence  do not overcome the issue of the inaccuracy of the snapshot set.

\subsection{Convex displacement interpolation}

\begin{figure}[H]
\centering

\subfloat[] 
{\includegraphics[width=0.5\textwidth]
 {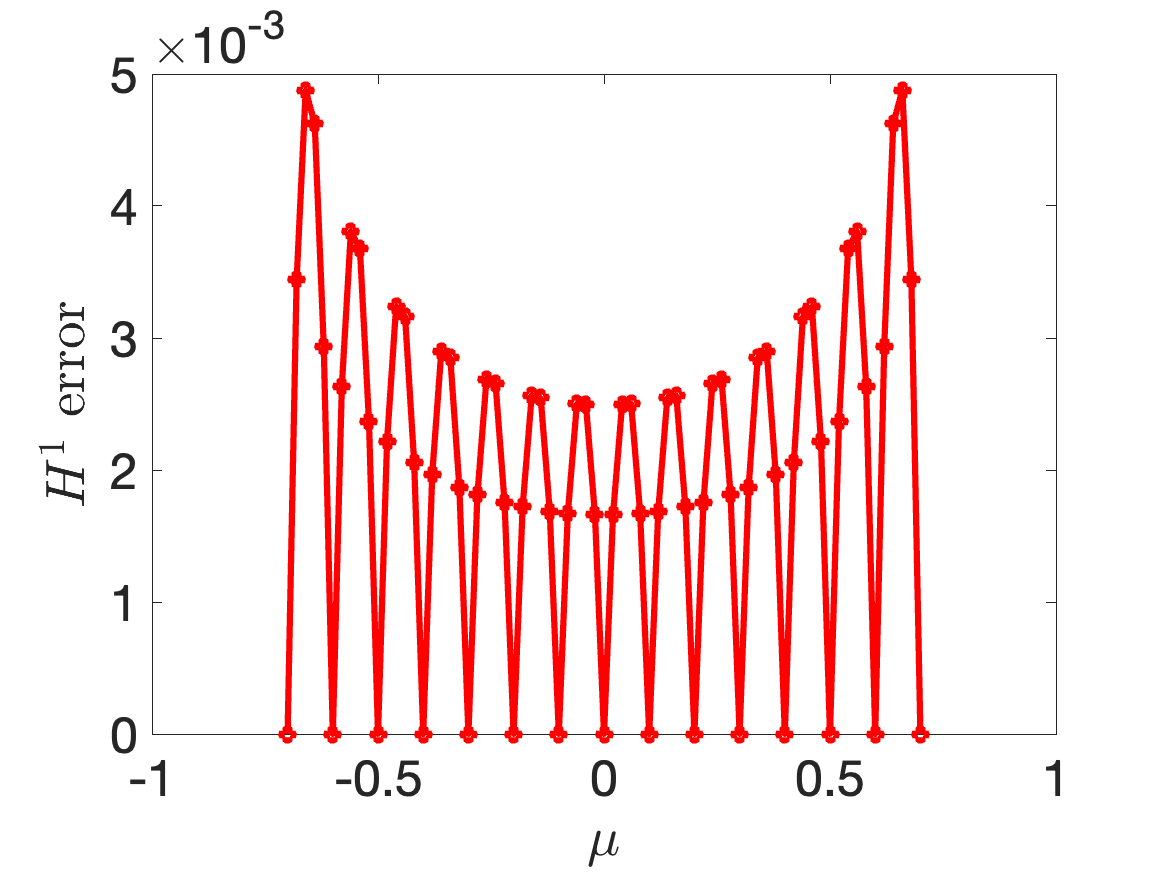}}
 ~~
   \subfloat[] 
{\includegraphics[width=0.5\textwidth]
 {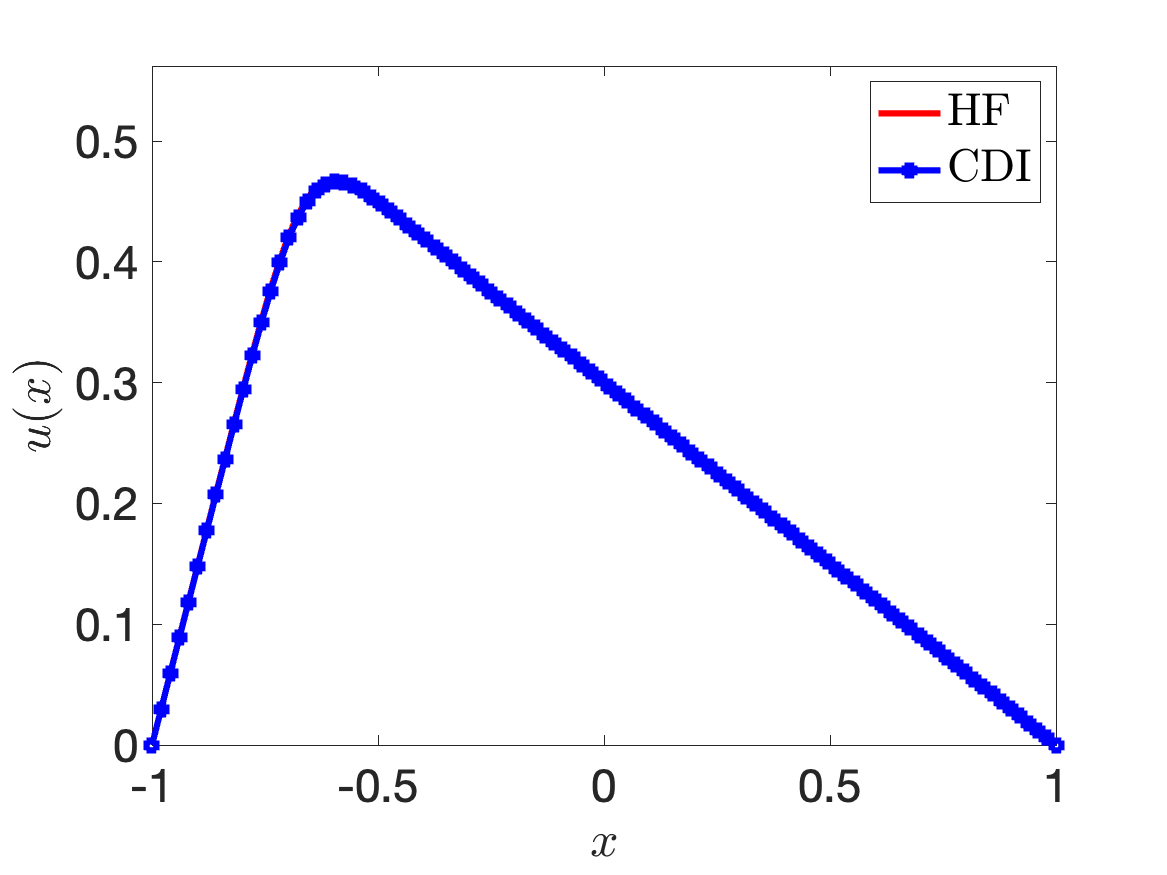}}
 
\caption{motivating example; CDI for $\sigma=10^{-1}$. 
(a) $H^1$ relative error.
(b) comparison of HF and ROM solutions for the parameter that maximizes the prediction error.
}
 \label{fig:motivating_CDI_sigmam1}
 \end{figure}

\begin{figure}[H]
\centering

\subfloat[] 
{\includegraphics[width=0.5\textwidth]
 {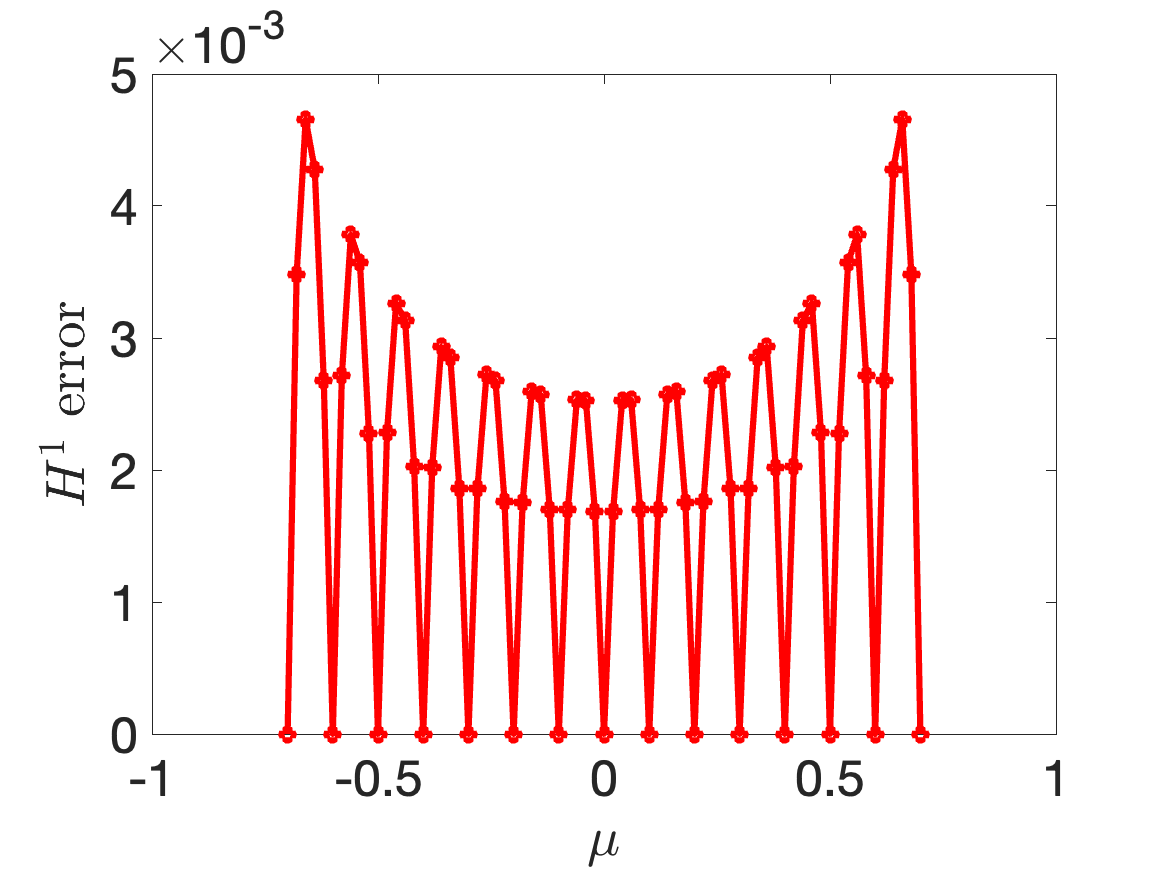}}
 ~~
   \subfloat[] 
{\includegraphics[width=0.5\textwidth]
 {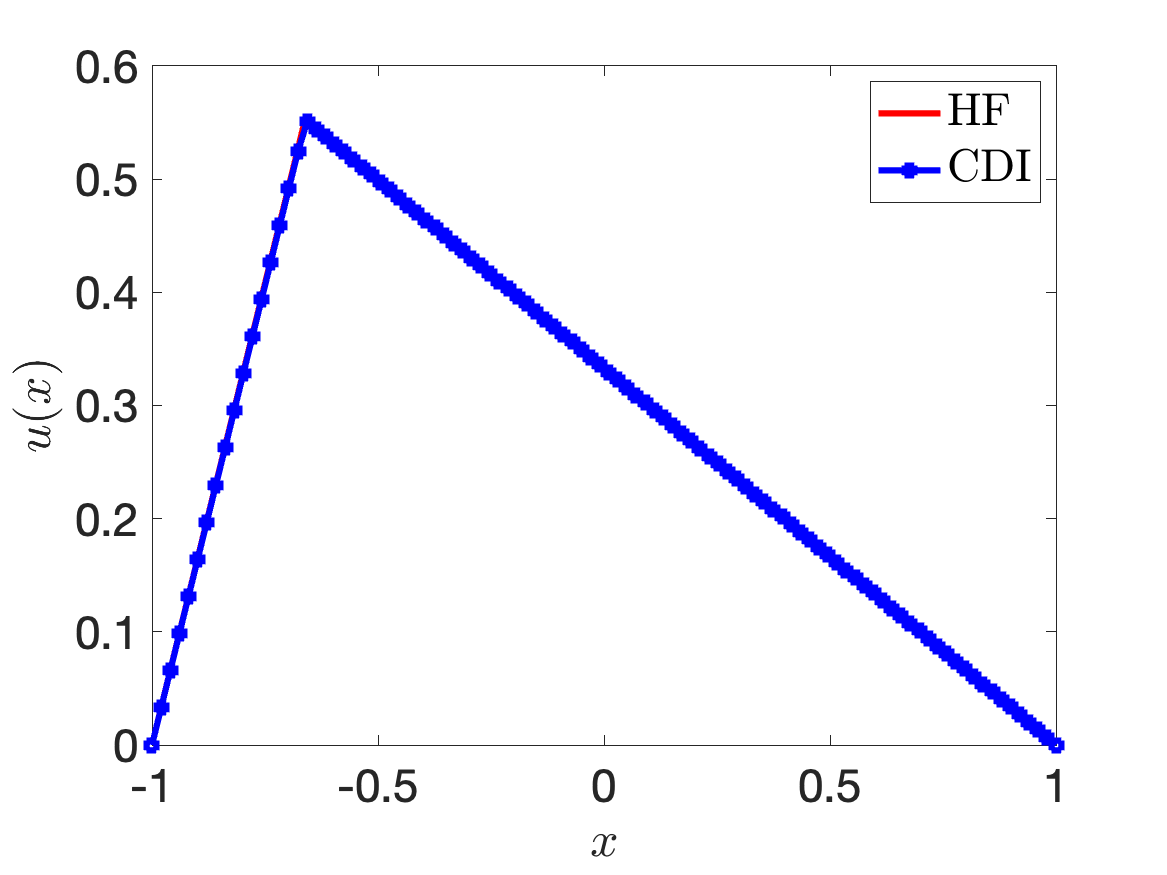}}
 
\caption{motivating example; CDI for $\sigma=10^{-3}$. 
(a) $H^1$ relative error.
(b) comparison of HF and ROM solutions  for the parameter that maximizes the prediction error.
}
 \label{fig:motivating_CDI_sigmam3}
 \end{figure}

Figures \ref{fig:motivating_CDI_sigmam1} and 
\ref{fig:motivating_CDI_sigmam3} show the performance of CDI for the model problem \eqref{eq:laplace_equation_variational} for $\sigma=10^{-1}$  and $\sigma=10^{-3}$:
Figures \ref{fig:motivating_CDI_sigmam1}(a) and 
\ref{fig:motivating_CDI_sigmam3}(a)  show the relative $H^1$ error, while 
Figures \ref{fig:motivating_CDI_sigmam1}(b) and 
\ref{fig:motivating_CDI_sigmam3}(b) show the behavior of the solution and its estimate for the   parameter that maximizes the prediction error.
We here set $X_{\mu}=\{ \mu \}$,
we consider $\kappa=2$ and 
 we define the bijective maps $\Phi_{\nu}$ for 
 $\nu\in \mathcal{P}_{\rm nn}^{\mu}$ so that 
 $\Phi_{\nu}(\Omega)=\Omega$,  $\Phi_{\nu}$ is piecewise linear in the intervals  $(-1,\mu)$ and $(\mu,1)$ and $\Phi_{\nu}(\mu)=\nu$.
We notice that the error is below $0.5\%$ for both values of $\sigma$.

\subsection{Data augmentation}
\begin{figure}[H]
\centering

\subfloat[] 
{\includegraphics[width=0.5\textwidth]
 {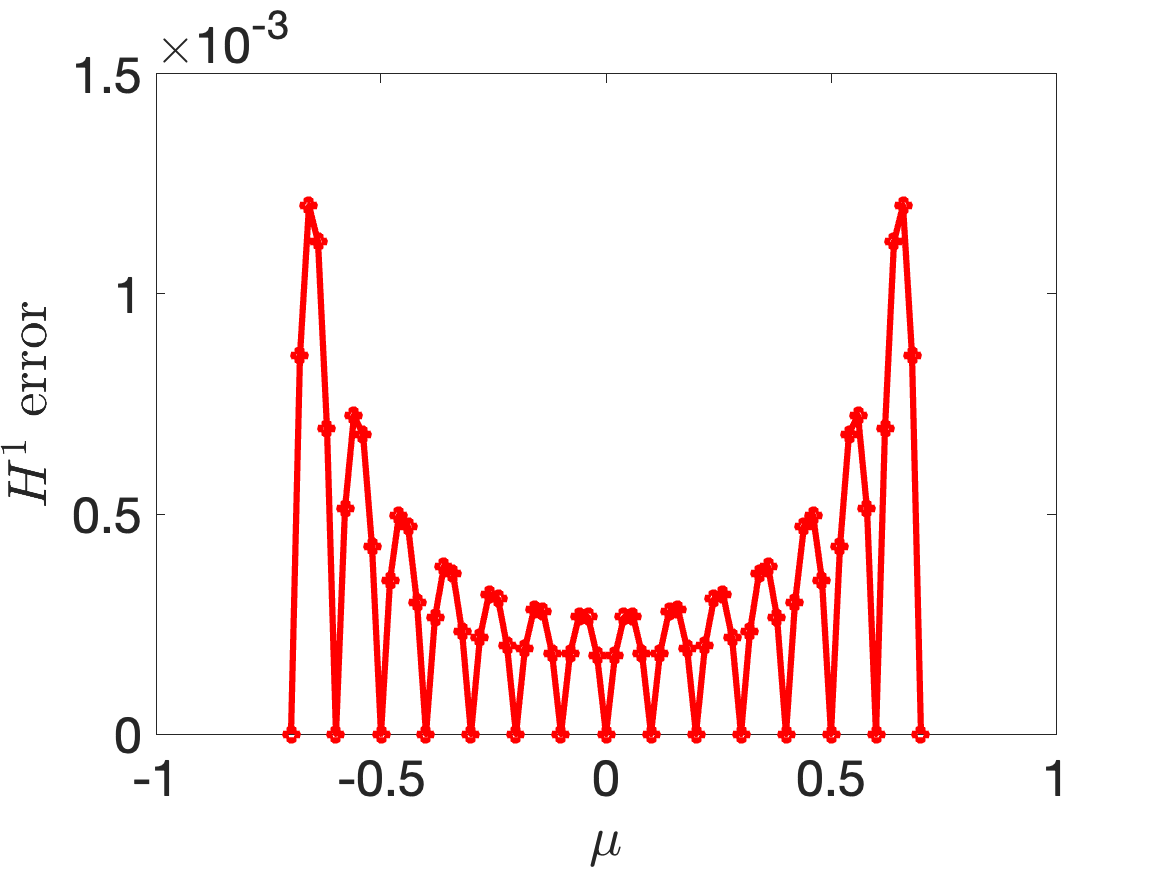}}
 ~~
   \subfloat[] 
{\includegraphics[width=0.5\textwidth]
 {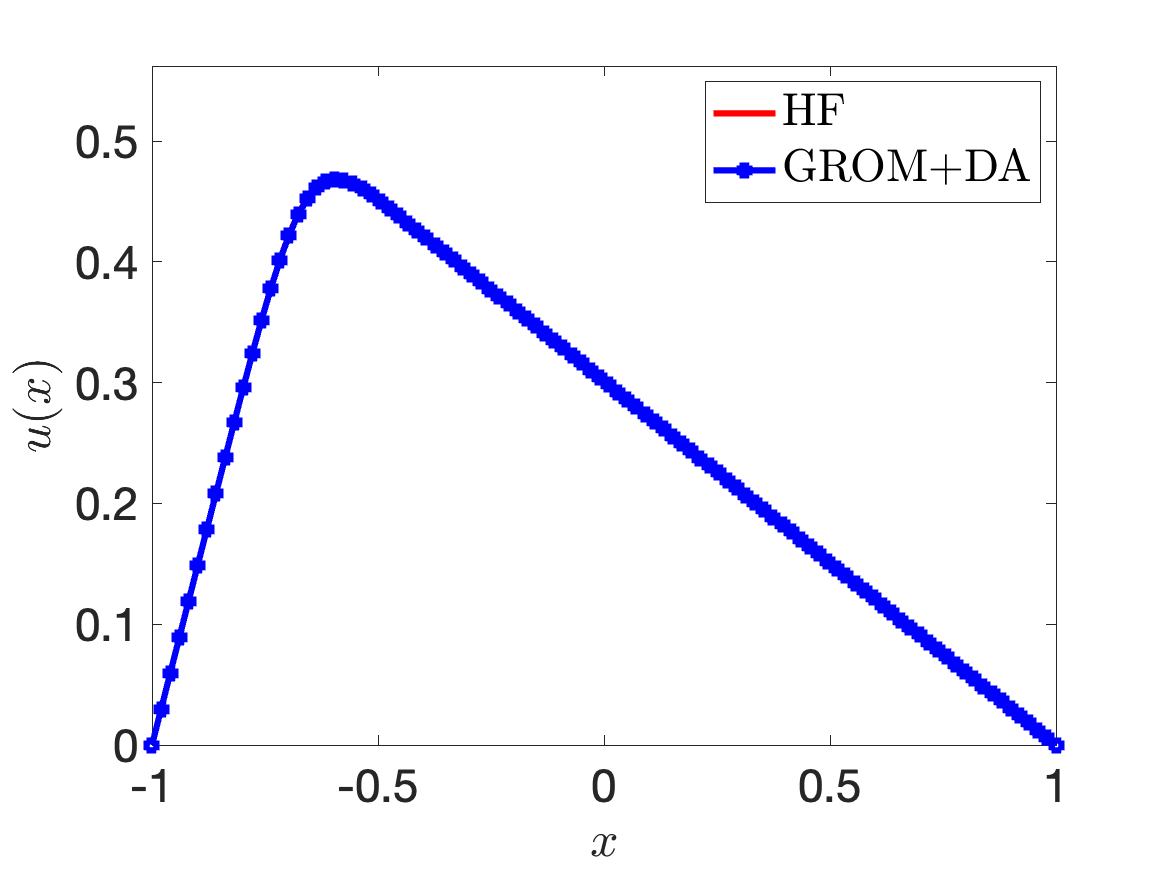}}
 
\caption{motivating example; Galerkin ROM with data augmentation for $\sigma=10^{-1}$. 
(a) $H^1$ relative error.
(b) comparison of HF and ROM solutions for the parameter that maximizes the prediction error.
}
 \label{fig:motivating_DAROM_sigmam1}
 \end{figure}

\begin{figure}[H]
\centering

\subfloat[] 
{\includegraphics[width=0.5\textwidth]
 {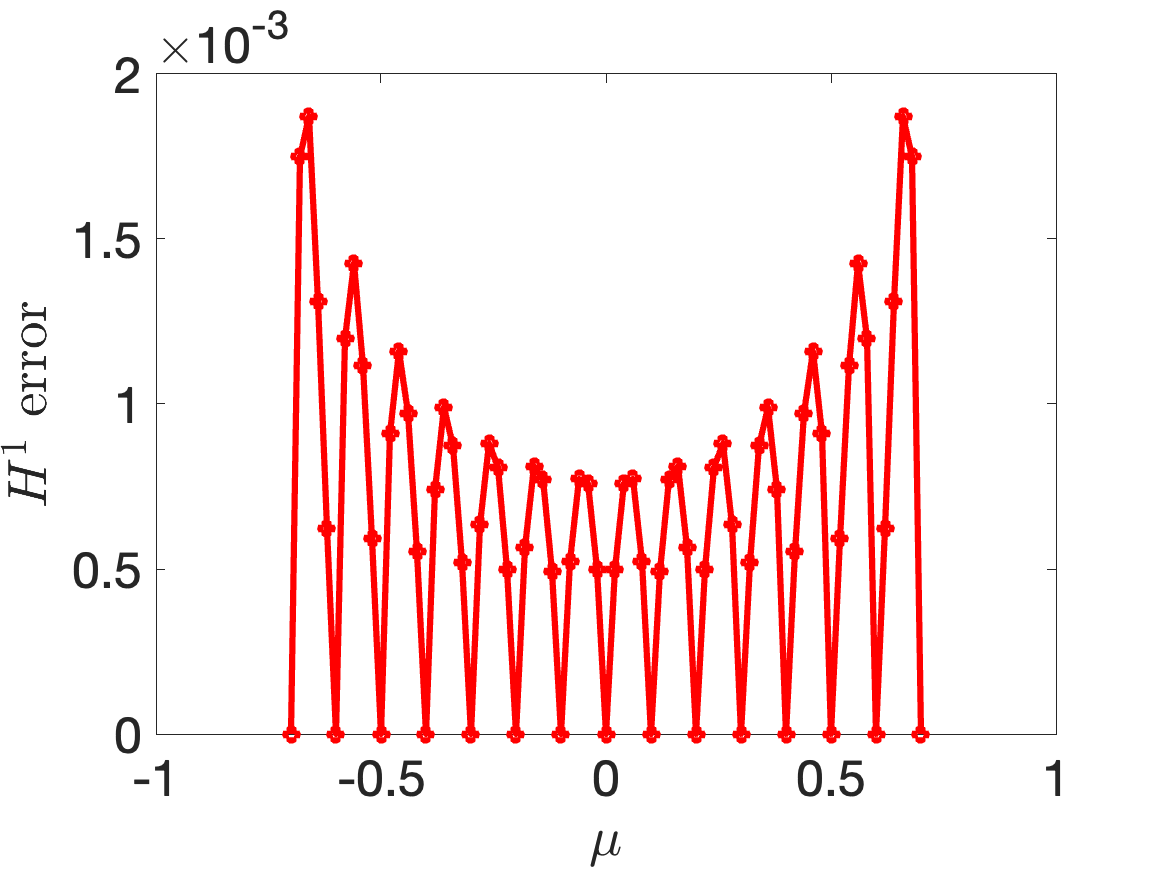}}
 ~~
   \subfloat[] 
{\includegraphics[width=0.5\textwidth]
 {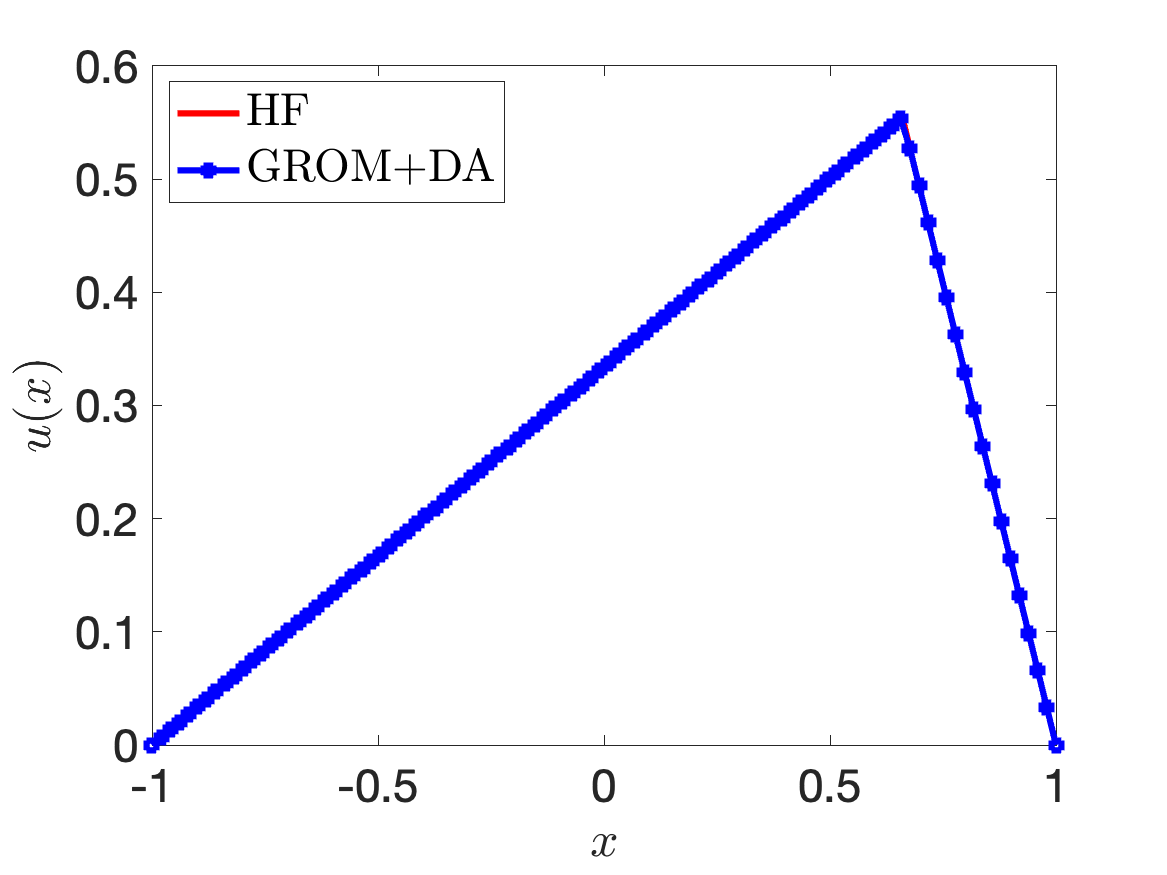}}
 
\caption{motivating example; Galerkin ROM with data augmentation for $\sigma=10^{-3}$. 
(a) $H^1$ relative error.
(b) comparison of HF and ROM solutions for the parameter that maximizes the prediction  error.
}
 \label{fig:motivating_DAROM_sigmam3}
 \end{figure}

As final test, we perform Galerkin projection based on the parameter-dependent reduced space
$$
\mathcal{Z}_{\mu} = {\rm span} \left\{ \widehat{u}_{\mu}^{\rm cdi}, 
{u}_{\nu^1}, \; 
{u}_{\nu^2}  \right\},
\quad
{\rm with} \;
\mathcal{P}_{\rm nn}^{\mu} = \{\nu^1,\nu^2\}.
$$
Figures \ref{fig:motivating_DAROM_sigmam1} and \ref{fig:motivating_DAROM_sigmam3} show the results: we observe that the error is roughly three times lower than the error associated with CDI. This result suggests the potential of data augmentation for Galerkin ROMs: nonlinear interpolation is employed to rapidly enhance the accuracy of the reduced space, while Galerkin projection is employed to ensure stable and near-optimal estimates.
We further investigate the combination of CDI and data augmentation in the next section.

\section{Nonlinear interpolation of parametric flow fields}
\label{sec:nonlinear_interpolation_numerics}

In this section, we illustrate the performance of our nonlinear interpolation procedure for three  model problems.
We rely on the elasticity-based  approach for registration  (cf. section \ref{sec:elasticity_based}):
in all our tests, we set $\epsilon=10^{-8}$, $\delta=50$, and $\eta=10^{-2}$ in  
\eqref{eq:corrected_velocity_field_Aeta} and \eqref{eq:corrected_velocity_field_b};
we further consider the time step $\Delta t=5 \cdot 10^{-3}$ for the integration of \eqref{eq:mappingPhi}.
We study the sensitivity of the interpolation parameter $s$, for datasets with two parameters $\mu^0,\mu^1$ and $P=1$. 
Given the point clouds
$X_i = \{ x_j^i  \}_{j=1}^N\subset \Omega$ for $i=0,1$, we define 
$X(s) = \{ (1-s) x_j^0  + s x_j^1   \}_{j=1}^N$ and the mappings 
$\Phi_i(s)$ such that
$\Phi_i(X(s), s    ) \approx X_i$ for $i=0,1$; finally, we define the two-field CDI
\begin{equation}
\label{eq:CDI_2field}
\widehat{u}(s) = 
(1-s) u_0\circ \Phi_0(s) + s  u_1\circ \Phi_1(s).
\end{equation}
and the linear counterpart (convex interpolation, CI)
\begin{equation}
\label{eq:CI_2field}
\widehat{u}(s) = 
(1-s) u_0 + s  u_1.
\end{equation}
In the absence of prior information, we here set $s=\frac{\mu-\mu_0}{\mu_1-\mu_0}$ (``linear relation'', cf. Remark \ref{remark:equivalenceCDI}); however, we also investigate the behavior of the optimal choice of $s$ as a function of $\mu$, 
\begin{equation}
\label{eq:relation_sopt}
s_{\mu}^{\rm opt} = {\rm arg}\min_{s\in [0,1]} \| \widehat{u}(s)  - u_{\mu} \|_{\star},
\end{equation}
 where the norm $ \| \cdot  \|_{\star}$ is introduced below for each model problem.

\subsection{Compressible flow past a RAE2822 airfoil}
\label{sec:RAE}

\subsubsection{Model problem}

We  study the compressible viscous flow past a RAE2822 airfoil for varying inflow Mach number  ${\rm Ma}_{\infty}$ and angle of attack $\alpha$.
We model the flow using the compressible Reynolds-averaged Navier-Stokes equations with  Menter's SST closure model 
\cite{menter1994two}; we rely on the 
open-source solver SU2 \cite{economon2016su2} for snapshot generation. In the remainder, we apply the nonlinear interpolation procedure to estimate the pressure coefficient 
\begin{equation}
\label{eq:pressure_coefficient}
C_{\rm p} = \frac{p-p_{\infty}}{\frac{1}{2}   \rho_{\infty} \| v_{\infty} \|_2^2},
\end{equation}
where $p,\rho, v$ denote flow pressure, density and velocity, and the subscript $(\bullet)_{\infty}$ denotes the free-stream value of the flow quantity $\bullet$.

\subsubsection{One-dimensional parameterization}

\begin{figure}[H]
\centering

\subfloat[] 
{\includegraphics[width=0.33\textwidth]
 {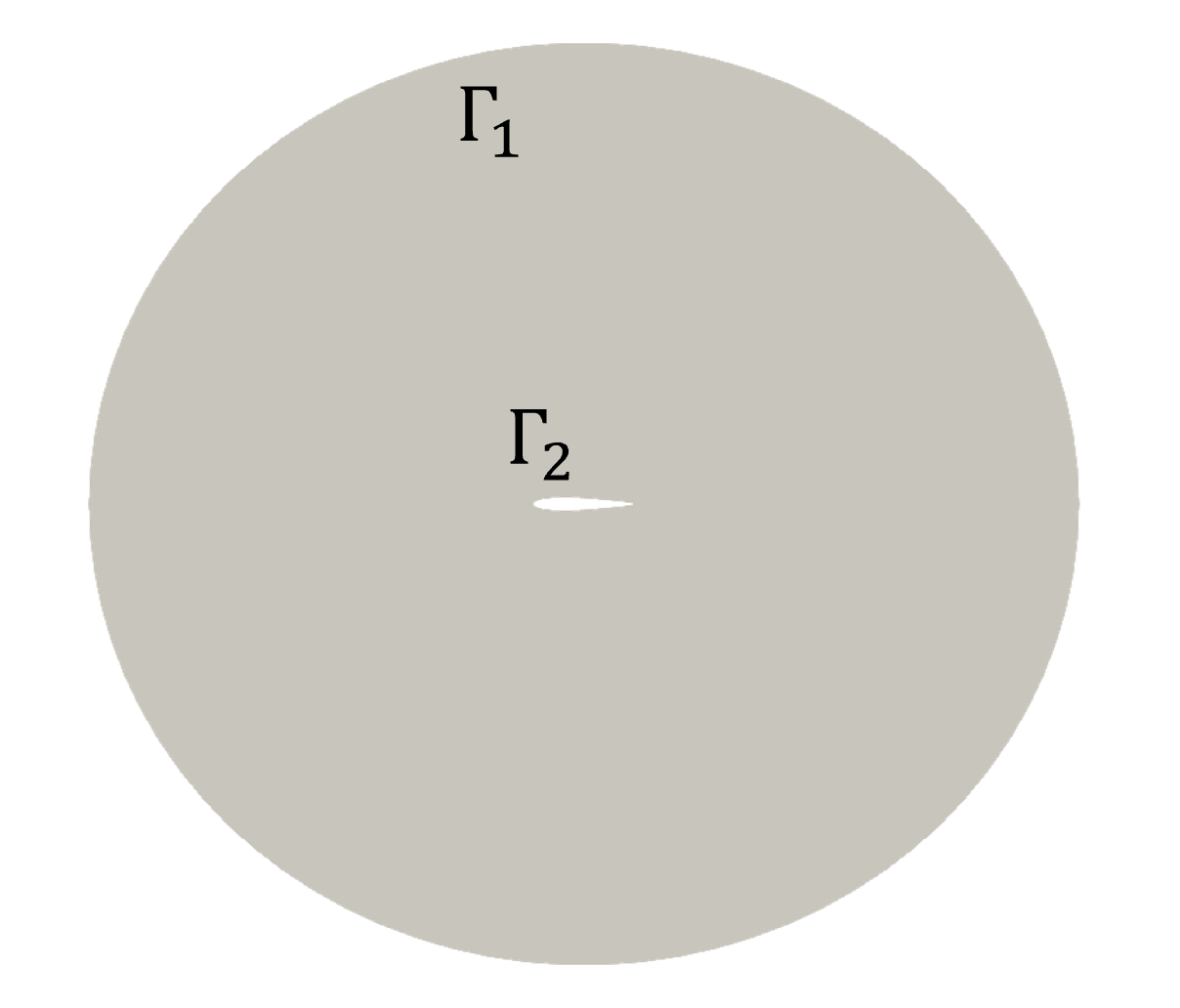}}
 ~~
   \subfloat[] 
{\includegraphics[width=0.33\textwidth]
 {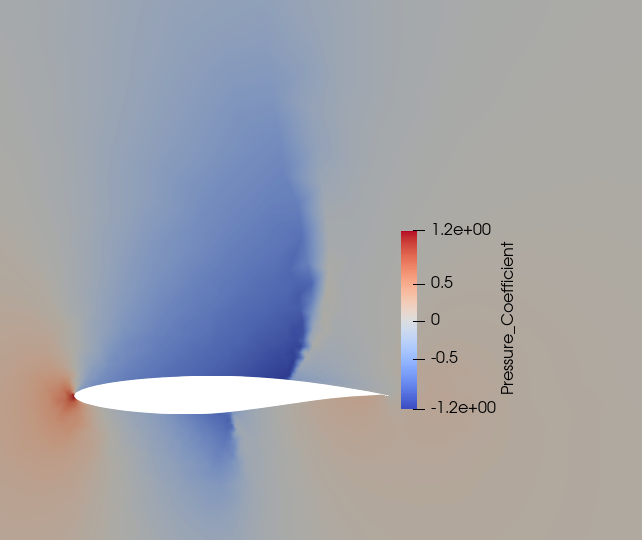}}
  ~~
   \subfloat[] 
{\includegraphics[width=0.33\textwidth]
 {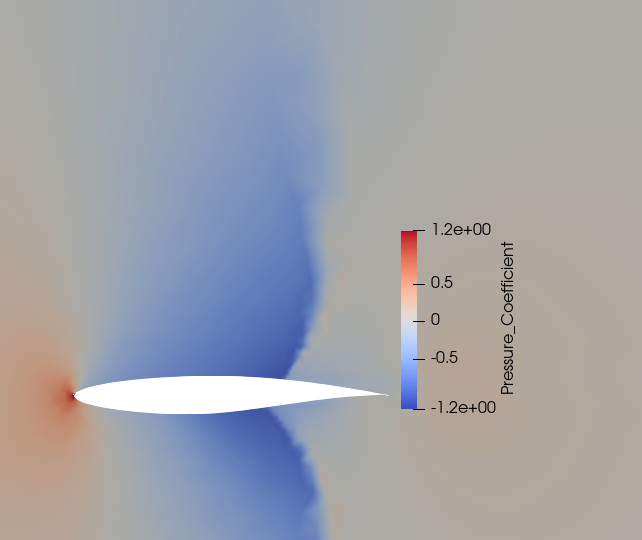}}
 
 \caption{Compressible flow past a RAE2822 airfoil. 
(a) computational domain.
(b)-(c) pressure coefficients for ${\rm Ma}_{\infty}=0.81$ and ${\rm Ma}_{\infty}=0.85$, $\alpha=2^o$.
}
 \label{fig:RAE_vis}
 \end{figure}

We first study variations of the free-stream Mach number in the range $[0.81,0.85]$ for $\alpha=2^o$.
Figure \ref{fig:RAE_vis}(a) shows the computational domain; 
Figures \ref{fig:RAE_vis}(b)  and (c) illustrate the behavior of the pressure coefficients for ${\rm Ma}_{\infty}=0.81$ and  ${\rm Ma}_{\infty}=0.85$. We observe that the solution develops two shocks on the lower and upper sides of the airfoil  whose position and intensity strongly depend on the parameter value.

\begin{figure}[H]
\centering

\subfloat[] 
{\includegraphics[width=0.5\textwidth]
 {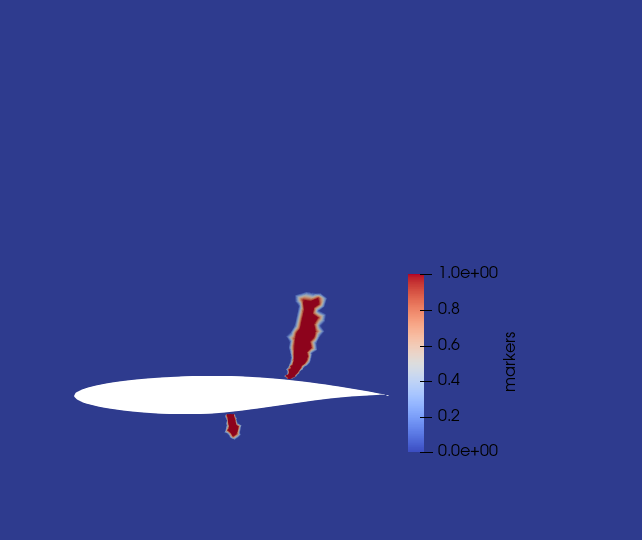}}
 ~~
\subfloat[] 
{\includegraphics[width=0.5\textwidth]
 {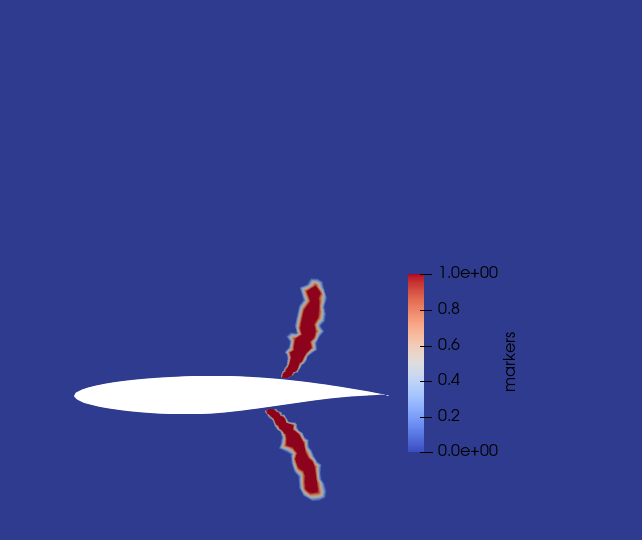}}
 
\caption{compressible flow past a RAE2822 airfoil. 
Marked cells for ${\rm Ma}_{\infty}=0.81$ and ${\rm Ma}_{\infty}=0.85$, $\alpha=2^o$.
}
\label{fig:RAE_markers}
\end{figure}

We rely on the Ducros sensor \eqref{eq:ducros_sensor}
to identify the markers, with 
$\gamma_{\rm thr}= 99.6\%$ (cf. \eqref{eq:quantile}). 
Figure \ref{fig:RAE_markers} shows the marked cells of the mesh for two parameter values.
Since the flow exhibits two distinct coherent structures (shocks), we apply the procedure  of section \ref{sec:PSR} to the lower and upper sides of the domain separately (see also \cite[section 5.2.2]{iollo2022mapping}).

\begin{figure}[H]
\centering

\subfloat[${\rm Ma}_{\infty}=0.82$] 
{\includegraphics[width=0.5\textwidth]
 {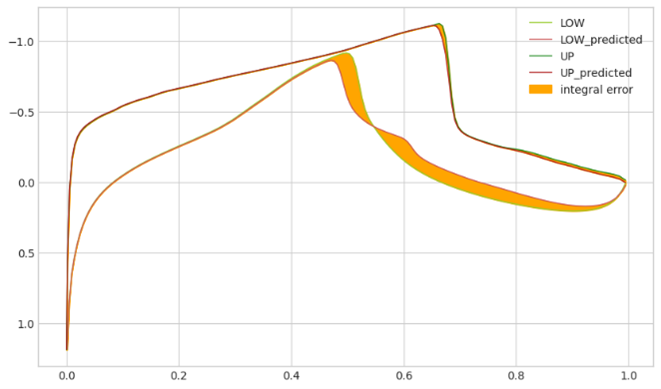}}
 ~~
\subfloat[${\rm Ma}_{\infty}=0.83$] 
{\includegraphics[width=0.5\textwidth]
 {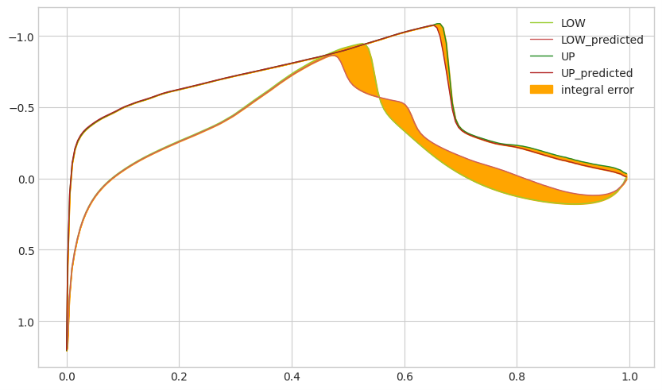}}

\subfloat[${\rm Ma}_{\infty}=0.82$] 
{\includegraphics[width=0.5\textwidth]
 {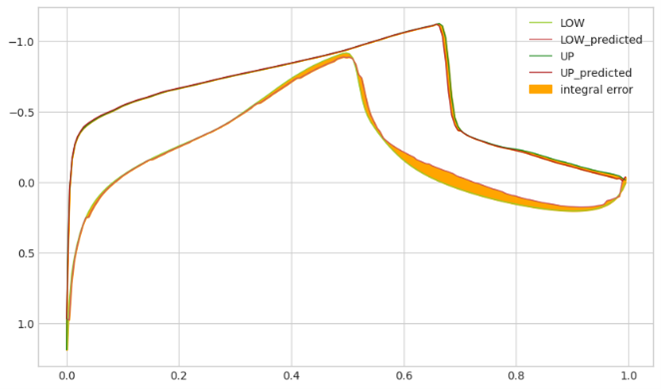}}
 ~~
\subfloat[${\rm Ma}_{\infty}=0.83$] 
{\includegraphics[width=0.5\textwidth]
 {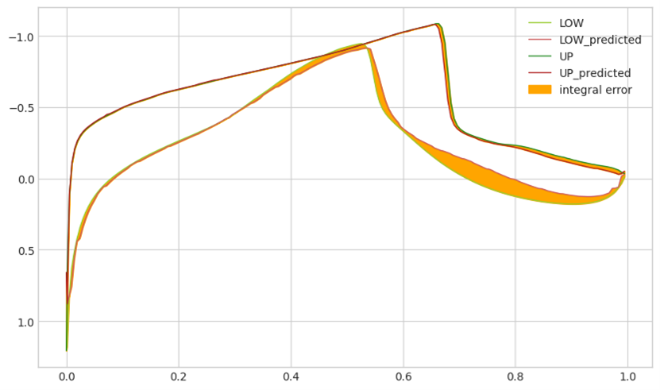}}

\caption{compressible flow past a RAE2822 airfoil, one-dimensional parameterization.
Behavior of truth and predicted pressure coefficient on lower (LOW) and upper (UP) sides of the blade, for two values of ${\rm Ma}_{\infty}$ and $\alpha=2^o$.
(a)-(b)  convex interpolation.
(c)-(d)  convex displacement interpolation.
}
 \label{fig:RAE_cp_2snaps}
 \end{figure}

\begin{figure}[H]
\centering

\subfloat[] 
{\includegraphics[width=0.5\textwidth]
 {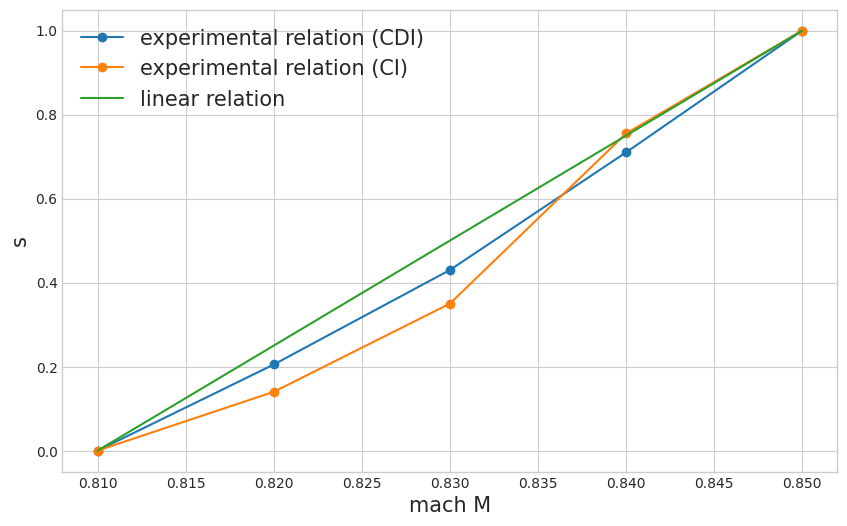}}
 ~~
\subfloat[] 
{\includegraphics[width=0.5\textwidth]
 {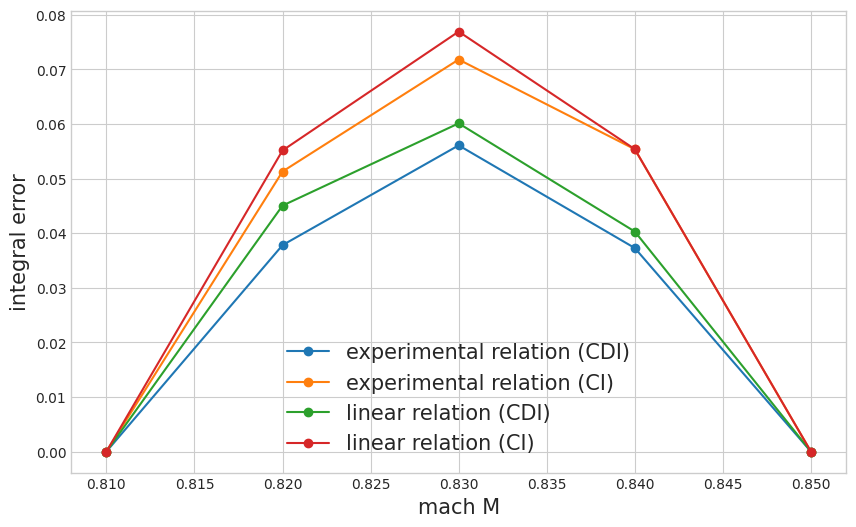}}
 
\caption{compressible flow past a RAE2822 airfoil,  one-dimensional parameterization.
(a) behavior of the optimal value of $s$ in 
\eqref{eq:relation_sopt} for CDI and CI.
(b) behavior of the integral error over the airfoil.
}
 \label{fig:RAE_error_2snaps}
 \end{figure}

Figure \ref{fig:RAE_cp_2snaps} investigates the performance of CDI \eqref{eq:CDI_2field} and CI \eqref{eq:CI_2field} over ${\rm Ma}_{\infty}\in [0.81,0.85]$ for linear choice of the interpolation parameter $s$ on the boundary $\Gamma_2$. 
We observe that CDI preserves the structure of the solution and is also significantly more accurate. Figure \ref{fig:RAE_error_2snaps} compares the behavior of the optimal choice of $s$ (``experimental'') with the linear relationship --- optimality is measured with respect to the $L^2(\Gamma_2)$  norm over the airfoil (``integral error''). We observe that for this model problem $s_{\mu}^{\rm opt}$
varies linearly with respect to the parameter $\mu={\rm Ma}_{\infty}$.

\subsubsection{Two-dimensional parameterization}
\begin{figure}[H]
\centering

\subfloat[linear] 
{\includegraphics[width=0.5\textwidth]
 {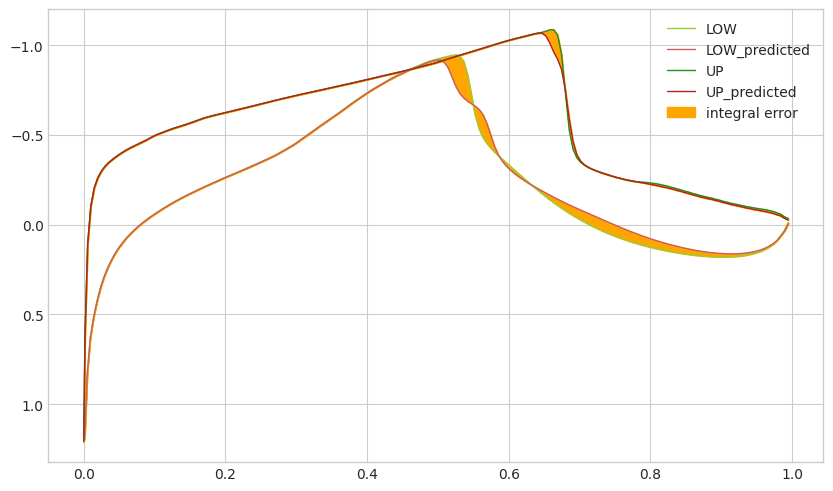}}
 ~~
\subfloat[nonlinear] 
{\includegraphics[width=0.5\textwidth]
 {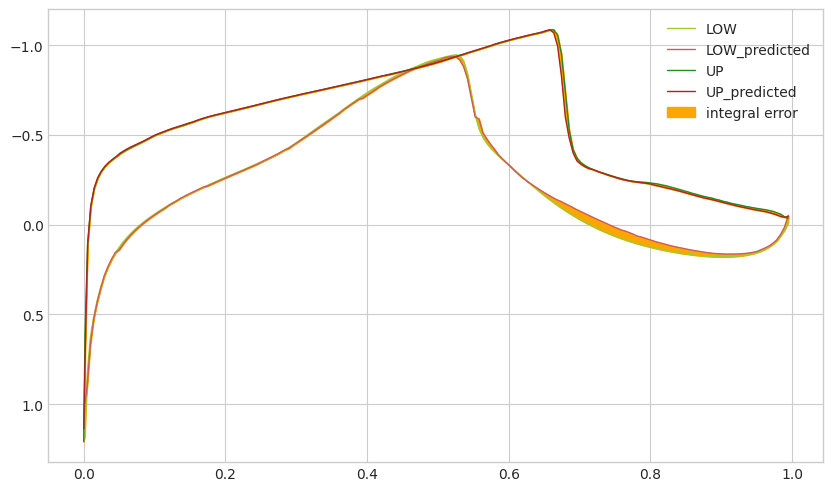}}
 
\caption{compressible flow past a RAE2822 airfoil, two-dimensional parameterization.
Behavior of truth and predicted pressure coefficient for an out-of-sample configuration.
(a) convex interpolation.
(b) convex displacement interpolation.
}
 \label{fig:RAE_cp_3snaps}
 \end{figure}

We consider the problem of estimating the solution for $\mu=({\rm Ma}_{\infty}, \alpha)=[0.83,2]$, based on the training points
$\mu^1=[0.82,2]$, 
$\mu^2=[0.84,1.5]$ and 
$\mu^3=[0.84,2.5]$.
We compare the performance of convex interpolation and CDI. Given the (very) modest amount of datapoints, we rely on linear interpolation for the regression step (cf. Algorithm \ref{alg:CDI_overview}, Line 1, Online stage). 
Figure \ref{fig:RAE_cp_3snaps} compares the predictions of truth and predicted pressure coefficient for the out-of-sample parameter $\mu$.
We observe that   nonlinear interpolation  
is qualitatively more accurate, especially in the neighborhood of the shock. 
We observe that the $L^2(\Omega)$ error is $5.62\%$ for CI and $4.10\%$ for CDI; the 
$L^2(\Gamma_2)$ error is $2.69\%$ for CI and $2.14\%$ for CDI. We emphasize that the  $L^2$ metric does not  fully capture local improvements in the regions where the shocks develop.

\subsection{Incompressible flow past a backward facing step}
\label{sec:backwardstep}

\subsubsection{Model problem}

We study the incompressible flow past a backflow facing step for varying inflow Reynolds number ${\rm Re}$ and step's inclination angle $\alpha$; we model the flow using the incompressible Navier-Stokes equations and we rely on the open-source software  OpenFOAM 
\cite{jasak2009openfoam}. Mesh morphing is based on the Optimad proprietary software \texttt{mimic}\footnote{See \url{https://optimad.github.io/mimmo/}.}. The coherent structures of interest are recirculation bubbles (vortices), whose extensions and positions vary as a function of the parameters. We apply the nonlinear interpolation procedure to the velocity field.

\subsubsection{One-dimensional parameterization}
 
\begin{figure}[H]
\centering

\subfloat[${\rm Re}=50$] 
{\includegraphics[width=0.69\textwidth]
 {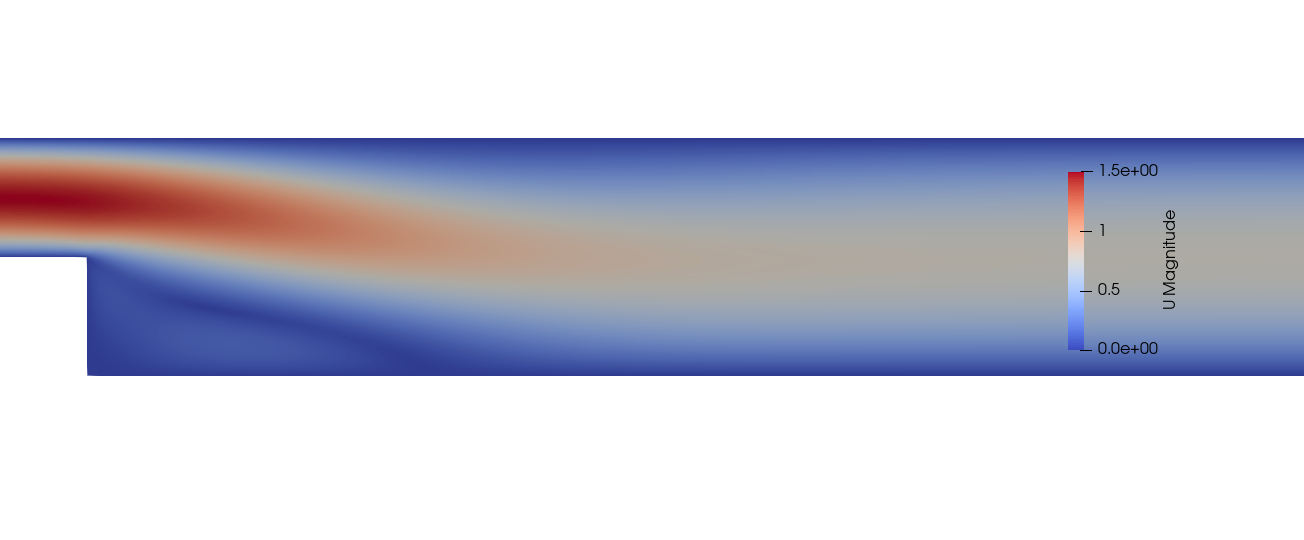}}

\subfloat[${\rm Re}=200$] 
{\includegraphics[width=0.69\textwidth]
 {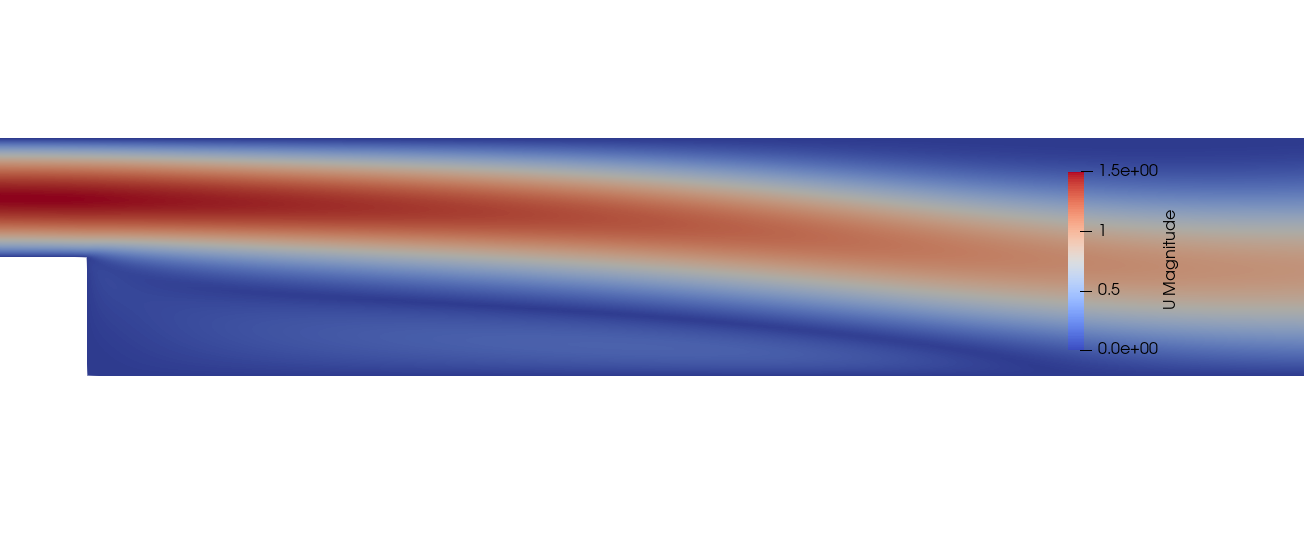}}
 
\caption{incompressible flow past a backward facing step. 
Velocity magnitude for two  values of ${\rm Re}$ and $\alpha=90^o$.
}
 \label{fig:BFS_vis_snaps}
 \end{figure}

\begin{figure}[H]
\centering

\subfloat[${\rm Re}=50$] 
{\includegraphics[width=0.69\textwidth]
 {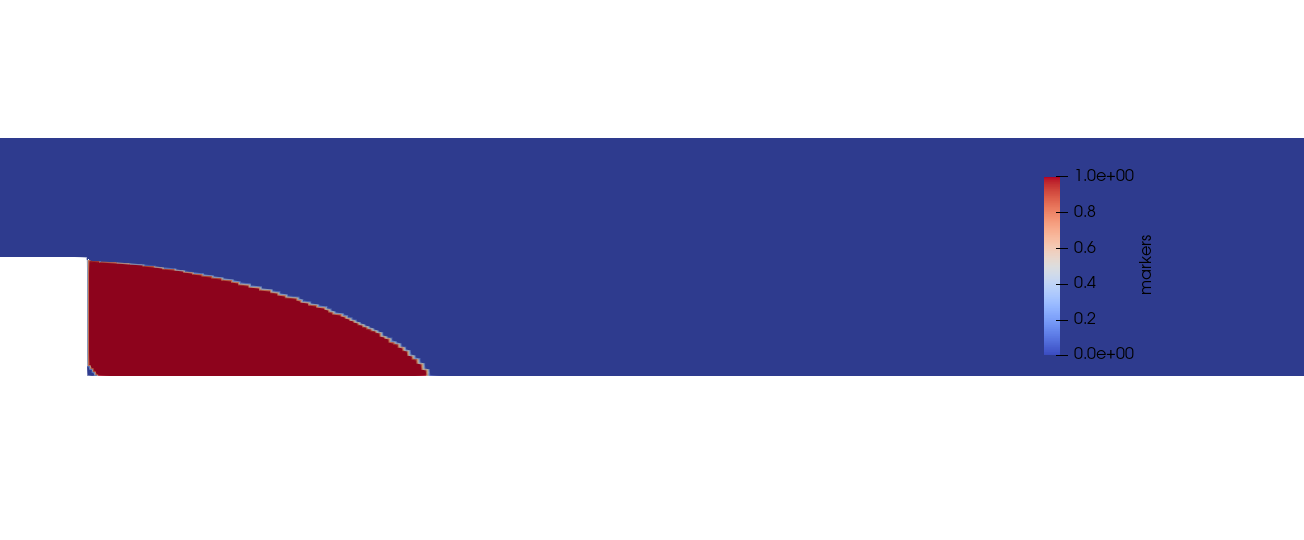}}

\subfloat[${\rm Re}=200$] 
{\includegraphics[width=0.69\textwidth]
 {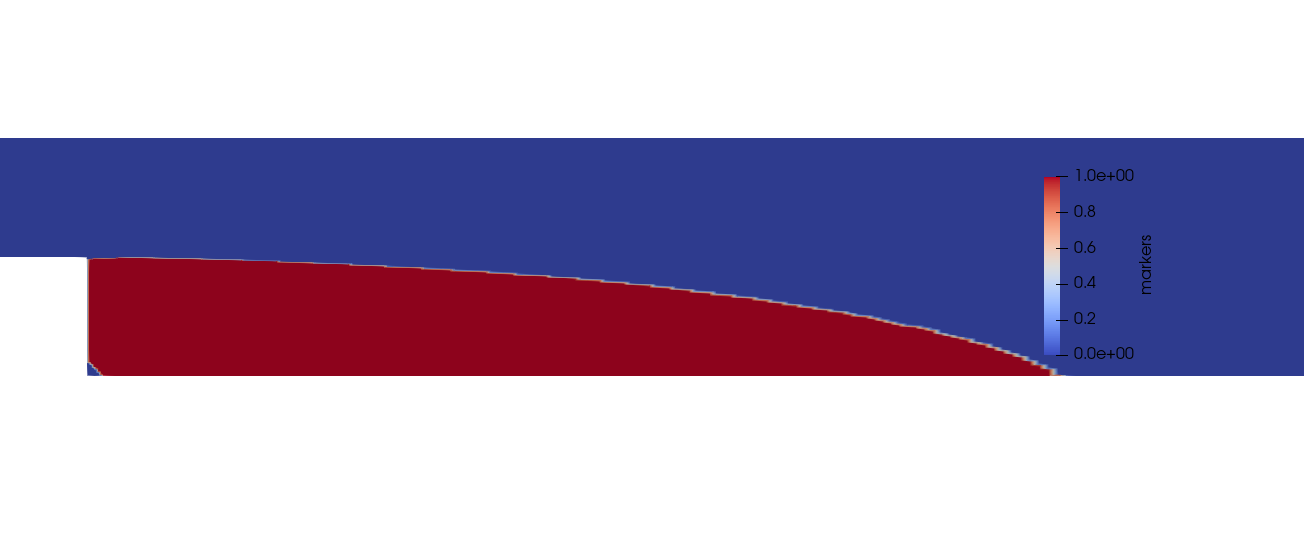}}
 
\caption{incompressible flow past a backward facing step. 
Marked cells for two  values of ${\rm Re}$ and $\alpha=90^o$.
}
 \label{fig:BFS_vis_markers}
 \end{figure}

 We first study variations of the inflow Reynolds number in the range $[50,200]$ for $\alpha=90^o$.
 Figure \ref{fig:BFS_vis_snaps} shows the velocity magnitude for two parameter values, while Figure \ref{fig:BFS_vis_markers} shows the recirculation area as predicted by the indicator 
\eqref{eq:streamline_function}. We observe that the recirculation area is very sensitive to the value of the 
inflow Reynolds number.

\begin{figure}[H]
\centering

\subfloat[high-fidelity model] 
{\includegraphics[width=0.73\textwidth]
 {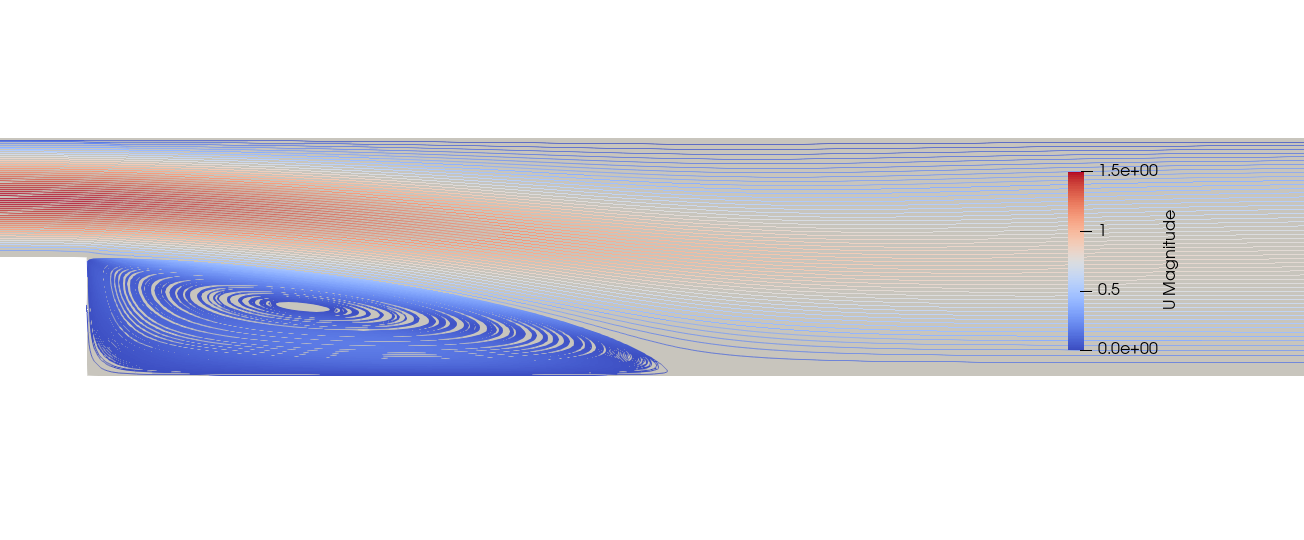}}

\subfloat[nonlinear interpolation] 
{\includegraphics[width=0.73\textwidth]
 {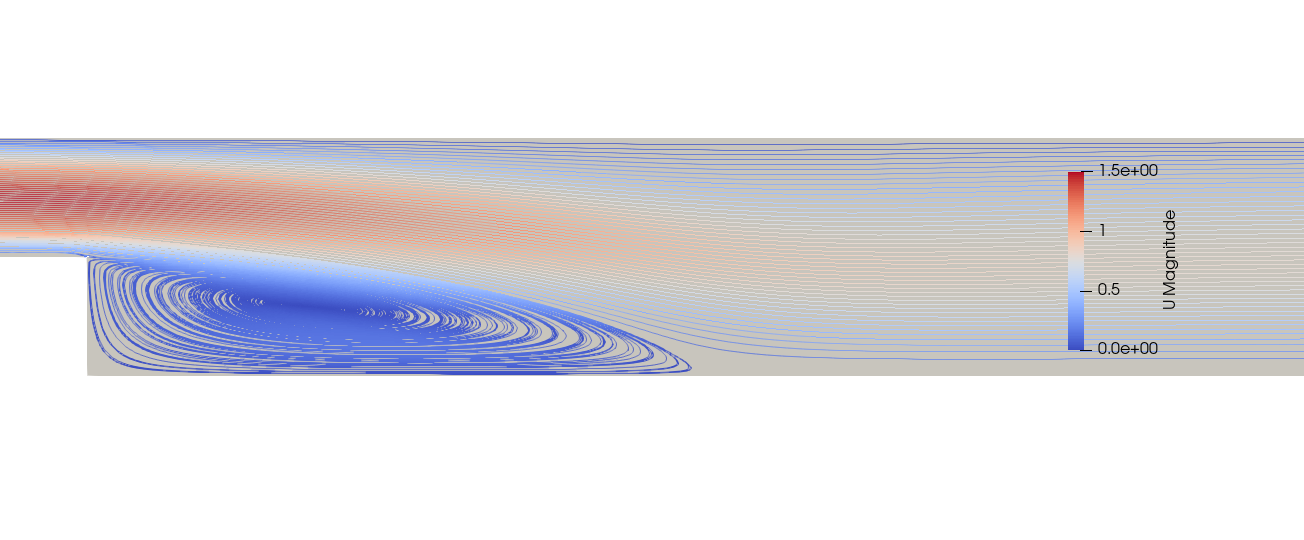}}
 
\subfloat[linear interpolation] 
{\includegraphics[width=0.73\textwidth]
 {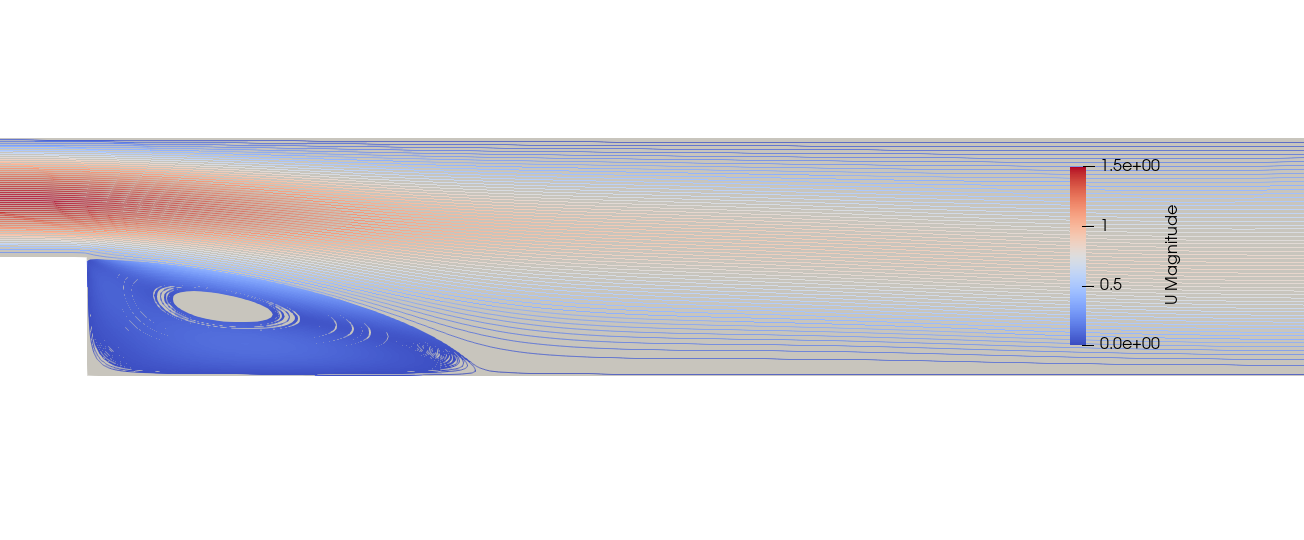}}
 
\caption{incompressible flow past a backward facing step, one-dimensional parameterization.
Predictions of velocity streamlines for ${\rm Re}=100$ based on the high-fidelity model, nonlinear interpolation and linear interpolation.}
 \label{fig:BFS_prediction_streamlines}
 \end{figure}

\begin{figure}[H]
\centering

\subfloat[] 
{\includegraphics[width=0.5\textwidth]
 {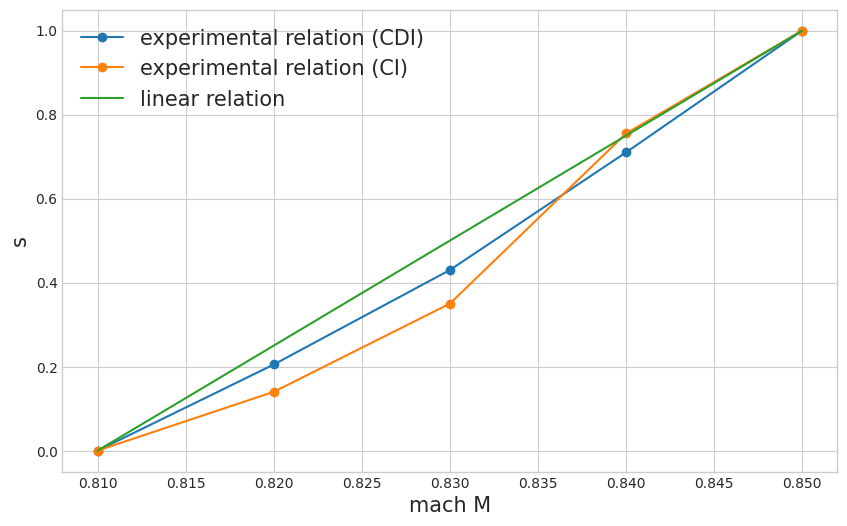}}
 ~~
\subfloat[] 
{\includegraphics[width=0.5\textwidth]
 {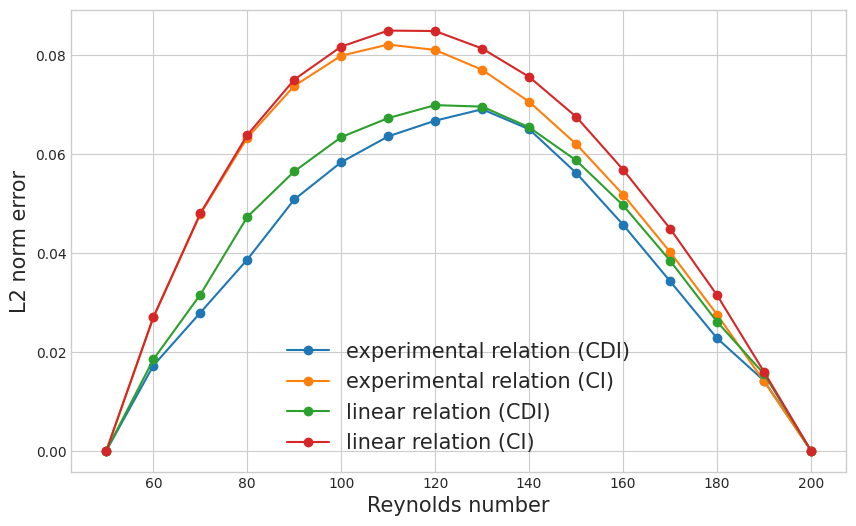}}
 
\caption{incompressible flow past a backward facing step, one-dimensional parameterization.
(a) behavior of the optimal value of $s$ in 
\eqref{eq:relation_sopt} for CDI and CI.
(b) behavior of the relative $L^2(\Omega)$ error for the velocity  field.
}
 \label{fig:BFS_error_2snaps}
 \end{figure}

Figure \ref{fig:BFS_prediction_streamlines}  compares of CDI \eqref{eq:CDI_2field} and CI \eqref{eq:CI_2field} for  
${\rm Re}= 100$.
 We observe that the nonlinear interpolation is qualitatively more accurate to represent the recirculation area.
Figure \ref{fig:BFS_error_2snaps} investigates the behavior of the optimal choice of $s$ (``experimental'') and compares it with the linear relationship --- optimality is measured with respect to the $L^2(\Omega)$ norm.
We observe that $s_{\mu}^{\rm opt}$ depends almost linearly on ${\rm Re}$; we further observe that the improvement due to the nonlinear interpolation is less significant than for the previous case; we again emphasize that the global $L^2$ norm does not properly capture  local-in-space features   of the error.

\subsubsection{Two-dimensional parameterization}
\begin{figure}[H]
\centering

\subfloat[high-fidelity model] 
{\includegraphics[width=0.48\textwidth]
 {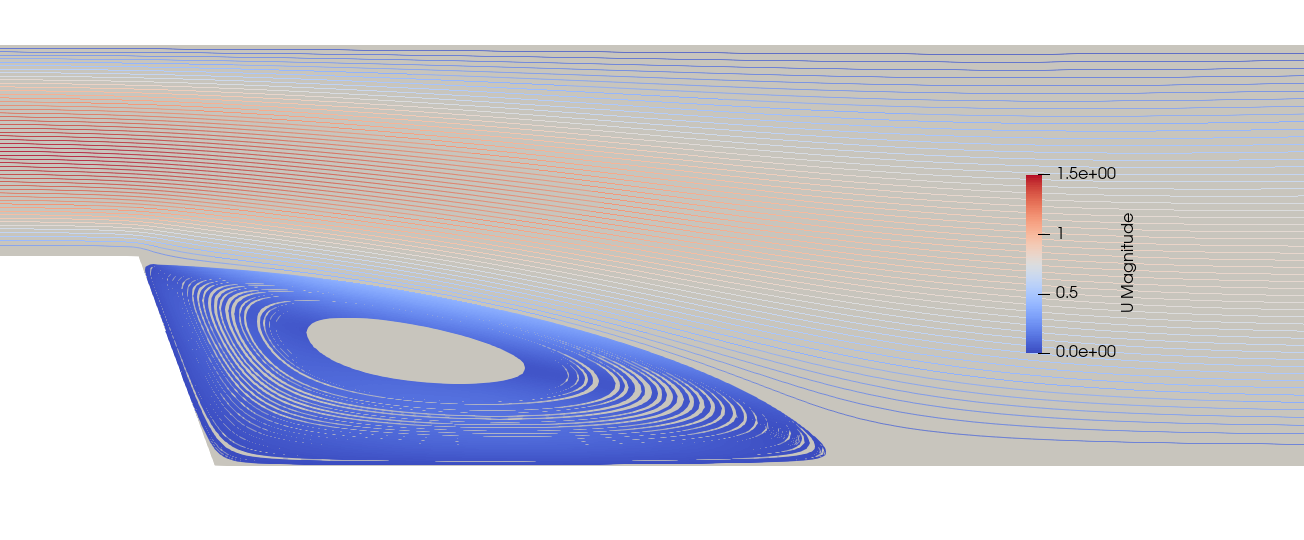}}
 ~~
\subfloat[nonlinear interpolation] 
{\includegraphics[width=0.48\textwidth]
 {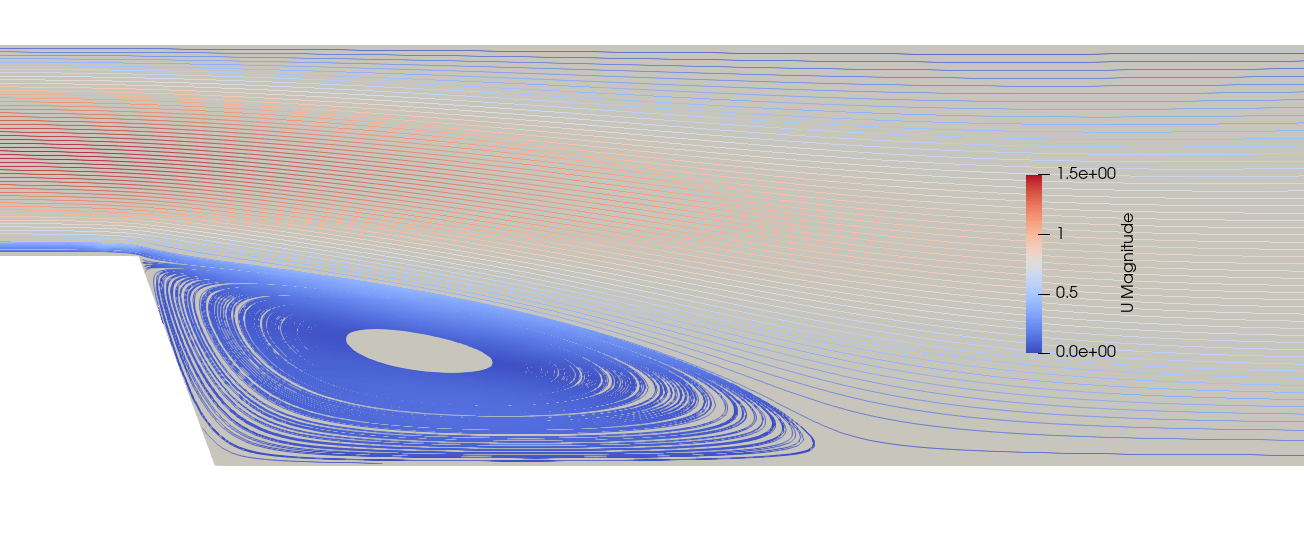}}

\subfloat[linear interpolation] 
{\includegraphics[width=0.48\textwidth]
 {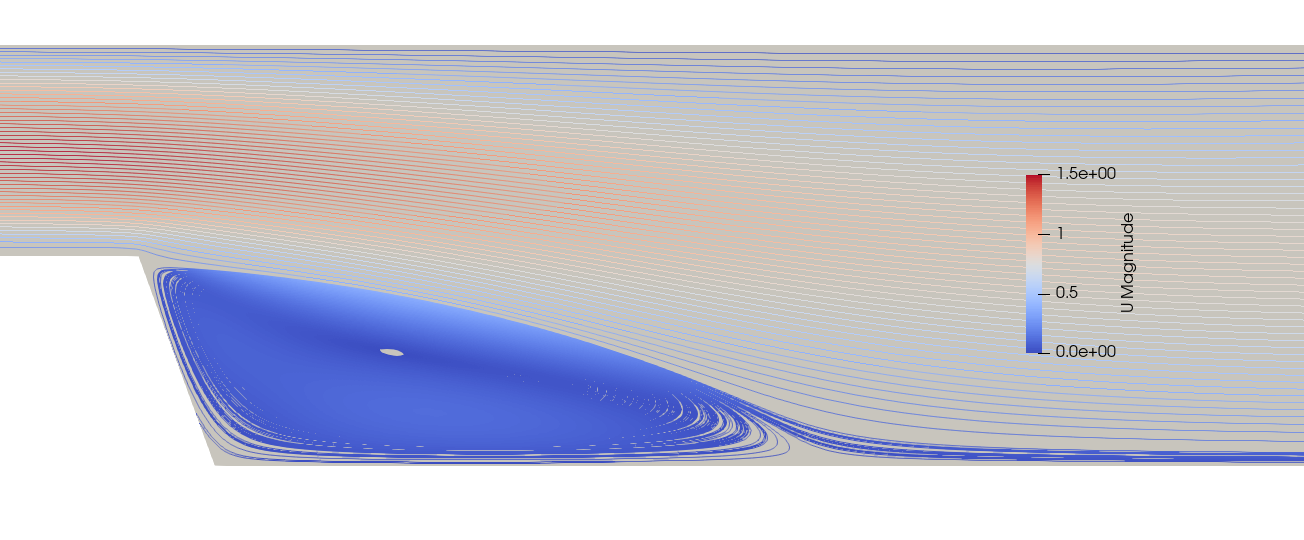}}
  ~~
\subfloat[] 
{\includegraphics[width=0.48\textwidth]
 {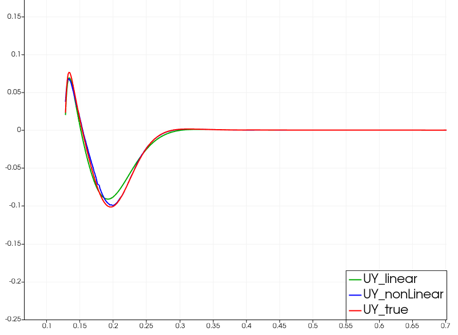}}
 
\caption{incompressible flow past a backward facing step two-dimensional parameterization. Comparison of 
high-fidelity, nonlinear interpolation, and linear interpolation predictions.
(a)-(b)-(c) streamlines of the  velocity field.
(d)  horizontal velocity over the horizontal slice $\{x: x_1\in (0.1,0.7), x_2=-0.01\}$.
}
 \label{fig:BFS_2D_vis}
 \end{figure}

We consider the problem of estimating the solution for 
$\mu=({\rm Re}, \alpha)=[120,70^o]$ based on the training points
$\mu^1=[100,50^o]$, 
$\mu^2=[100,90^o]$ and
$\mu^3=[180,70^o]$.
We  compare performance  of  convex interpolation and CDI in terms of streamlines and horizontal velocity prediction.
As in the previous test case, prediction  of the markers' locations for out-of-sample configurations is performed through linear interpolation. Figure \ref{fig:BFS_2D_vis} shows the results.

We  observe that the improvement due to nonlinear interpolation in terms of global $L^2(\Omega)$ relative error in velocity prediction is very modest ($1.31\%$ vs $1.29\%$).
However, we notice that the streamlines of the flow field that are predicted through  nonlinear interpolation are qualitatively more accurate than the ones associated with   linear interpolation; similarly,  the prediction of the horizontal velocity  based on the nonlinear model over the horizontal slice 
$\{x: x_1\in (0.1,0.7), x_2=-0.01\}$ are 
significantly  more accurate in the recirculation region, especially in the proximity of the minimum of the horizontal  velocity.

\subsection{Compressible flow past a sphere}
\label{sec:sphere}

\subsubsection{Model problem}
\begin{figure}[H]
\centering

\subfloat[${\rm Ma}_{\infty}=2$] 
{\includegraphics[width=0.5\textwidth]
 {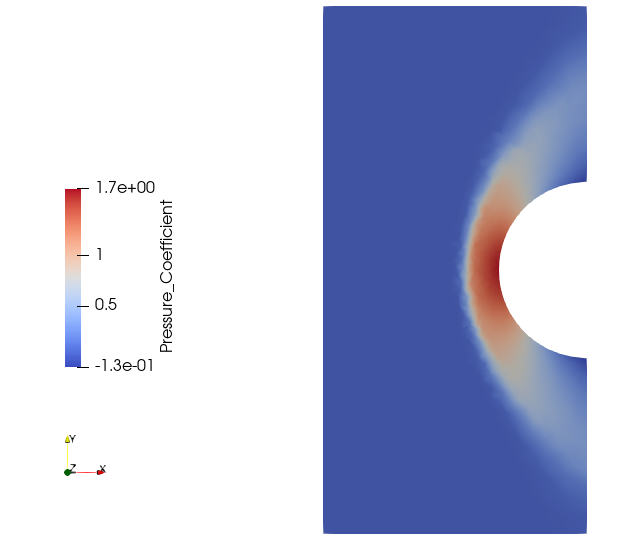}}
 ~~
\subfloat[${\rm Ma}_{\infty}=5$] 
{\includegraphics[width=0.5\textwidth]
 {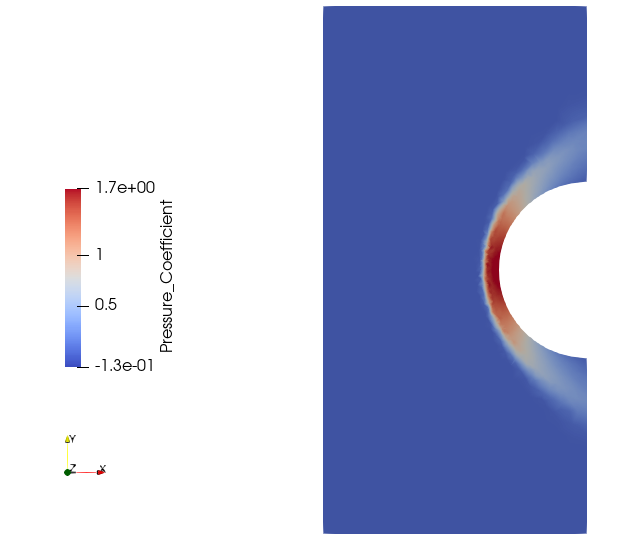}}
 
\caption{compressible flow past a sphere. 
Pressure coefficients for two values of ${\rm Ma}_{\infty}$.
}
 \label{fig:sphere_vis}
 \end{figure}

\begin{figure}[H]
\centering

\subfloat[${\rm Ma}_{\infty}=2$] 
{\includegraphics[width=0.5\textwidth]
 {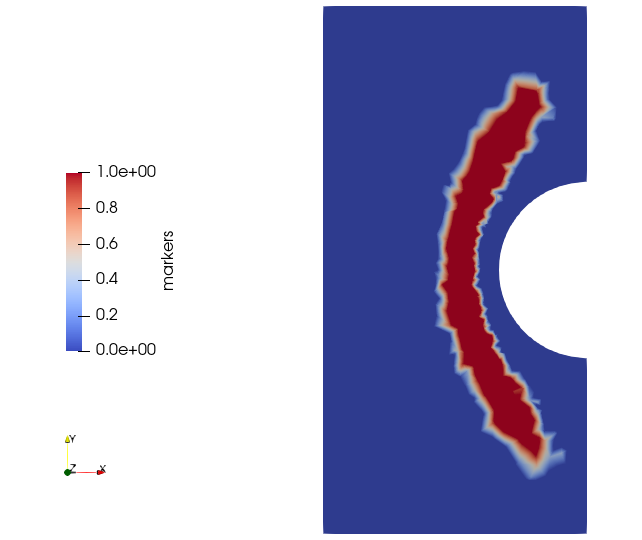}}
 ~~
\subfloat[${\rm Ma}_{\infty}=5$] 
{\includegraphics[width=0.5\textwidth]
 {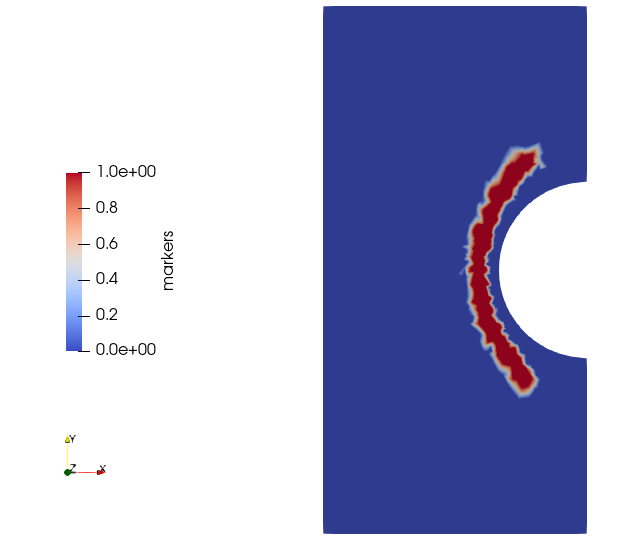}}
 
\caption{compressible flow past a sphere. 
Marked cells for ${\rm Ma}_{\infty}=2$ and ${\rm Ma}_{\infty}=5$.
}
 \label{fig:sphere_markers}
 \end{figure}

We study the compressible supersonic inviscid flow past a sphere for inflow Mach number ${\rm Ma}_{\infty}\in [2,5]$. We model the flow using the compressible Euler equations for ideal gases and we rely on SU2 for snapshot generation. As in section \ref{sec:RAE}, we apply the nonlinear interpolation strategy to estimate the pressure coefficient \eqref{eq:pressure_coefficient} and we rely on the Ducros sensor with $\gamma_{\rm thr}= 99.6\%$ to identify the coherent structure.
Figure \ref{fig:sphere_vis} shows the behavior of $C_{\rm p}$ for ${\rm Ma}_{\infty}=2$ and ${\rm Ma}_{\infty}=5$: we observe that the solution exhibits a bow shock ahead of the sphere whose location depends on the parameter. 
Figure \ref{fig:sphere_markers} shows the marked cells for the same parameter values.

\begin{figure}[H]
\centering

\subfloat[] 
{\includegraphics[width=0.5\textwidth]
 {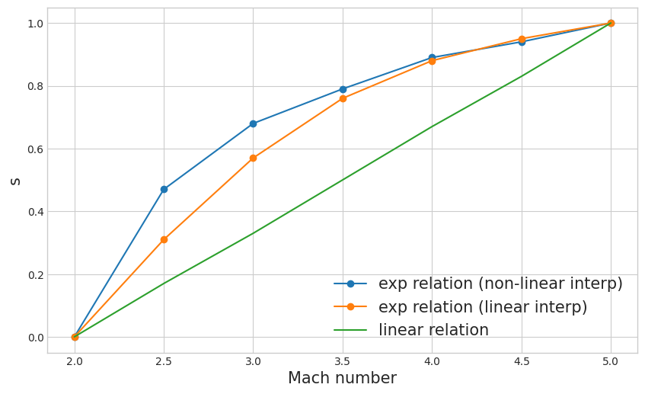}}
 ~~
\subfloat[] 
{\includegraphics[width=0.5\textwidth]
 {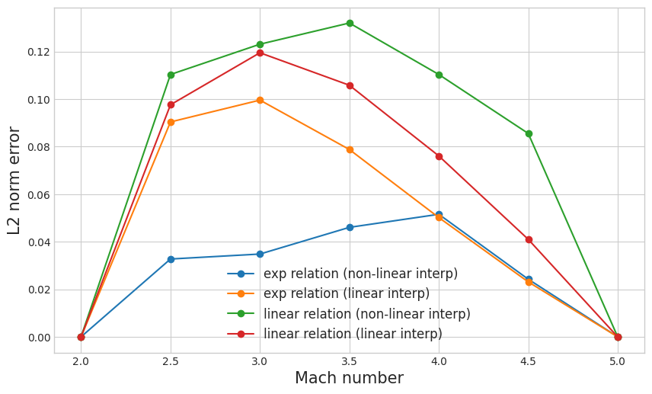}}
 
\caption{compressible flow past a sphere. 
(a) behavior of the optimal value of $s$ in 
\eqref{eq:relation_sopt} for CDI and CI.
(b) behavior of the relative $L^2(\Omega)$ error for the pressure coefficient.
}
 \label{fig:sphere_error_analysis}
 \end{figure}

Figure \ref{fig:sphere_error_analysis} illustrates the performance   of CDI \eqref{eq:CDI_2field} and CI \eqref{eq:CI_2field} for ${\rm Ma}_{\infty}\in [2,5]$:
Figure 
\ref{fig:sphere_error_analysis}(a) investigates the  behavior of the optimal value of $s$ in 
\eqref{eq:relation_sopt} for CDI and CI;
Figure 
\ref{fig:sphere_error_analysis}(b) shows the behavior of the relative $L^2(\Omega)$ error for the pressure coefficient.
Interestingly, we observe that the behavior of $s_{\mu}^{\rm opt}$ is strongly nonlinear with respect to the Mach number. As for the previous cases, nonlinear interpolation  is superior to linear interpolation, especially   when the parameter $s$ in \eqref{eq:CDI_2field} is chosen in an optimal manner.

\section{Data augmentation}
\label{sec:data_augmentation_numerics}

\subsection{Model problem: inviscid flow past a bump}
\begin{figure}[H]
\centering

\subfloat[$\mu=(0.05,0.78)$] 
{\includegraphics[width=0.5\textwidth]
 {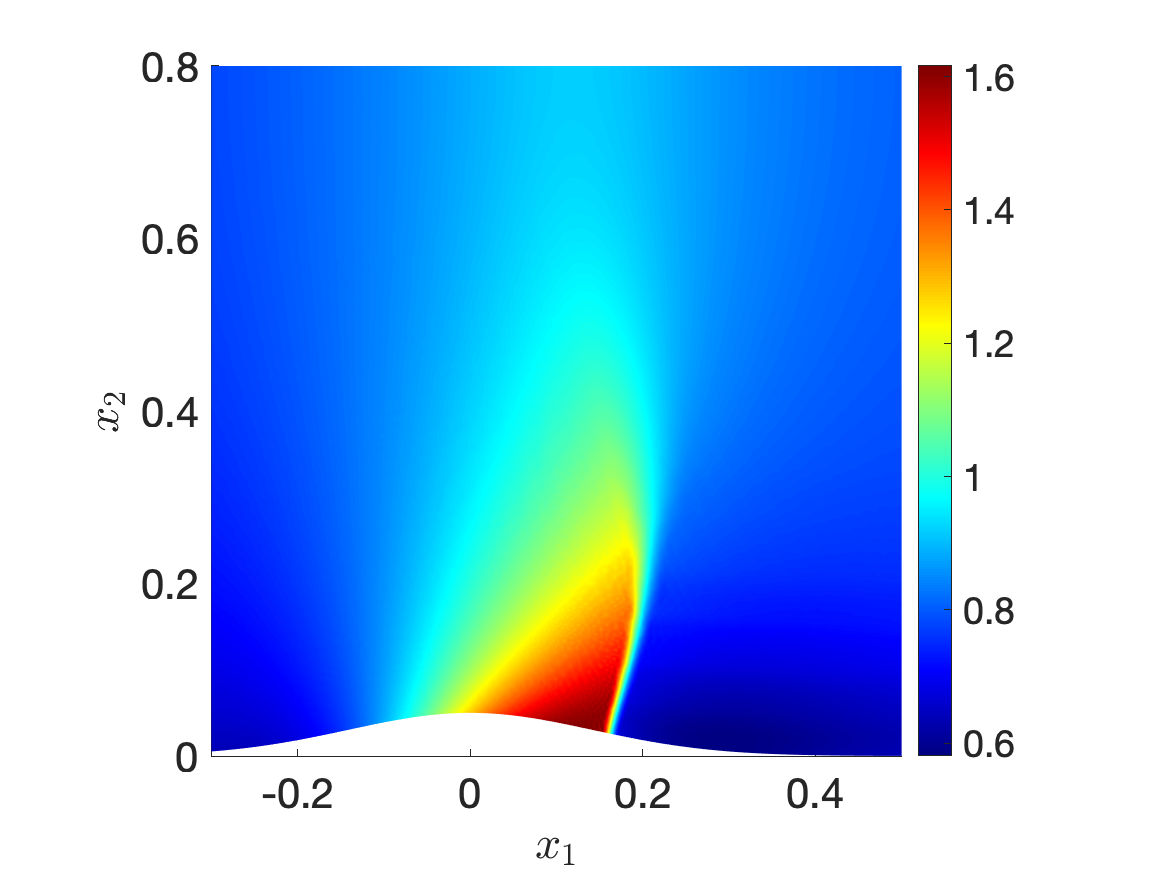}}
 ~~
   \subfloat[$\mu=(0.065,0.78)$] 
{\includegraphics[width=0.5\textwidth]
 {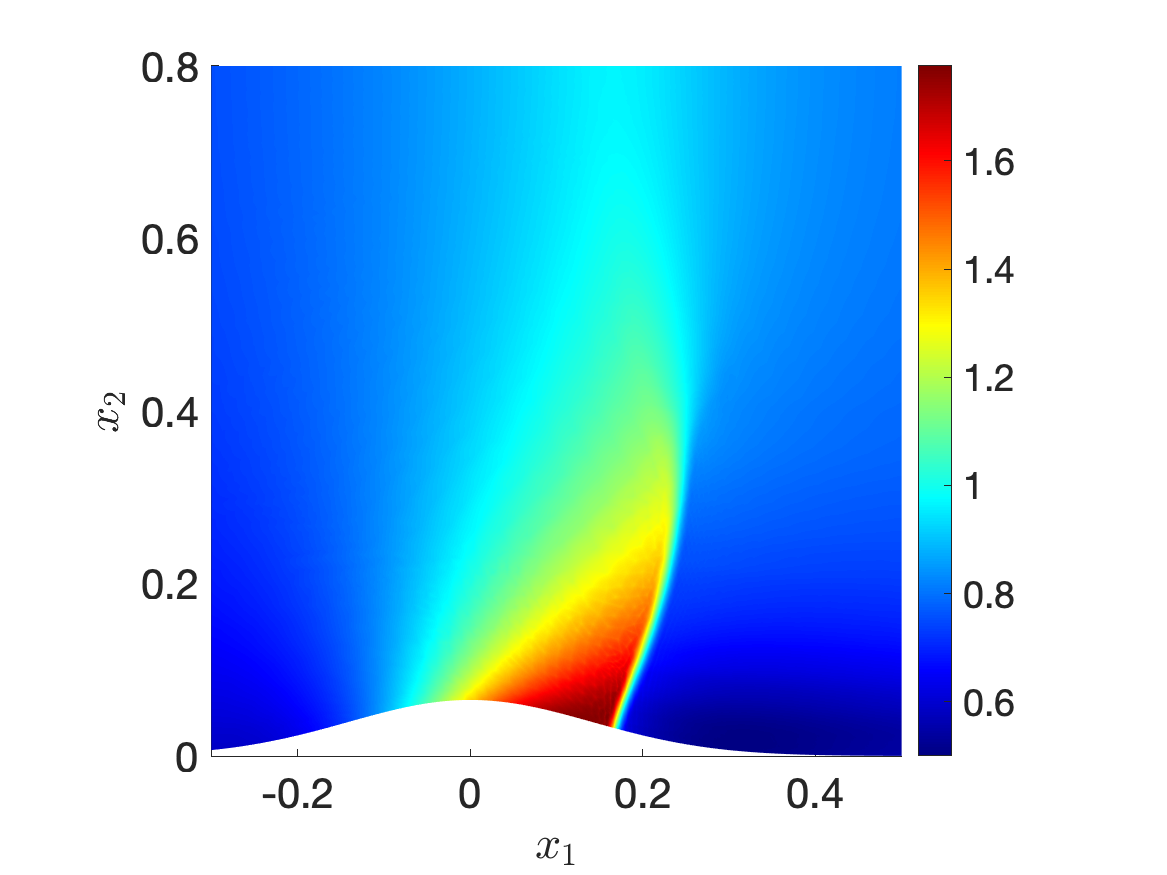}}
 
 \caption{inviscid flow over a Gaussian bump. 
Visualization of the Mach field for  two parameter values.
}
 \label{fig:transbump_vis}
 \end{figure}

We consider the problem of estimating the solution to the two-dimensional Euler equations in the transonic regime in a channel with a Gaussian bump. 
We denote by $\rho$ the density, by $v$ the velocity field, by $E$ the total energy, and by $p=
 (\gamma - 1) \left( E  - \frac{1}{2} \rho \| v  \|_2^2\right)$ 
 the pressure, where $\gamma>0$ is the ratio of specific heats;
 we further denote  by 
$u=[\rho,\rho v, E]$ 
the vector of conserved variables.
We introduce the computational domain $\Omega_{\mu} = \{
x\in (-1.5,1.5)\times (0,0.8): x_2 > h e^{-25 x_1^2}
\}$ where $h>0$ is a given parameter.
We consider the conservation law:
\begin{equation}
\label{eq:transbump_flow}
\nabla \cdot F(u^{\rm true}) = 0,
\;  {\rm where} 
\;\;
F(u) = \left[
\begin{array}{l}
\rho v^{\top}  \\
 \rho v v^{\top}  \\
(E+p) v^{\top}  \\
\end{array}
\right],
\end{equation}
completed with wall boundary conditions on top and bottom boundaries, subsonic inlet condition (total temperature, total pressure and flow direction) at the left boundary and  subsonic outlet condition (static pressure) at the right boundary. We express the free-stream field $u_{\infty}$ in terms of the Mach number ${\rm Ma}_{\infty}$,
$$
T_{\infty}=1, \;\;
p_{\infty}=\frac{1}{\gamma}, \;\;
\rho_{\infty}= 1, \;\;
u_{\infty}= [{\rm Ma}_{\infty},0],
$$
where $T=\frac{p}{R \rho}$ is the temperature and $R=\gamma-1$ is the (non-dimensional) specific gas constant. 
Finally, we introduce the parameter vector $\mu=[h, {\rm Ma}_{\infty}]$ and the parameter region $\mathcal{P}=[0.05,0.065] \times [0.68,0.78]$.
Figures \ref{fig:transbump_vis}(a) and (b) show the behavior of the Mach field for two values of the parameters in $\mathcal{P}$.

\subsection{Geometric parameterization and high-fidelity solver}

The computational  domain $\Omega_{\mu}$ depends on the geometric parameter $h$; therefore,  we should introduce a geometric mapping to recast the problem over a parameter-independent configuration. 
We here resort to a Gordon-Hall transformation: 
we denote by $\Phi_{\rm geo,\mu}$ the geometric mapping from the rectangle $\widetilde{\Omega}=(-1.5,1.5) \times (0,0.8)$ to the physical domain $\Omega_{\mu}$ 
(see, e.g.,  \cite[section 2]{ferrero2022registration}).
We further define the 
domain $\Omega:=\Omega_{\bar{\mu}}$, where $\bar{\mu}$ is the centroid of $\mathcal{P}$ and   the additional geometric 
mapping
$\Psi_{\rm geo,\mu} := \Phi_{\rm geo,\mu}
\circ \Phi_{\rm geo,\bar{\mu}}^{-1}:\Omega\to \Omega_{\mu}$. 

{The reference configuration $\widetilde{\Omega}$ is used below to define the sensors and the regression procedure; on the other hand, the HF mesh is defined in the domain $\Omega$. The latter choice is justified by the observation that for small variations of $h$ the mapping $\Psi_{\rm geo}$ is (significantly) closer to the identity map than 
$\Phi_{\rm geo}$;  therefore, defining the mesh in $\Omega$ offers significantly enhanced control over the quality of the deformed mesh.
}

Simulations are
  performed in Matlab 2022a
\cite{MATLAB:2022} based on an in-house code, and executed over a commodity Linux  workstation (RAM 32 GB, Intel i7 CPU 3.20 GHz x 12).
We consider a P2 mesh with $N_{\rm e}=12748$ elements; given a new value of the parameter $\mu$, we define the deformed mesh using the geometric mapping $\Psi_{\rm geo}$,
and then we call the DG solver to compute the solution field.
{We measure performance and we compute reduced spaces in the reference domain $\Omega$; we denote by $(\cdot, \cdot)$ the $[L^2(\Omega)]^4$ inner product and by $\| \cdot \| = \sqrt{(\cdot, \cdot)}$ the induced norm.
}

\subsection{Model order reduction}
We rely on a 
projection-based ROM based on
a least-square Petrov-Galerkin (LSPG, \cite{carlberg2013gnat}) formulation and 
empirical quadrature 
(EQ, \cite{farhat2015structure,yano2019lp}) hyper-reduction; we refer to \cite{ferrero2022registration} for a detailed presentation of the method.

\subsection{Convex displacement interpolation}
We apply Algorithm \ref{alg:CDI_overview} to construct the CDI estimate.

\noindent
\emph{Feature selection.}
We rely on the Ducros sensor \eqref{eq:ducros_sensor} to identify the elements of the mesh that contain the shock. 
We set $\gamma_{\rm thr}=99\%$ in \eqref{eq:quantile}   to determine the threshold $\mathcal{t}_{\mu}$.
In order to take into account geometric variations, we define the ``mapped'' point clouds
$\widetilde{X}_{\mu}^{\rm raw} = \{  \Phi_{\rm geo,\mu}^{-1}(   x_{i,\mu}^{\rm raw}   )  \}_{i=1}^{N_{\mu}} \subset \widetilde{\Omega}$,
where 
$\Phi_{\rm geo,\mu}: \widetilde{\Omega} \to \Omega_{\mu}$ is the Gordon-Hall map.
Then, we set 
 $\widetilde{X}^{\rm ref} = \widetilde{X}_{\mu^{\rm ref} }^{\rm raw}$ where $\mu^{\rm ref}$ is the parameter in $\mathcal{P}_{\rm train}$ that minimizes the distance from the centroid of $\mathcal{P}$, and we rely on  Gaussian-based PSR
to determine the sorted point clouds
$\{\widetilde{X}_{\mu}: \mu \in \mathcal{P}_{\rm train} \}$.

\noindent 
\emph{Regression.}
We apply the regression procedure described 
in section \ref{sec:regression} to find the RBF model
$\mu \in \mathcal{P} \mapsto  \widehat{X}_{\mu}^{\star} = 
\{ \widehat{x}_{i,\mu}^{\star} \}_{i=1}^N     \subset \widetilde{\Omega}$.
Then, we define the point cloud in physical domain
using the geometric mapping $\Phi_{\rm geo}$: 
$\widehat{X}_{\mu} = 
\{ \widehat{x}_{i,\mu}   = \Phi_{\rm geo,\mu}(\widehat{x}_{i,\mu}^{\star} )\}_{i=1}^N$.
The choice of applying RBF and PSR in the reference (parameter-independent) domain $\widetilde{\Omega}$ robustifies the learning process and ultimately improves prediction performance.

\noindent 
\emph{Nearest neighbors, registration, and choice of the weights.}
We build the CDI based on $\kappa=4$ neighbors; we consider the Euclidean distance to identify the neighbors and also in the IDW strategy for the weights (cf. section \ref{sec:weights}).
Finally, we resort to   optimization-based registration   (cf. \eqref{eq:optimization_based_registration}) to determine the mapping $\Phi$: we consider the 
 penalty term $\mathfrak{f}_{\rm jac}$ 
 (\cite[section 4.3]{taddei2023compositional}) 
 that penalizes small values of the Jacobian determinant, 
and we rely on a piecewise-polynomial search space with $M=198$ degrees of freedom.

\subsection{Data augmentation procedure}

We denote by  $\mathcal{P}_{\rm train} = \{  \mu^i \}_{i=1}^{n_{\rm train}} \subset \mathcal{P}$ the training set of parameters for which HF simulations are available and by
 $\mathcal{P}_{\rm train,cdi} = \{  \nu^i \}_{i=1}^{n_{\rm train,cdi}} \subset \mathcal{P}$ an additional set of parameters.
 We first define the Lagrangian space
 $\mathcal{Z}_0 = {\rm span} \{ u_{ \mu^i}^{\rm hf}  \}_{i=1}^{n_{\rm train}}$
and the projection operator
 $\Pi_{  \mathcal{Z}_0  }: [L^2(\Omega)]^4 \to \mathcal{Z}_0$; 
 then, we define the POD space
 \begin{equation}
 \label{eq:POD_data_augmentation}
\mathcal{Z}_n =
 \mathcal{Z}_0 \oplus
 \mathcal{Z}_{\rm lf},\quad
  \mathcal{Z}_{\rm lf}
  =\texttt{POD} \left(
  \mathcal{D}_{\rm train,lf}, n-n_{\rm train}, (\cdot,\cdot)
\right),
 \end{equation}
where $\oplus$ indicates the direct sum of orthogonal spaces,
$\mathcal{D}_{\rm train,lf} = 
\{ \widehat{u}_{\mu} - \Pi_{  \mathcal{Z}_0  }
\widehat{u}_{\mu}: \mu \in \mathcal{P}_{\rm train,cdi} \}$ and  
$  \mathcal{Z}_{\rm lf}$ is the POD space of size 
$ n-n_{\rm train}$ associated with the snapshots 
$\mathcal{D}_{\rm train,lf}$ and the $L^2(\Omega)$ inner product $(\cdot,\cdot)$.
We observe that the space $\mathcal{Z}_n$ contains 
the HF training set: potential inaccuracies of the CDI estimates do not pollute the performance of the reduced space for $\mu\in \mathcal{P}_{\rm train}$.

\subsection{Numerical results: nonlinear interpolation}

\begin{figure}[H]
\centering
\includegraphics[width=0.6\textwidth]
 {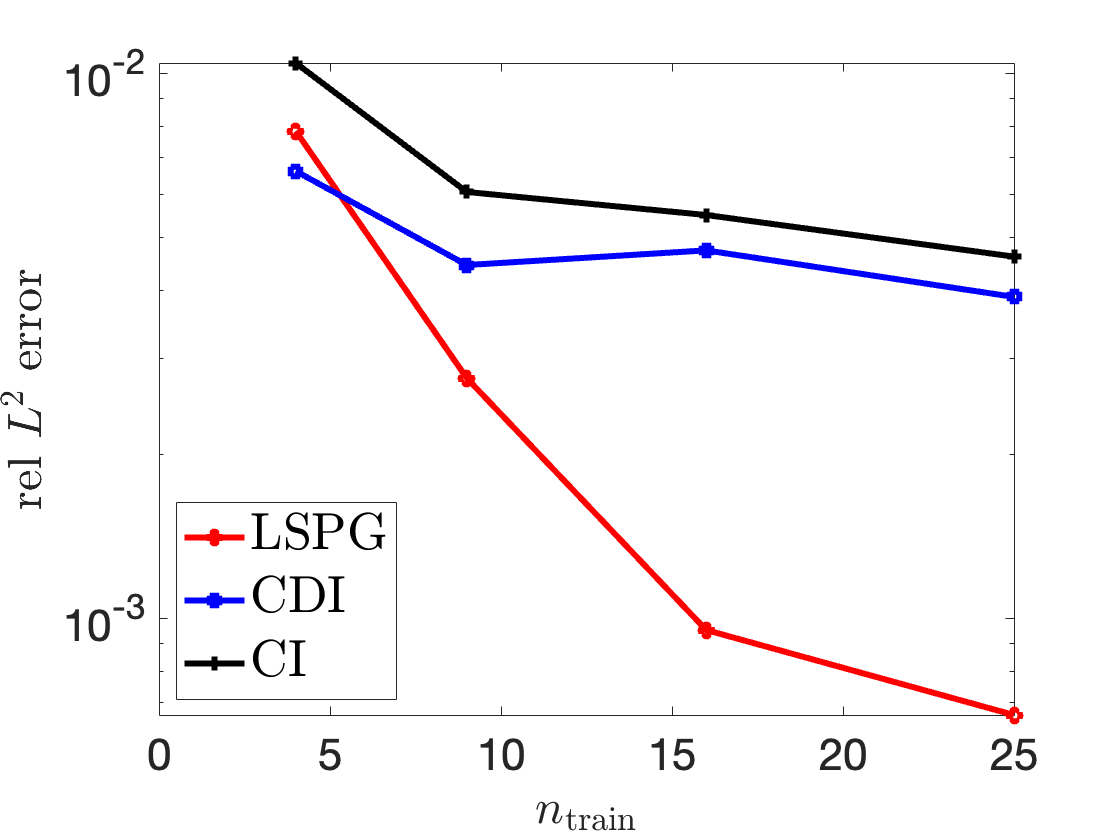}
  
 \caption{nonlinear interpolation. Behavior of the average relative $L^2$ error for LSPG, CDI and CI for snapshot sets of increasing size $n_{\rm train}$.
}
 \label{fig:transbump_onemore}
 \end{figure}

We first investigate the performance of CDI for this model problem. Towards this end, we study the behavior of the average relative $L^2(\Omega)$ state prediction error of
the projection-based linear ROM (LSPG),
CDI, 
and convex interpolation (CI), for tensorized grids of increasing size
$n_{\rm train}$; in all our experiments we set $\kappa=4$ for CDI and CI, while the reduced space of  LSPG includes all the training snapshots ($n=n_{\rm train}$).
Figure \ref{fig:transbump_onemore} shows the results.

We observe that CDI outperforms CI for all parameter values and all training set sizes  $n_{\rm train}$. 
 Interestingly, CDI is more accurate than LSPG for $n_{\rm train}=4$: 
this can be explained by observing that, unlike the LSPG ROM, CDI provides estimates that do not belong to the span of the dataset of HF simulations; 
 for small values of  $n_{\rm train}$, 
CDI can hence  compensate the fact that 
 it does not exploit the knowledge of the PDE model for prediction.

\subsection{Numerical results: data augmentation}
\begin{figure}[H]
\centering
 \subfloat[] 
{  \includegraphics[width=0.5\textwidth]
 {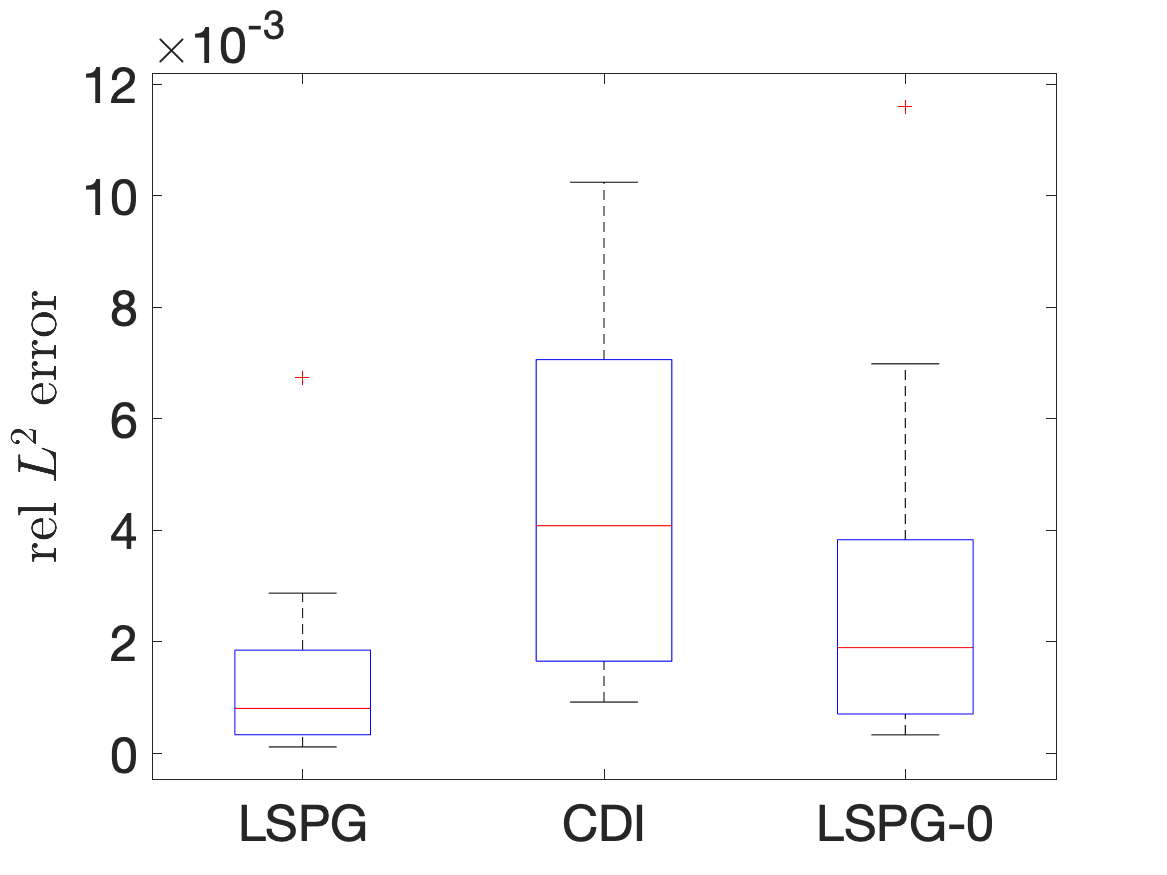}}
   ~~
 \subfloat[] 
{  \includegraphics[width=0.5\textwidth]
 {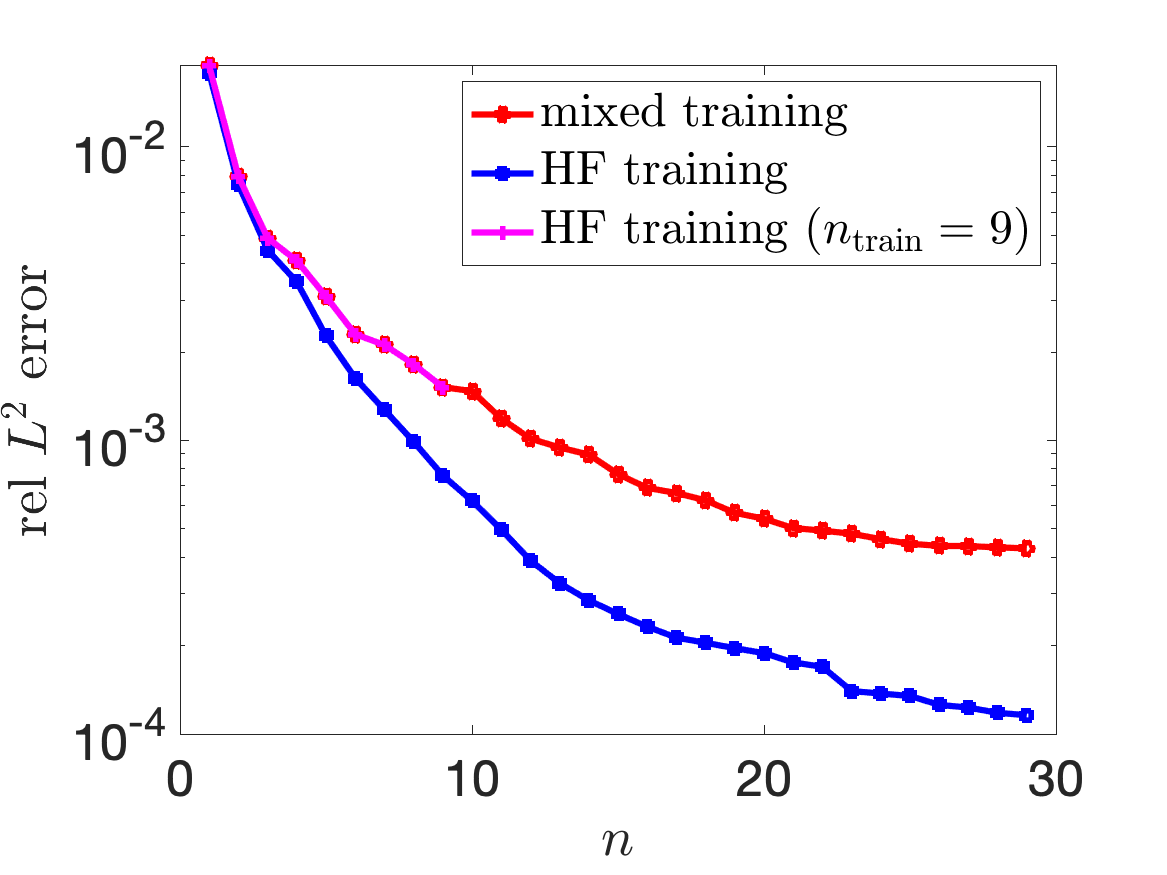}}
 
 \caption{data augmentation $n_{\rm train}=9$, $n_{\rm train,cdi}=20$. 
Behavior of LSPG ROM with CDI-based data augmentation.
(a) boxplots of relative $L^2$ error for LSPG ROM, CDI and projection error ($n=29$).
(b) projection error.
}
 \label{fig:transbump_ntrain9}
 \end{figure}

Figure \ref{fig:transbump_ntrain9} illustrates the performance for $n_{\rm train}=9$, $n_{\rm train,cdi}=20$.  We consider  a tensorized $3\times 3$ grid for $\mathcal{P}_{\rm train}$ and we generate $\mathcal{P}_{\rm train,cdi}$  using latin hyper-cube sampling \cite{mckay2000comparison}.
We consider a  POD  space in \eqref{eq:POD_data_augmentation} of dimension $n=29$;
as for the previous experiment,
we assess performance based on a dataset of $n_{\rm test}=20$ out-of-sample parameters.
Figure \ref{fig:transbump_ntrain9}(a)
shows boxplots of the relative $L^2(\Omega)$ error 
for LSPG ROM, CDI and LSPG ROM based on $n=9$ modes (i.e., without data augmentation) (LSPG-0).
We observe that CDI is not as accurate as LSPG-0;
however, it can be used to augment the dataset of simulations and ultimately contribute to reduce  the prediction error.
Figure \ref{fig:transbump_ntrain9}(b) shows the projection error over the test for the POD space built using \eqref{eq:POD_data_augmentation} (``mixed training''), the POD space based on 
 $n_{\rm train} + n_{\rm train,cdi}$ HF simulations 
(``HF training'') and 
the POD space based 
 $n_{\rm train}$ HF simulations 
(``HF training ($n_{\rm train}=9$)''):
we observe that the mixed training strategy based on CDI is not as accurate as  HF training, particularly for large values of $n$; however, it enables significant reduction of the projection error without resorting to new simulations.

\begin{figure}[H]
\centering
 \subfloat[] 
{  \includegraphics[width=0.45\textwidth]
 {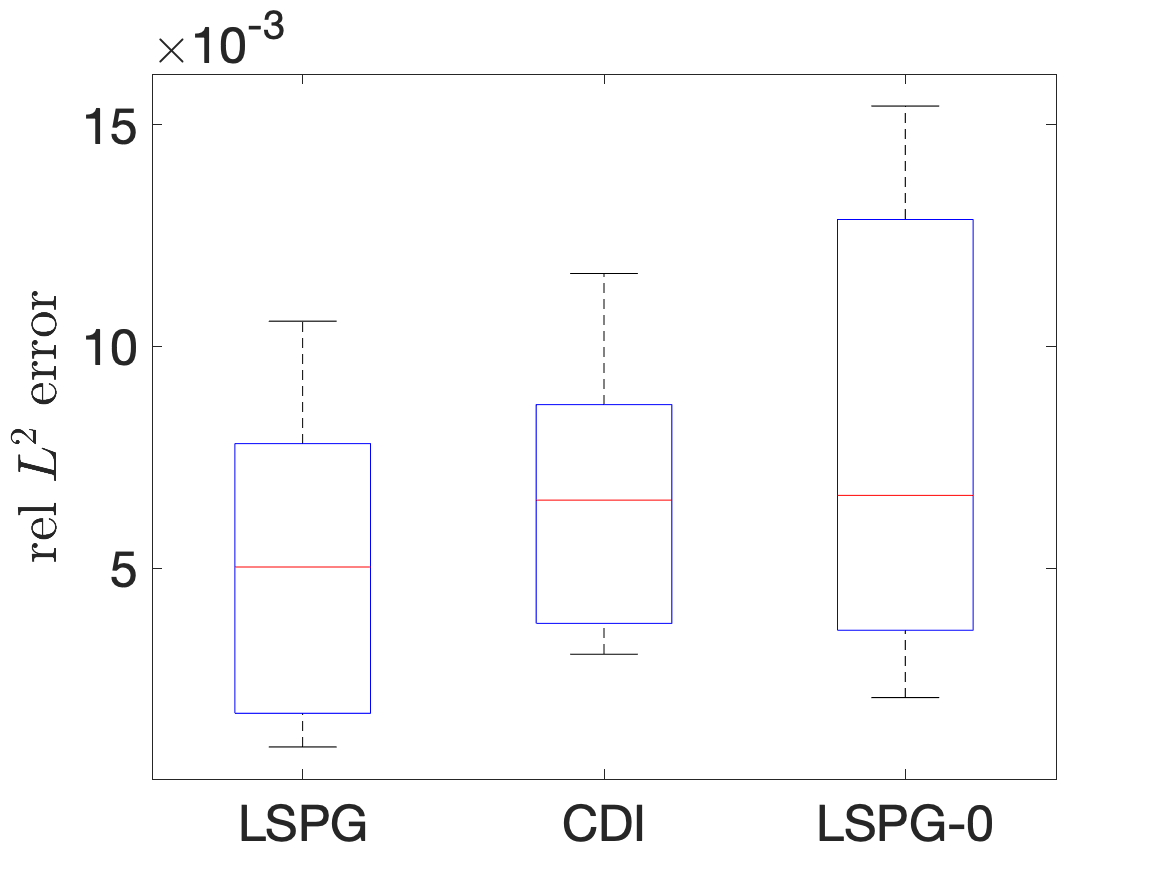}}
   ~~
 \subfloat[] 
{  \includegraphics[width=0.45\textwidth]
 {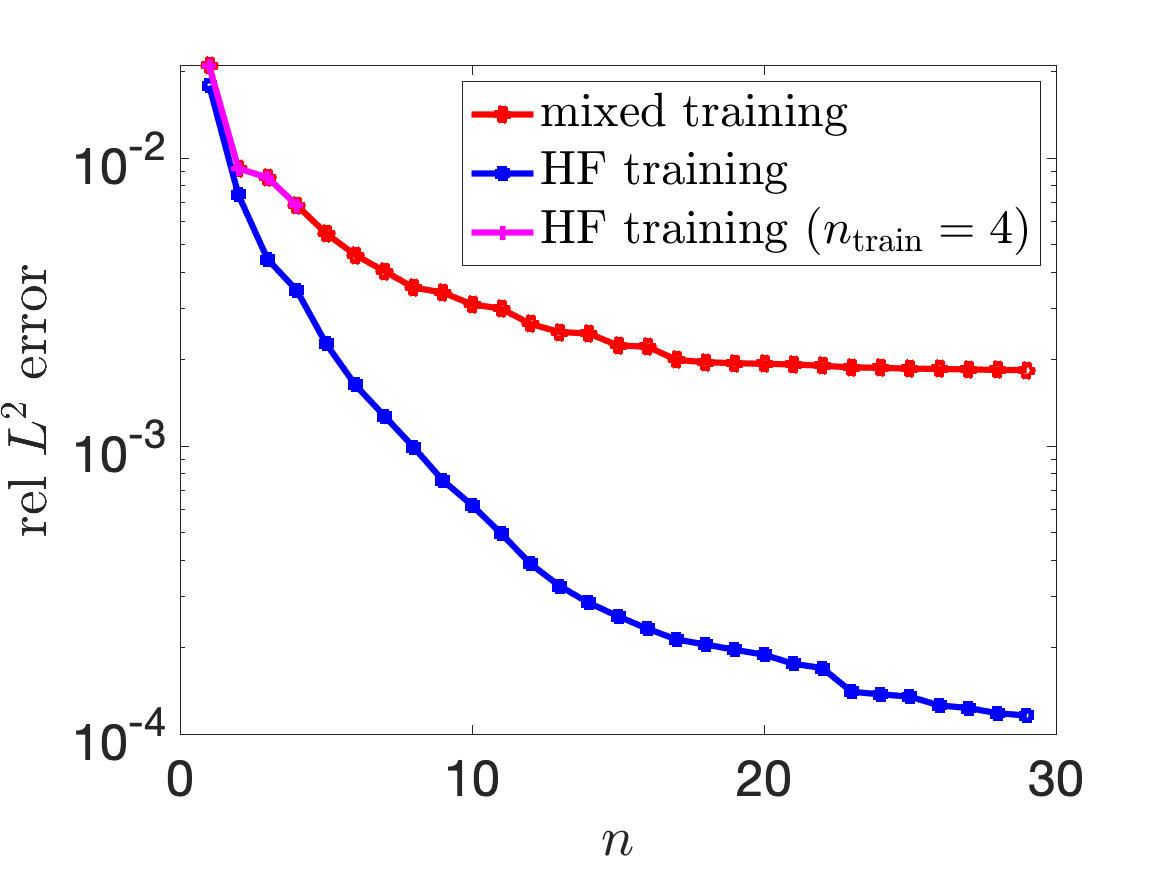}}
 
 \caption{data augmentation $n_{\rm train}=4$, $n_{\rm train,cdi}=25$. 
Behavior of LSPG ROM with CDI-based data augmentation.
(a) boxplots of relative $L^2$ error for LSPG ROM, CDI and projection error ($n=29$).
(b) projection error.
}
 \label{fig:transbump_ntrain4}
 \end{figure}

Figure \ref{fig:transbump_ntrain4} 
replicates the results of Figure \ref{fig:transbump_ntrain9} for  
 $n_{\rm train}=4$, $n_{\rm train,cdi}=25$.  
 Note that CDI is more accurate than LSPG-0, which is consistent with the results of Figure
 \ref{fig:transbump_onemore}.

Table  \ref{tab:comp_costs_transbump} provides an overview of the computational costs. 
The HF solver is roughly $14$ times slower than CDI and $625$ times slower than the LSPG ROM; the construction of the LSPG ROM,  which comprises
the construction of the trial and test spaces and hyper-reduction, is roughly as expensive as a single full-order solve.
We remark that approximately  $99\%$ of the total cost of CDI is associated with   mesh interpolation: we hence  envision that more sophisticated implementations of this step 
(see, e.g., \cite{luke2012fast})
might lead to a significant reduction of the overall CDI cost.

\begin{table}[H]
\centering
\begin{tabular}{||l| l||} 
\hline 
& cost [s] \\[0.5ex]  \hline \hline 
HF solver (avg) & 
$481.21$
  \\[0.5ex] \hline 
ROM generation (data compression+hyperreduction) &   
$431.24$
\\[0.5ex] \hline 
Convex displacement interpolation (avg)  &
$33.49$
 \\[0.5ex]   \hline 
LSPG ROM  ($n=29$, avg) &   
$0.77$ \\[0.5ex]  \hline 
\end{tabular}
\caption{data augmentation. Wall-clock costs of HF and MOR procedures ($n_{\rm train}=9, n_{\rm train,cdi}=20$).}
\label{tab:comp_costs_transbump}
\end{table}

\section{Conclusions}
\label{sec:conclusions}

We presented a general (i.e., independent of the underlying equations) nonlinear interpolation procedure for parametric fields with local compactly-supported coherent structures; the approach is dubbed convex displacement interpolation (CDI) and extends the previous proposal of   \cite{iollo2022mapping}.
The methodology comprises five main building blocks that were all discussed in section \ref{sec:methods}: 
sensor selection, point cloud matching, regression of sorted point clouds, boundary-aware registration, and piecewise-multivariate interpolation. 
CDI assumes that coherent structures of the solution field --- such as shocks, shear layers, vortices --- vary smoothly with the parameter. The numerical experiments suggest that this assumption holds for a broad and diverse range of parametric problems in continuum mechanics.

We illustrated the application of CDI to data augmentation.
Since CDI generates estimates that do not belong to the span of the training snapshots, we can exploit   predictions for out-of-sample parameters  to augment the dataset  of simulations that are later used for the construction of projection-based ROMs.
We envision that combination of lower and higher fidelity snapshots for MOR training has the potential to significantly reduce offline costs, and ultimately enable the application of MOR techniques to complex engineering problems for which extensive parameter explorations based on HF solves are not practical.
 
\section*{Acknowledgements}
The authors acknowledge the support by European Union’s Horizon 2020 research and innovation program  under the Marie Skłodowska-Curie Actions, grant agreement 872442 (ARIA).

\bibliographystyle{abbrv}	
\bibliography{all_refs}

\begin{thebibliography}{10}

\bibitem{amsallem2008interpolation}
D.~Amsallem and C.~Farhat.
\newblock Interpolation method for adapting reduced-order models and
  application to aeroelasticity.
\newblock {\em AIAA journal}, 46(7):1803--1813, 2008.

\bibitem{amsallem2012nonlinear}
D.~Amsallem, M.~J. Zahr, and C.~Farhat.
\newblock Nonlinear model order reduction based on local reduced-order bases.
\newblock {\em International Journal for Numerical Methods in Engineering},
  92(10):891--916, 2012.

\bibitem{angot1999penalization}
P.~Angot, C.-H. Bruneau, and P.~Fabrie.
\newblock A penalization method to take into account obstacles in
  incompressible viscous flows.
\newblock {\em Numerische Mathematik}, 81(4):497--520, 1999.

\bibitem{barnett2022quadratic}
J.~Barnett and C.~Farhat.
\newblock Quadratic approximation manifold for mitigating the {K}olmogorov
  barrier in nonlinear projection-based model order reduction.
\newblock {\em Journal of Computational Physics}, 464:111348, 2022.

\bibitem{barnett2022neural}
J.~L. Barnett, C.~Farhat, and Y.~Maday.
\newblock Neural-network-augmented projection-based model order reduction for
  mitigating the kolmogorov barrier to reducibility of {CFD} models.
\newblock {\em arXiv preprint arXiv:2212.08939}, 2022.

\bibitem{bernard2018reduced}
F.~Bernard, A.~Iollo, and S.~Riffaud.
\newblock Reduced-order model for the {BGK} equation based on {POD} and optimal
  transport.
\newblock {\em Journal of Computational Physics}, 373:545--570, 2018.

\bibitem{camion2001geodesic}
V.~Camion and L.~Younes.
\newblock Geodesic interpolating splines.
\newblock In {\em International workshop on energy minimization methods in
  computer vision and pattern recognition}, pages 513--527. Springer, 2001.

\bibitem{cao2005large}
Y.~Cao, M.~I. Miller, R.~L. Winslow, and L.~Younes.
\newblock Large deformation diffeomorphic metric mapping of vector fields.
\newblock {\em IEEE transactions on medical imaging}, 24(9):1216--1230, 2005.

\bibitem{carlberg2013gnat}
K.~Carlberg, C.~Farhat, J.~Cortial, and D.~Amsallem.
\newblock The {GNAT} method for nonlinear model reduction: effective
  implementation and application to computational fluid dynamics and turbulent
  flows.
\newblock {\em Journal of Computational Physics}, 242:623--647, 2013.

\bibitem{chantalat2009level}
F.~Chantalat, C.-H. Bruneau, C.~Galusinski, and A.~Iollo.
\newblock Level-set, penalization and cartesian meshes: A paradigm for inverse
  problems and optimal design.
\newblock {\em Journal of Computational Physics}, 228(17):6291--6315, 2009.

\bibitem{ching2022model}
D.~S. Ching, P.~J. Blonigan, F.~Rizzi, and J.~A. Fike.
\newblock Model reduction of hypersonic aerodynamics with residual minimization
  techniques.
\newblock In {\em AIAA SCITECH 2022 Forum}, page 1247, 2022.

\bibitem{economon2016su2}
T.~D. Economon, F.~Palacios, S.~R. Copeland, T.~W. Lukaczyk, and J.~J. Alonso.
\newblock {SU2}: An open-source suite for multiphysics simulation and design.
\newblock {\em Aiaa Journal}, 54(3):828--846, 2016.

\bibitem{eftang2010hp}
J.~L. Eftang, A.~T. Patera, and E.~M. R{\o}nquist.
\newblock An ``hp'' certified reduced basis method for parametrized elliptic
  partial differential equations.
\newblock {\em SIAM Journal on Scientific Computing}, 32(6):3170--3200, 2010.

\bibitem{farhat2015structure}
C.~Farhat, T.~Chapman, and P.~Avery.
\newblock Structure-preserving, stability, and accuracy properties of the
  energy-conserving sampling and weighting method for the hyper reduction of
  nonlinear finite element dynamic models.
\newblock {\em International journal for numerical methods in engineering},
  102(5):1077--1110, 2015.

\bibitem{ferrero2022registration}
A.~Ferrero, T.~Taddei, and L.~Zhang.
\newblock Registration-based model reduction of parameterized two-dimensional
  conservation laws.
\newblock {\em Journal of Computational Physics}, 457:111068, 2022.

\bibitem{guo2018reduced}
M.~Guo and J.~S. Hesthaven.
\newblock Reduced order modeling for nonlinear structural analysis using
  {G}aussian process regression.
\newblock {\em Computer methods in applied mechanics and engineering},
  341:807--826, 2018.

\bibitem{horaud2010rigid}
R.~Horaud, F.~Forbes, M.~Yguel, G.~Dewaele, and J.~Zhang.
\newblock Rigid and articulated point registration with expectation conditional
  maximization.
\newblock {\em IEEE Transactions on Pattern Analysis and Machine Intelligence},
  33(3):587--602, 2010.

\bibitem{iollo2014advection}
A.~Iollo and D.~Lombardi.
\newblock Advection modes by optimal mass transfer.
\newblock {\em Physical Review E}, 89(2):022923, 2014.

\bibitem{iollo2022mapping}
A.~Iollo and T.~Taddei.
\newblock Mapping of coherent structures in parameterized flows by learning
  optimal transportation with {G}aussian models.
\newblock {\em Journal of Computational Physics}, 471:111671, 2022.

\bibitem{jasak2009openfoam}
H.~Jasak.
\newblock Open{FOAM}: Open source {CFD} in research and industry.
\newblock {\em International Journal of Naval Architecture and Ocean
  Engineering}, 1(2):89--94, 2009.

\bibitem{jeong1995identification}
J.~Jeong and F.~Hussain.
\newblock On the identification of a vortex.
\newblock {\em Journal of fluid mechanics}, 285:69--94, 1995.

\bibitem{kadeethum2022enhancing}
T.~Kadeethum, D.~O’Malley, F.~Ballarin, I.~Ang, J.~N. Fuhg, N.~Bouklas, V.~L.
  Silva, P.~Salinas, C.~E. Heaney, and C.~C. Pain.
\newblock Enhancing high-fidelity nonlinear solver with reduced order model.
\newblock {\em Scientific Reports}, 12(1):20229, 2022.

\bibitem{kolavr2007vortex}
V.~Kol{\'a}{\v{r}}.
\newblock Vortex identification: New requirements and limitations.
\newblock {\em International journal of heat and fluid flow}, 28(4):638--652,
  2007.

\bibitem{krah2023front}
P.~Krah, S.~B{\"u}chholz, M.~H{\"a}ringer, and J.~Reiss.
\newblock Front transport reduction for complex moving fronts.
\newblock {\em Journal of scientific computing}, 96(28), 2023.

\bibitem{lee2020model}
K.~Lee and K.~T. Carlberg.
\newblock Model reduction of dynamical systems on nonlinear manifolds using
  deep convolutional autoencoders.
\newblock {\em Journal of Computational Physics}, 404:108973, 2020.

\bibitem{luke2012fast}
E.~Luke, E.~Collins, and E.~Blades.
\newblock A fast mesh deformation method using explicit interpolation.
\newblock {\em Journal of Computational Physics}, 231(2):586--601, 2012.

\bibitem{ma2015robust}
J.~Ma, W.~Qiu, J.~Zhao, Y.~Ma, A.~L. Yuille, and Z.~Tu.
\newblock Robust {$L_2$ $E$} estimation of transformation for non-rigid
  registration.
\newblock {\em IEEE Transactions on Signal Processing}, 63(5):1115--1129, 2015.

\bibitem{ma2018nonrigid}
J.~Ma, J.~Wu, J.~Zhao, J.~Jiang, H.~Zhou, and Q.~Z. Sheng.
\newblock Nonrigid point set registration with robust transformation learning
  under manifold regularization.
\newblock {\em IEEE transactions on neural networks and learning systems},
  30(12):3584--3597, 2018.

\bibitem{maiseli2017recent}
B.~Maiseli, Y.~Gu, and H.~Gao.
\newblock Recent developments and trends in point set registration methods.
\newblock {\em Journal of Visual Communication and Image Representation},
  46:95--106, 2017.

\bibitem{MATLAB:2022}
MATLAB.
\newblock {\em R2022a}.
\newblock The MathWorks Inc., Natick, Massachusetts, 2022.

\bibitem{mccann1997convexity}
R.~J. McCann.
\newblock A convexity principle for interacting gases.
\newblock {\em Advances in mathematics}, 128(1):153--179, 1997.

\bibitem{mckay2000comparison}
M.~D. McKay, R.~J. Beckman, and W.~J. Conover.
\newblock A comparison of three methods for selecting values of input variables
  in the analysis of output from a computer code.
\newblock {\em Technometrics}, 42(1):55--61, 2000.

\bibitem{menter1994two}
F.~R. Menter.
\newblock Two-equation eddy-viscosity turbulence models for engineering
  applications.
\newblock {\em AIAA journal}, 32(8):1598--1605, 1994.

\bibitem{mirhoseini2023model}
M.~A. Mirhoseini and M.~J. Zahr.
\newblock Model reduction of convection-dominated partial differential
  equations via optimization-based implicit feature tracking.
\newblock {\em Journal of Computational Physics}, 473:111739, 2023.

\bibitem{modesti2017low}
D.~Modesti and S.~Pirozzoli.
\newblock A low-dissipative solver for turbulent compressible flows on
  unstructured meshes, with {Open-FOAM} implementation.
\newblock {\em Computers \& Fluids}, 152:14--23, 2017.

\bibitem{mojgani2017arbitrary}
R.~Mojgani and M.~Balajewicz.
\newblock Arbitrary {L}agrangian {E}ulerian framework for efficient
  projection-based reduction of convection dominated nonlinear flows.
\newblock In {\em APS Division of Fluid Dynamics Meeting Abstracts}, pages
  M1--008, 2017.

\bibitem{mojgani2021low}
R.~Mojgani and M.~Balajewicz.
\newblock Low-rank registration based manifolds for convection-dominated
  {PDEs}.
\newblock In {\em Proceedings of the AAAI Conference on Artificial
  Intelligence}, volume~35, pages 399--407, 2021.

\bibitem{myronenko2010point}
A.~Myronenko and X.~Song.
\newblock Point set registration: Coherent point drift.
\newblock {\em IEEE transactions on pattern analysis and machine intelligence},
  32(12):2262--2275, 2010.

\bibitem{nicoud1999subgrid}
F.~Nicoud and F.~Ducros.
\newblock Subgrid-scale stress modelling based on the square of the velocity
  gradient tensor.
\newblock {\em Flow, turbulence and Combustion}, 62(3):183--200, 1999.

\bibitem{ohlberger2013nonlinear}
M.~Ohlberger and S.~Rave.
\newblock Nonlinear reduced basis approximation of parameterized evolution
  equations via the method of freezing.
\newblock {\em Comptes Rendus Mathematique}, 351(23-24):901--906, 2013.

\bibitem{peherstorfer2020model}
B.~Peherstorfer.
\newblock Model reduction for transport-dominated problems via online adaptive
  bases and adaptive sampling.
\newblock {\em SIAM Journal on Scientific Computing}, 42(5):A2803--A2836, 2020.

\bibitem{persson2006sub}
P.-O. Persson and J.~Peraire.
\newblock Sub-cell shock capturing for discontinuous {G}alerkin methods.
\newblock In {\em 44th AIAA aerospace sciences meeting and exhibit}, page 112,
  2006.

\bibitem{quarteroni2009numerical}
A.~Quarteroni.
\newblock {\em Numerical models for differential problems}, volume~2.
\newblock Springer, 2009.

\bibitem{reiss2018shifted}
J.~Reiss, P.~Schulze, J.~Sesterhenn, and V.~Mehrmann.
\newblock The shifted proper orthogonal decomposition: A mode decomposition for
  multiple transport phenomena.
\newblock {\em SIAM Journal on Scientific Computing}, 40(3):A1322--A1344, 2018.

\bibitem{rice2006mathematical}
J.~A. Rice.
\newblock {\em Mathematical statistics and data analysis}.
\newblock Cengage Learning, 2006.

\bibitem{shepard1968two}
D.~Shepard.
\newblock A two-dimensional interpolation function for irregularly-spaced data.
\newblock In {\em Proceedings of the 1968 23rd ACM national conference}, pages
  517--524, 1968.

\bibitem{sirovich1987turbulence}
L.~Sirovich.
\newblock Turbulence and the dynamics of coherent structures. {I. C}oherent
  structures.
\newblock {\em Quarterly of applied mathematics}, 45(3):561--571, 1987.

\bibitem{taddei2020registration}
T.~Taddei.
\newblock A registration method for model order reduction: data compression and
  geometry reduction.
\newblock {\em SIAM Journal on Scientific Computing}, 42(2):A997--A1027, 2020.

\bibitem{taddei2023compositional}
T.~Taddei.
\newblock Compositional maps for registration in complex geometries.
\newblock {\em arXiv preprint arXiv:2308.15307}, 2023.

\bibitem{taddei2015reduced}
T.~Taddei, S.~Perotto, and A.~Quarteroni.
\newblock Reduced basis techniques for nonlinear conservation laws.
\newblock {\em ESAIM: Mathematical Modelling and Numerical Analysis},
  49(3):787--814, 2015.

\bibitem{tezduyar1992new}
T.~E. Tezduyar, M.~Behr, S.~Mittal, and J.~Liou.
\newblock A new strategy for finite element computations involving moving
  boundaries and interfaces—the deforming-spatial-domain/space-time
  procedure: {II.} computation of free-surface flows, two-liquid flows, and
  flows with drifting cylinders.
\newblock {\em Computer methods in applied mechanics and engineering},
  94(3):353--371, 1992.

\bibitem{tonon2021linear}
P.~Tonon, R.~A.~K. Sanches, K.~Takizawa, and T.~E. Tezduyar.
\newblock A linear-elasticity-based mesh moving method with no cycle-to-cycle
  accumulated distortion.
\newblock {\em Computational Mechanics}, 67:413--434, 2021.

\bibitem{volkwein2011model}
S.~Volkwein.
\newblock Model reduction using proper orthogonal decomposition.
\newblock {\em Lecture Notes, Institute of Mathematics and Scientific
  Computing, University of Graz. see http://www. uni-graz.
  at/imawww/volkwein/POD. pdf}, 1025, 2011.

\bibitem{wendland2004scattered}
H.~Wendland.
\newblock {\em Scattered data approximation}, volume~17.
\newblock Cambridge university press, 2004.

\bibitem{yano2019lp}
M.~Yano and A.~T. Patera.
\newblock An {LP} empirical quadrature procedure for reduced basis treatment of
  parametrized nonlinear {PDEs}.
\newblock {\em Computer Methods in Applied Mechanics and Engineering},
  344:1104--1123, 2019.

\bibitem{zitova2003image}
B.~Zitova and J.~Flusser.
\newblock Image registration methods: a survey.
\newblock {\em Image and vision computing}, 21(11):977--1000, 2003.

\end{thebibliography}

\end{document}